\documentclass[UTF-8,reqno]{amsart}
\usepackage{enumerate}
\usepackage{mhequ}
\usepackage{amssymb,url,color, booktabs,nccmath}
\usepackage[left=3.2cm,right=3.2cm,top=4cm,bottom=4cm]{geometry}
\usepackage{mathrsfs}
\usepackage{enumitem,dsfont}

\usepackage{color}
\usepackage[colorlinks=true]{hyperref}
\hypersetup{
    linkcolor=blue,          
    citecolor=red,        
    filecolor=blue,      
    urlcolor=cyan}

\definecolor{darkergreen}{rgb}{0.0, 0.5, 0.0}

\newtheorem{theorem}{Theorem}[section]
\newtheorem{lemma}[theorem]{Lemma}
\newtheorem{proposition}[theorem]{Proposition}

\newtheorem{corollary}[theorem]{Corollary}

\newtheorem{definition}{Definition}
\theoremstyle{definition}
\newtheorem{remark}{Remark}

\newcommand{\mathd}{\mathrm{d}}

\allowdisplaybreaks

\usepackage{esint}
\usepackage{cite}
\usepackage{verbatim}
\usepackage{times}
\usepackage{color}
\usepackage{bm}
\usepackage{fancyhdr}
\numberwithin{equation}{section}

\begin{document}
\title{Martingale Suitable Weak Solutions of $3$-D Stochastic Navier-Stokes Equations with Vorticity Bounds}
\author{Weiquan Chen}
\address[Weiquan Chen]{Academy of Mathematics and Systems Science,
Chinese Academy of Sciences, Beijing 100190, China}
\address[Weiquan Chen]{University of Chinese Academy of Science, Beijing 100190, China}
\email{chenweiquan@amss.ac.cn}

\author{Zhao Dong}
\address[Zhao Dong]{Academy of Mathematics and Systems Science,
Chinese Academy of Sciences, Beijing 100190, China}
\address[Zhao Dong]{University of Chinese Academy of Science, Beijing 100190, China}
\email{dzhao@amt.ac.cn}
\thanks{Research  supported by National Key R\&D Program of China (No. 2020YFA0712700) and the NSFC (No. 11931004, 12090010, 12090014)}

\begin{abstract}
	In this paper, we construct martingale suitable weak solutions for $3$-dimensional incompressible stochastic Navier-Stokes equations with generally non-linear noise. In deterministic setting, as widely known, ``suitable weak solutions'' are Leray-Hopf weak solutions enjoying two different types of local energy inequalities (LEIs). In stochastic setting, we apply the idea of ``martingale solution", avoid transforming to random system, and show new stochastic versions of the two local energy inequalities. In particular, in additive and linear multiplicative noise case, OU-processes and the exponential formulas DO NOT play a role in our formulation of LEIs. This is different to \cite{FR02,Rom10} where the additive noise case is dealt. Also, we successfully apply the concept of ``a.e. super-martingale'' to describe this local energy behavior. To relate the well-known ``dissipative weak solutions" come up with in \cite{DR00}, we derive a local energy equality and extend the concept onto stochastic setting naturally. For further regularity of solutions, we are able to bound the $L^\infty\big([0,T];L^1(\Omega\times\mathbb T^3)\big)$ norm of the vorticity and $L^{\frac{4}{3+\delta}}\big(\Omega\times[0,T]\times\mathbb T^3\big)$ norm of the gradient of the vorticity, in case that the initial vorticity is a finite regular signed measure.  
\end{abstract}\
\maketitle

\section{Introduction}\label{sec. introd.}
\setcounter{section}{1}

\par 
The Navier-Stokes equations play a pivotal role in the field of fluid mechanics. To date, many fundamental and profound mathematical issues related to these equations remain unresolved. Nevertheless, persistent research efforts have led to series of substantial results in this area. It's well known that there exists a global weak solution $u\in C\big([0,T];L^2_{w}\big)\bigcap L^2\big([0,T];H^1\big)$ to incompressible Navier-Stokes equations satisfying an energy inequalty, which is named after J. Leray and E. Hopf in honor of their fundamental works \cite{Ler34, Hop51}. For the case $d=2$, it's also famously known that the system (\ref{sto. NS}) is well-posed in the class of Leray-Hopf solutions and these weak solutions are in fact strong solutions. In contrast, whether the uniqueness of Leray-Hopf solutions continues to hold for $d\geq3$ remains open to date and makes up the main content of the $Clay~Millennium~Prize$ problem for N-S equations. A remarkable attempt in the history to solve the uniqueness problem was the famous $weak$-$strong~uniqueness$ result which is gradually extended to the version known today by \cite{Pro59, Ser62, Lad69, FJR72, Kat84, FLRT00,LM01, ESS03}, see also \cite{CL22} for the current version and an alternative proof. In principle, the $weak$-$strong~uniqueness$ result transfers the problem of uniqueness of Leray-Hopf solutions to that of existence of global strong solutions, which suggests that to verify the smoothness (or integrability within L-P-S scales) of Leray-Hopf solutions is also a possible way to solve the uniqueness problem. In his pioneering works \cite{Sch77,Sch80,Sch87}, V. Scheffer constructs weak solutions with space-time singularity set of finite $5/3$-dimensional Hausdorff measure and spatial singularity set (at each time slice) of finite $1$-dimensional Hausdorff measure. This inspires a series of remarkable works now known as Partial regularity results. In this direction, one of the most important results is the so-called CKN Theorem named after L. Caffarelli, R. Kohn and L. Nirenberg. In \cite{CKN82} they construct a Leray-Hopf weak solution which satisfies a local energy inequality and show that such solutions possess space-time singularity set of zero $1$-dimensional Hausdorff measure. See \cite{Koh82,Lin98,LS99,Vas07,Kuk09a} for instructive summary and alternative proofs of CKN, and \cite{Kuk09b,KP12} for recent improvements on bounds for the box-counting dimension of space-time singularity set. 
\par Leray-Hopf weak solutions satisfying a local energy inequality are called suitable weak solutions in literatures. And such solutions are constructed in \cite{CKN82} via a regularization of the original Navier-Stokes equations, where the non-linear term $u\cdot\nabla u$ is replaced by $(\Psi_\varepsilon\ast u)\cdot\nabla u$. Here, differently from the Leray regularization used in \cite{Ler34}, $\Psi_\varepsilon$ is a well-chosen space-time mollifier. It's then natural to ask whether Leray regularization would also lead to existence of suitable weak solutions? And what about the Galerkin approximation used in \cite{Hop51}? An affirmative proof of the first question is given in \cite{BCI07} while the other one remains open to date. Even so, it's proven in \cite{Gue06,Gue07} that the Feado-Galerkin strategy still works if one carefully choose the finite-dimensional spaces for approximations. It's well worth mentioning that there are many other methods which also succeed in constructing suitable weak solutions. For example, see \cite{Bei85a,Bei85b,BCI07} for methods by using vanishing super-viscosity, \cite{Lio94} for a similar method to Leray regularization, \cite{BS16} for semi-discretization method and \cite{BS18} for an artificial compressibility method.
\par In the stochastic setting, due to the fact that the equation is not differentiable in time, the compactness which follows by uniform bounds on time derivatives of approximating solutions does not work as in deterministic PDEs. And people are forced to look for new compactness method to construction solutions. In classical literatures existence of martingale weak solutions is achieved by either estimates on the modulus of continuity of approximating solutions \cite{Vio76,Mef88}, or the factorization of stochastic integrals \cite{DaZ92,Gat93,GD94}. In \cite{FG95} the result is extended to the case of regular and coercive diffusion coefficients of noise by introducing a new compactness method based on fractional Sobolev space. This method also yields  the existence of martingale Leray-Hopf weak solutions for additive noise case, see \cite{FR08}. We should note here that linear multiplicative noise became popular due to its stabilizing effect: as shown in \cite{GHV14,RZZ14}, under such noise local in time strong solution can be extended globally with high probability. The existence of probabilistically strong and analytically weak solutions was open until very recently a stochastic version of the convex integration method is achieved in \cite{HZZ19,HZZ21a,HZZ21b} where the authors are able to construct non-unique probabilistically strong $C_tH^\gamma$-weak solutions for small $\gamma>0$. And non-unique probabilistically strong $L_t^pL^\infty$-weak solutions ($1\leq p<2$) for $d$-dimensional case ($d\geq2$) is obtained in \cite{CDZ23} via a new stochastic convex integration method. However, these solutions are not known to satisfy an energy inequality in the sense of Leray-Hopf. And existence of probabilistically strong Leray-Hopf weak solutions remains open to date. Still and all, it seems possible to achieve martingale Leray-Hopf weak solutions with more delicate energy balance and further space-time regularity. In \cite{Rom10}, M. Romito successfully constructed martingale suitable weak solution for the additive noise case and arrived at local energy inequality involving the OU-process which solves the corresponding stochastic Stokes system. And in \cite{FR02} they were able to show a partial regularity result to the stochastic N-S in spirit of CKN.
\par In this paper, we construct martingale suitable weak solutions to stochastic Navier-Stokes with generally non-linear noise. In particular, in additive and linear multiplicative noise case, the OU-process and the exponential formulas (which are widely used to transform the stochastic PDEs system to a random one) do not enter our formulations of local energy inequality. This is different to \cite{FR02,Rom10}. We derive two typical types of local energy inequalities where a resulting semi-martingale term shows up as the effect of the stochastic noise, see (\ref{LEI-I}) and (\ref{LEI. II}) below. Comparing with the deterministic side, a new and interesting angle is that energy behavior of stochastic Navier-Stokes can be described in terms of ``$a.e.$ super-martingale". The concept of $a.e.$ super-martingale was introduced for the first time by Flandoli and Romito \cite{FR08} to describe global (in space) energy behavior of Leray-Hopf martingale weak solutions for additive noise case. Our result successfully describe local energy behavior of martingale suitable weak solutions for stochastic noise under generally non-linear setting, see Theorem \ref{Thm. energy supermartingale} below. Due to non-linearity of noise and the goal of achieving a natural local energy inequality, we are not allowed to transform the system to random PDEs. The key idea is to apply Leray regularization, make use of It\^o's formula in appropriate ways and try to apply the idea of constructing martingale solutions in deriving local energy inequality. Passing to the limit requires further tightness properties than that in \cite{FG95,FR08}. And to derive the resulting semi-martingale term in our stochastic local energy inequalities, we have to investigate corresponding martingale convergence. To our knowledges, the results and treatments are new in existing literature.
\par To see a clearer picture of the local energy balance, we also apply a mollification procedure to derive a local energy equality which extends the concept of ``dissipative weak solutions" of J. Duchon and R. Robert \cite{DR00} onto stochastic Navier-Stokes equations. See Theorem \ref{Thm. LEE}, remark \ref{Remark-LEE} and Corollary \ref{Corollary-LEE} below.      
\par On the other hand, it is well-known that a sufficient condition for global regularity of deterministic N-S is the boundedness of $\int_0^T\left(\int_{\mathbb T^3}|\Delta u|^2{\rm d}x\right)^{2/3}{\rm d}t$, see \cite{CF89} for example. However, without further assumptions than the initial data $u_0\in L^2_\sigma$, Leray-Hopf weak solutions are not yet known to satisfy this condition but only admit certain bounds regarding derivatives of high orders, for example, the bound $\int_0^T\left(\int_{\mathbb T^3}|\Delta u|^2{\rm d}x\right)^{1/3}{\rm d}t$ (c.f. \cite{FGT81}). In \cite{Con90}, P. Constantin applies the method of change of dependent variables to show the a priori bound $\int_0^T\int_{\mathbb T^3}|\Delta u|^{\frac{4}{3+\delta}}{\rm d}x{\rm d}t$ for suitable weak solutions, assuming only that the initial vorticity is a finite regular signed measure. Another benefit of our method mentioned above, is that we are able to apply the method of change of dependent variables to the Leray regularised solution in the context of stochastic calculus. And paralleled to the deterministic side \cite{Con90}, we are able to bound the $L^\infty\big([0,T];L^1(\Omega\times\mathbb T^3)\big)$ norm of the vorticity of martingale suitable weak solutions, and the $L^{\frac{4}{3+\delta}}\big(\Omega\times[0,T]\times\mathbb T^3\big)$ norm of the gradient of the vorticity, in case that the initial vorticity is a finite regular signed measure. Such results are also new in the stochastic analysis literature.
\par We thank the referee for very useful suggestions especially about generallizing the conditions on noise. 

\subsection{Main results}
\par In this paper, we consider Cauchy problem of the $3$-dimensional 
incompressible stochastic Navier-Stokes equations, posed on the torus $\mathbb{T}^3\subset\mathbb R^3$ with periodic boundary conditions, with a generally non-linear noise:
\\
\begin{equation}\label{sto. NS}
	\left\{\begin{aligned}
		&{\rm d}u(t)=\Big(-{\rm div}\big(u(t)\otimes u(t)\big)+\nu \Delta u(t)-\nabla p(t)\Big){\rm d}t+\sum_{k=1}^\infty\varsigma_k(u(t)){\rm dB}_t^k,\\
		&{\rm div}\ u(t)=0,\\
		&u(0)=u_0\in L^2_\sigma,
	\end{aligned}\right.
\end{equation}
\\
where $\nu>0$ is constant, $\{\varsigma_k\}$ is a family of (deterministic) continuous operators from $L^2(\mathbb{T}^3;\mathbb R^3)$ to $L^2(\mathbb{T}^3;\mathbb R^3)$, $\left\{{\rm B}^k\right\}$ an independent family of one dimensional standard Brownian motions on some probability basis $\left(\Omega, \mathscr{F}, \{\mathscr{F}_t\}_{t\geq0}, \mathbb{P}\right)$, and $L^2_\sigma=\left\{v\in L^2: {\rm div}v=0\right\}$. 
\par For technical reason we need to impose some conditions on the operators $\{\varsigma_k\}$:
\begin{align}
    &Linear~Growth:\ \sum_{k=1}^\infty\big\|\varsigma_k(u)\big\|_{L^2}^2\lesssim1+\|u\|_{L^2}^2~,~~\forall u\in L^2, \label{con. linear growth}\\
	&Uniform~Decay~of~the~Tails:\ \lim_{N\rightarrow\infty}\sup_{u\in L^2_\sigma}\frac{\sum_{k>N}\big\|\varsigma_k(u)\big\|_{L^2}^2}{1+\|u\|_{L^2}^2}=0, \label{con. uniform decay of the tails}\\
	&Control~on~Vorticity:\ \sum_{k=1}^\infty\big\|\nabla\times\varsigma_k(u)\big\|_{L^2}^2\lesssim1+\|\nabla u\|_{L^2}^2~,~~\forall u\in H^1. \label{con. control on vorticity}
\end{align}
We leave two remarks on the conditions.
\begin{remark}
	Conditions (\ref{con. linear growth}-\ref{con. uniform decay of the tails}) are sufficient to show existence of martingale suitable weak solutions and the local energy inequalities. (\ref{con. control on vorticity}) is only required for the vorticity bounds (\ref{further bounds for vorticity}) below.
\end{remark} 
\begin{remark}
	We give some examples of noise satisfying conditions (\ref{con. linear growth}-\ref{con. control on vorticity}) :\\
	{$a)$} Additive Noise : $\varsigma_k(u)\equiv\varsigma_k\in H^1(\mathbb T^3;\mathbb R^3)$ with
	\begin{align}
		\sum_{k=1}^\infty\big\|\varsigma_k\big\|_{L^2}^2<\infty~for~(\ref{con. linear growth}-\ref{con. uniform decay of the tails})\quad and\quad\sum_{k=1}^\infty\big\|\nabla\times\varsigma_k\big\|_{L^2}^2<\infty~for~(\ref{con. control on vorticity});~\nonumber
	\end{align}
	{$b)$} Linear Multiplicative Noise : $\varsigma_k(u)=\varsigma_k u$ with $\varsigma_k\in C^1(\mathbb T^3;\mathbb R)$ satisfying
	\begin{align}
		\sum_{k=1}^\infty\|\varsigma_k\|_{C^0}^2<\infty~for~(\ref{con. linear growth}-\ref{con. uniform decay of the tails})\quad and\quad\sum_{k=1}^\infty\|\varsigma_k\|_{C^1}^2<\infty~for~(\ref{con. control on vorticity})~\nonumber
	\end{align}
	\quad(the verification of (\ref{con. control on vorticity}) would need to use Poincar\'e's inequality);\\
	\\
	{$c)$} Cosine coefficients : $\varsigma_k(u)=f_k\cos\left(k\sqrt{1+|u|^2}\right)$ with $f_k\in C^1(\mathbb T^3;\mathbb R^3)$ satisfying
	\begin{align}
		\sum_{k=1}^\infty\big\|f_k\big\|_{L^2}^2<\infty~for~(\ref{con. linear growth}-\ref{con. uniform decay of the tails})\quad and\quad\sum_{k=1}^\infty\big\|\nabla\times f_k\big\|_{L^2}^2+\sum_{k=1}^\infty k^2\big\|f_k\big\|_{C^0}^2<\infty~for~(\ref{con. control on vorticity});~\nonumber
	\end{align}
\end{remark}
To state our results, we first recall the definition of martingale weak solutions.

\begin{definition}\label{def. weak solu.}
We say that there exists a martingale weak solution $\Big\{\left(\Omega, \mathscr{F}, \{\mathscr{F}_t\}_{t\in[0,T]}, \mathbb{P}\right),u\Big\}$ to (\ref{sto. NS}) on $[0,T]$ with initial data $u_0\in L^2_\sigma(\mathbb{T}^3)$, if there exists a probability basis $\left(\Omega, \mathscr{F}, \{\mathscr{F}_t\}_{t\in[0,T]}, \mathbb{P}\right)$ and a progressively measurable process $u: [0,T]\times\Omega\longrightarrow L^2_\sigma(\mathbb{T}^3)$ with $\mathbb{P}$-$a.s.$ paths: 
\begin{align}
	u\in C\big([0,T]; H^{-\alpha}\big)\cap L^2\big([0,T]; H^{1}\big)\ \,for~some\ \,\alpha\in\mathbb R
\end{align}
such that :\\
(i) The paths of process $p(t):=(-\Delta)^{-1}{\rm divdiv}\big(u(t)\otimes u(t)\big)$, $t\in[0,T]$ which corresponds to the gradient pressure lay in $L^{5/3}\left([0,T]\times\mathbb T^3\right)$,~$\mathbb{P}$-$a.s.$ ; \\\\
(ii) For any test function $\varphi\in C^\infty(\mathbb{T}^3;\mathbb R^3)$, the process
\begin{align}\label{weak solu.} 
	\mathbf M_t(\varphi):=\big\langle u(t),\varphi\big\rangle - \big\langle u_0,\varphi\big\rangle\ - \int^{t}_{0}\Big\langle-{\rm div}\big(u(s)\otimes u(s)\big)+\nu \Delta u(s)-\nabla p(s) ~,~ \varphi\Big\rangle~{\rm d}s,
	\quad t\in[0,T]
\end{align}
is a real-valued continuous square-integrable $\{\mathscr{F}_t\}$-martingale with quadratic variation
\begin{align}
	\big\langle\mathbf M(\varphi)\big\rangle_t=\sum_{k=1}^\infty\int^{t}_{0}\big\langle\varsigma_k(u(s)),\varphi\big\rangle^2{\rm d}s~.
\end{align}
\end{definition}
\
\par Among martingale weak solutions, those obeying a local energy inequality are called {\bf martingale suitable weak solutions}. We define them as follows. 
\\
\begin{definition}[Martingale Suitable Weak Solutions]\label{def. martingale suitable weak solutions}
	A martingale weak solution $\Big\{\left(\Omega, \mathscr{F}, \{\mathscr{F}_t\}_{t\in[0,T]}, \mathbb{P}\right),u\Big\}$ is said to be suitable if, for every non-negative test function $\varphi\in C^\infty_c\left((0,T)\times\mathbb{T}^3\right)$, there exists a real-valued continuous square-integrable $\{\mathscr{F}_t\}$-martingale $\mathbf N(\varphi)=\mathbf N(u;\varphi)$ with quadratic variation
\begin{align}\label{quadratic variation of energy martingale}
	\big\langle\mathbf N(\varphi)\big\rangle_t=4\sum_{k=1}^\infty\int^{t}_{0}\left(\int_{\mathbb{T}^3}\varsigma_k(u(\tau))\cdot u(\tau)\varphi(\tau){\rm d}x\right)^2{\rm d}\tau
\end{align}
such that the local energy inequality
\begin{align}\label{LEI-I}
	\mathbb E\,\xi\,\mathcal E_T(u;\varphi)\leq\mathbb E\,\xi\left(\sum_{k=1}^{\infty}\int_0^T\int_{\mathbb{T}^3}\big|\varsigma_k(u(\tau))\big|^2\varphi(\tau){\rm d}x{\rm d}\tau+\mathbf N_T(\varphi)\right)
\end{align}
holds for all non-negative random variable $\xi\in L^\infty(\Omega, \mathscr{F}, \mathbb{P})$. 
Here,
\begin{align}\label{energy supermartingale}
	\mathcal E_t(u;\varphi):
	&=\int_{\mathbb{T}^3}\big|u(t)\big|^2\varphi(t){\rm d}x+2\nu\int_0^t\int_{\mathbb{T}^3}\big|\nabla u(\tau)\big|^2\varphi(\tau){\rm d}x{\rm d}\tau-\int_0^t\int_{\mathbb{T}^3}\big|u(\tau)\big|^2\Big(\partial_\tau\varphi(\tau)+\nu\Delta\varphi(\tau)\Big){\rm d}x{\rm d}\tau\nonumber\\
	&-\int_0^t\int_{\mathbb{T}^3}\Big(u(\tau)^2+2p(\tau)\Big)u(\tau)\cdot\nabla\varphi(\tau){\rm d}x{\rm d}\tau~,\quad t\in[0,T].
\end{align}
\end{definition}
\begin{remark}
	Note that in (\ref{LEI-I}) the local energy term $\displaystyle\mathbb E\,\xi\int_{\mathbb{T}^3}\big|u(T)\big|^2\varphi(T){\rm d}x$ vanishes, since ${\rm supp}_t\varphi\subset(0,T)$. And the inequality is in fact a bound for the local enstrophy term $\displaystyle 2\nu\mathbb E\,\xi\int_0^T\int_{\mathbb{T}^3}\big|\nabla u(\tau)\big|^2\varphi(\tau){\rm d}x{\rm d}\tau$. In deterministic case, (\ref{LEI-I}) is simply the known local energy inequality, i.e. $\mathcal E_T(\varphi)\leq0$.
\end{remark}

\par Now we are ready to state our results. 

\begin{theorem}\label{Thm. existence of martingale suitable solution}
	Let $\{\varsigma_k\}$ satisfy (\ref{con. linear growth}-\ref{con. uniform decay of the tails}). Then for any given $r\in[1,\infty),~T>0$ and initial data $u_0\in L^2_\sigma$, there exists a martingale suitable weak solution $\Big\{\left(\Omega, \mathscr{F}, \{\mathscr{F}_t\}_{t\in[0,T]}, \mathbb{P}\right),u\Big\}$ to (\ref{sto. NS}) on $[0,T]$ such that:\\\\
	i)~$u\in L^r\big(\Omega; L^\infty\big([0,T];L^2\big)\big)\cap L^2\big(\Omega; L^2\big([0,T];H^1\big)\big)$;\\\\
	ii)~$u$ has $\mathbb{P}$-$a.s.$ paths in $C\big([0,T];L^2_{w}\big)\cap L^{10/3}\big([0,T]\times\mathbb T^3\big)$.\\ 
	
	\noindent If additionally $\omega_0=\nabla\times u_0\in\mathcal M_r\big(\mathbb T^3,\mathscr B(\mathbb T^3)\big)$ and $\{\varsigma_k\}$ satisfy (\ref{con. control on vorticity}), the solution can be constructed such that $\omega=\nabla\times u$ satisfies for any pre-given $0<\delta\leq1/2$ :
	\begin{align}\label{further bounds for vorticity}
		\mathop{\rm ess~sup}_{t\in[0,T]}\mathbb E\|\omega(t)\|_{L^1}+\mathbb E\int_0^T\int_{\mathbb T^3}|\nabla\omega(t,x)|^{\frac{4}{3+\delta}}\mathd x\mathd t\leq C\left(\|u_0\|_{L^2},\|\omega_0\|_{var},\delta,T\right)
	\end{align}
	with constant $C>0$.
\end{theorem}

\par As in deterministic case, we can also deduce a local energy inequality in a different form where the localized energy $\displaystyle\mathbb E\,\xi\int_{\mathbb{T}^3}\big|u(t)\big|^2\varphi(t){\rm d}x$ is also bounded. And of course, this inequality is equivalent to (\ref{LEI-I}). 
 
\begin{theorem}\label{Thm. LEI. II} 
Let $\{\varsigma_k\}$ satisfy (\ref{con. linear growth}-\ref{con. uniform decay of the tails}). And let $\Big\{\left(\Omega, \mathscr{F}, \{\mathscr{F}_t\}_{t\in[0,T]}, \mathbb{P}\right),u\Big\}$ be the martingale suitable weak solution yielded by Theorem \ref{Thm. existence of martingale suitable solution}. Then for any non-negative test function $\varphi\in C^\infty_c\left((0,T)\times\mathbb{T}^3\right)$ and non-negative random variable $\xi\in L^\infty(\Omega, \mathscr{F}, \mathbb{P})$, we have for all $t\in[0,T]$ that
\begin{align}\label{LEI. II}
	\mathbb E\,\xi\,\mathcal E_t(u;\varphi)\leq\mathbb E\,\xi\left(\sum_{k=1}^{\infty}\int_0^t\int_{\mathbb{T}^3}\big|\varsigma_k(u(\tau))\big|^2\varphi(\tau){\rm d}x{\rm d}\tau+\mathbf N_t(\varphi)\right).
\end{align}
\end{theorem}

\par Typically in the literature on deterministic Navier-Stokes, (\ref{LEI. II}) is derived from (\ref{LEI-I}) by setting a well-chosen family of test functions $\varphi_\varepsilon$ with parameter $\varepsilon>0$ in (\ref{LEI-I}) and taking $\varepsilon\rightarrow0$. In the stochastic case here we need to deal with the convergence of the martingale term $\mathbf N_T(\varphi_\varepsilon)$. And this is done by investigating certain continuity of the map $\varphi\longmapsto\mathbf N(\varphi)$. 
\par Another interesting way to describe energy behavior of stochastic Navier-Stokes equations is by utilizing the concept of $a.e.$ super-martingales. The following result is a local version of the ``Leray-Hopf energy inequality" for stochastic Navier-Stokes given in \cite{FR08} (for additive noise). 

\begin{theorem}\label{Thm. energy supermartingale}
	Let $\{\varsigma_k\}$ satisfy (\ref{con. linear growth}-\ref{con. uniform decay of the tails}). And let $\Big\{\left(\Omega, \mathscr{F}, \{\mathscr{F}_t\}_{t\in[0,T]}, \mathbb{P}\right),u\Big\}$ be the martingale suitable weak solution yielded by Theorem \ref{Thm. existence of martingale suitable solution}. Then for any non-negative test function $\varphi\in C^\infty_c\left((0,T)\times\mathbb{T}^3\right)$, the process
\begin{align}\label{energy a.e. supermartingale}
	t\longmapsto\mathcal E_t(u;\varphi)
	-\sum_{k=1}^{\infty}\int_0^t\int_{\mathbb{T}^3}\big|\varsigma_k(u(\tau))\big|^2\varphi(\tau){\rm d}x{\rm d}\tau,\quad t\in[0,T]
\end{align}
is an $a.e.$ $\{\mathscr{F}_t\}$-super-martingale.
\end{theorem}

\begin{definition}[c.f. \cite{FR08}]\label{Remark a.e. sup-martingale}
We say that a stochastic process ${\rm X}:[0,T]\times\Omega\rightarrow\mathbb R$ is an $a.e.$ $\{\mathscr{F}_t\}$-super-martingale if :\\
(i) $\rm X$ is adapted to the filtration $\{\mathscr{F}_t\}$;\\
(ii) $\mathbb E\left|{\rm X}_t\right|<+\infty$ for all $t\in[0,T]$;\\
(iii) For all $t\in[0,T]$, there exists a set $\mathcal T_t\subset[0,t]$ of Lebesgue measure zero s.t. $\mathbb E\big(\left.{\rm X}_t-{\rm X}_s\right|\mathscr{F}_s\big)\leq0,\forall s\in[0,t]\backslash\mathcal T_t$.\\
An ``$a.e.$ $\{\mathscr{F}_t\}$-sub-martingale" is defined similarly by reversing the inequality in (iii).
\end{definition}

Finally, by a mollification procedure, we are able to derive a local energy equality in distribution sense, extending the concept of ``dissipative weak solution" in the sense of J. Duchon and R. Robert \cite{DR00} onto stochastic Navier-Stokes. In the following, we apply the convenient notation of spatial displacement of velocity field: $\delta_yu(t,x):=u(t,x+y)-u(t,x)$.

\begin{theorem}\label{Thm. LEE}
	Let $\{\varsigma_k\}$ satisfy (\ref{con. linear growth}-\ref{con. uniform decay of the tails}) and let $\Big\{\left(\Omega, \mathscr{F}, \{\mathscr{F}_t\}_{t\in[0,T]}, \mathbb{P}\right),u\Big\}$ be the martingale suitable weak solution yielded by Theorem \ref{Thm. existence of martingale suitable solution}. Then for any standard mollifier $\alpha_{\ell}$ (in space variables), the limit
	\begin{align}\label{dissipation distribution D(u)}
		\mathcal D_t(u;\varphi):=\lim_{\ell\rightarrow0}\int_0^t\int_{\mathbb T^3}\left(\underbrace{\frac{1}{4}\int_{\mathbb R^3}\nabla\alpha_{\ell}(y)\cdot\delta_yu(\tau,x)\big|\delta_yu(\tau,x)\big|^2{\rm d}y}_{=:\mathcal D^\ell(u)(\tau,x)}\right)\varphi(\tau,x){\rm d}x{\rm d}\tau~,\quad t\in[0,T]
	\end{align}
	exists in $L^1(\Omega)$ for every $\varphi\in C^\infty_c\big((0,T)\times\mathbb T^3\big)$ and $a.e.~t\in[0,T]$ (including $t=T$), and the following local energy equality holds true:
	\begin{align}\label{LEE 1}
		\mathcal E_t(u;\varphi)+2\mathcal D_t(u;\varphi)=\sum_{k=1}^\infty\int_0^t\int_{\mathbb T^3}\big|\varsigma_k(u(\tau))\big|^2\varphi(\tau){\rm d}x{\rm d}\tau+\mathcal N_t(u;\varphi)~,\quad a.e.~t\in[0,T].
	\end{align}
	Here $\mathcal N(u;\varphi)$ is a continuous square-integrable $\{\mathscr{F}_t\}$-martingale with the same quadratic variation as in (\ref{quadratic variation of energy martingale}), and the linear map $\varphi\mapsto\mathcal N(u;\varphi)$ is bounded from $L^5$ to $\mathcal M_{2,c}$ . Moreover, from Theorem \ref{Thm. energy supermartingale} we know that, for non-negative  test function $\varphi$, $\mathcal D(u;\varphi)$ is an $a.e.$ $\{\mathscr{F}_t\}$-sub-martingale.
\end{theorem}

Taking $t=T$ in (\ref{LEE 1}), it should be clear that for $\mathbb P$-$a.s.$ fixed $\omega\in\Omega$, both $\varphi\mapsto\mathcal N_T^\omega(u;\varphi)$ and $\varphi\mapsto\mathcal D_T^\omega(u;\varphi)$ define a distribution on $[0,T]\times\mathbb T^3$ which we denote by
\begin{align}
	\big\langle\mathcal N^\omega(u),\varphi\big\rangle=\mathcal N_T^\omega(u;\varphi),\quad \big\langle\mathcal D^\omega(u),\varphi\big\rangle=\mathcal D_T^\omega(u;\varphi)\nonumber
\end{align}
respectively. Similarly, each term in $\mathcal E_T^\omega(u;\varphi)$ defines a space-time distribution. And we are led to the following equality.
\begin{corollary}\label{Corollary-LEE}
	 The following energy equality holds :
	\begin{align}\label{LEE 2}
		(\partial_t-\nu\Delta)|u|^2+{\rm div}\Big[\big(|u|^2+2p\big)u\Big]+2\nu|\nabla u|^2+2\mathcal D(u)=\sum_{k=1}^\infty\big|\varsigma_k(u)\big|^2+\mathcal N(u)~,\quad\mathbb P-a.s.
	\end{align}
	where each term is understood as a random space-time distribution.
\end{corollary}

\begin{remark}\label{Remark-LEE}[Dissipative weak solution]
	 In \cite{DR00}, J. Duchon and R. Robert derive the local energy equality  in distribution sense:
	 \begin{align}
	 	(\partial_t-\nu\Delta)|u|^2+{\rm div}\Big[\big(|u|^2+2p\big)u\Big]+2\nu|\nabla u|^2+2\mathcal D(u)=0\nonumber
	 \end{align}
	 where $\displaystyle\mathcal D(u):=\lim_{\ell\rightarrow0}\mathcal D^\ell(u)$ is a non-negative distribution. Weak solutions satisfying the above energy equality with non-negative distribution term $\mathcal D(u)$ are classified as ``dissipative weak solutions" to deterministic incompressible Navier-Stokes equations. In this spirit, we may classify ``martingale weak solutions to (\ref{sto. NS}) satisfying energy equality (\ref{LEE 1}) with an $a.e.$ sub-martingale term $\mathcal D(u;\varphi)$" as ``{\bf martingale dissipative weak solutions}".
\end{remark}\


\subsection{Notations}
Throughout the paper, we employ the notation $a\lesssim b$ if there exists a constant $C>0$ such that $a\leq Cb$. And by $C(p,q)$ we mean that the constant $C$ depends on parameters $p$ and $q$. We would also use the notation $\mathbb N_0=\mathbb{N}\cup\{0\}$. 
\
\par {\bf Functional Spaces:} Let $\mathcal{S}'$ be the space of distributions on $\mathbb{T}^d$. For $s\in\mathbb R$ and $1\leq p\leq\infty$, we denote the Sobolev space $W^{s,p}=\big\{f\in\mathcal{S}':\big\|f\big\|_{W^{s,p}}:=\big\|(I-\Delta)^{s/2}f\big\|_{L^p}<\infty\big\}$. In case $p=2$, we use the more conventional notation $H^s=W^{s,2}$ and $\langle\cdot,\cdot\rangle_{H^s}$ for the corresponding inner product. Also, in case $s=0$, we have the usual Lebesgue spaces $L^p=W^{0,p}$. For simplicity, we always write $\langle\cdot,\cdot\rangle$ for the usual inner product for $L^2$ space. And we also use the notation $L^2_\sigma=\left\{v\in L^2: {\rm div}v=0~in~weak~sense\right\}$.
\par We use a convenient notation for Sobolev interpolation: let $0<p<s<q$ with $s=\theta p+(1-\theta)q$ for some $\theta\in(0,1)$, we denote by ``$H^s=\big[H^p,H^q\big]_\theta$" the interpolation $\|f\|_{H^s}\leq\|f\|_{H^p}^{\theta}\|f\|_{H^q}^{1-\theta}$.
\par Given a Banach space $\big(Y,~\|\cdot\|_Y\big)$ and $I\subset \mathbb{R}$, we write $L^p(I;Y)$ for the space of $L^p$-integrable functions from $I$ to $Y$, equipped with the usual $L^p$-norm. 
We also use $L^p_{loc}([0,\infty);Y)$ to denote the space of functions $f$ from $[0,\infty)$ to $Y$ satisfying $f\in L^p([0,T];Y)$ for all $T>0$. Similarly, we write $C(I;Y)$ for the space of continuous functions from $I$ to $Y$ equipped with the supremum norm in a bounded subset. 
We denote by $C(I;L^2_{w})$ the space of functions $u: I\rightarrow L^2$ such that $t\mapsto\langle u(t),v\rangle$ is a real continuous function for each $v\in L^2$.
\par For any Hilbert space $(H,\langle\cdot,\cdot\rangle_H)$, $r>1$ and $0<\alpha<1$, we recall from \cite{FG95} the fractional Sobolev space ${\mathcal W}^{\alpha,r}([0,T];H)=\big\{u\in L^r([0,T];H):\|u\|_{{\mathcal W}^{\alpha,r}}<\infty\big\}$ where the norm $\|\cdot\|_{{\mathcal W}^{\alpha,r}}$ is given by
\begin{align}
	\|u\|_{{\mathcal W}^{\alpha,r}}:=\int_0^T\|u(t)\|_H^r\mathd t+\int_0^T\int_0^T\frac{\|u(t)-u(s)\|_H^r}{|t-s|^{1+\alpha r}}\mathd t\mathd s~.\nonumber
\end{align}
\
\par {\bf Measure Space:} For a general Hausdorff space $X$, we denote by $\mathscr B(X)$ its Borel $\sigma$-algebra, and by $\mathcal M_r(X,\mathscr B(X))$ the space of all regular signed measures on $(X,\mathscr B(X))$ with finite total variation. Note that $\mathcal M_r(X,\mathscr B(X))$ is a Banach space endowed with the total variation norm
\begin{align}
	\|\mu\|_{var}:=\sup\left\{\displaystyle\sum_{n}|\mu(A_n)|:X=\sum_nA_n, \{A_n\}\subset\mathscr B(X)\right\}=\mu^+(X)+\mu^-(X)\nonumber
\end{align}
where $\displaystyle\sum_nA_n$ means pairwise disjoint union. It is well-known that the isomorphism
\begin{align}
	\mathcal M_r(X,\mathscr B(X))\cong \big(C_0(X)\big)^*\nonumber
\end{align}
is norm-preserving. Here
\begin{align}
	C_0(X)=\big\{continuous~function~f~on~X:~\forall\varepsilon>0,\exists~ compact~set~K_\varepsilon~s.t.~|f|\leq\varepsilon~on~K_\varepsilon^c\big\}\nonumber
\end{align}
and $\big(C_0(X)\big)^*$ denotes the space of all continuous linear functional on $C_0(X)$.
\\
\par {\bf Probability Spaces:} By $\left(\Omega, \mathscr{F}, \mathbb{P}\right)$ we denote a probability space where $\Omega$ is a sample space; $\mathscr{F}$ is a $\sigma$-field on $\Omega$ and $\mathbb{P}$ a probability measure on the measurable space $\left(\Omega, \mathscr{F}\right)$ such that $\mathscr{F}$ contains all $\mathbb{P}$-negligible sets. By $\left(\Omega, \mathscr{F}, \{\mathscr{F}_t\}_{t\in[0,T]}, \mathbb{P}\right)$ we denote a stochastic basis where a complete filtration $\{\mathscr{F}_t\}_{t\in[0,T]}$ is given on $\left(\Omega, \mathscr{F}, \mathbb{P}\right)$. We denote by $\mathcal M_{2,c}\big(\{\mathscr{F}_t\}_{t\in[0,T]}\big)$ the space of all continuous square-integrable martingales adapted to $\{\mathscr{F}_t\}_{t\in[0,T]}$. When the underlying filtration is obvious and no ambiguities should arise, we simply write $\mathcal M_{2,c}$. 
\par Given a Banach space $\big(Y,~\|\cdot\|_Y\big)$, for $1\leq p\leq\infty$, we denote by $L^p\left(\Omega, \mathscr{F}, \mathbb{P};Y\right)$ the space of $Y$-valued random variables $F$ on $\left(\Omega, \mathscr{F}, \mathbb{P}\right)$ such that $\mathbb{E}^\mathbb{P}\|F\|_Y^p<\infty$ when $p<\infty$ and $\mathop{\rm ess~sup}\limits_{\Omega}\|F\|_Y<\infty$ when $p=\infty$. When the underlying $\sigma$-field and probability measure are obvious and no ambiguities should arise, we simply write $L^p\left(\Omega;Y\right)$.\\

\subsection{Organization of the Paper}
In Section \ref{Sec. main}, we state the main idea of the paper and prove Theorem \ref{Thm. existence of martingale suitable solution} via two propositions, which are proved in Section \ref{Sec. proof of Leray regularised solution}. Section \ref{Sec. LEIs of various types} is divided into two parts where in subsection \ref{Subsec. Local energy inequality of type II} we prove Theorem \ref{Thm. LEI. II} and in subsection \ref{Subsec. Super-martingale statement of local energy inequality} Theorem \ref{Thm. energy supermartingale}. The Appendix \ref{Appen. Intuitive computation} contains an intuitive deviation of local energy equality for stochastic Navier-Stokes (\ref{sto. NS}), assuming good enough regularity.

\
\section{Main Idea and Proof of Theorem \ref{Thm. existence of martingale suitable solution}}\label{Sec. main}
\subsection{Main ideas and propositions}
The key properties that identify the suitable weak solutions are the local energy inequalities (\ref{LEI-I}) and (\ref{LEI. II}). For clarity, we give in Appendix an intuitive derivation for local energy inequalities. In short, one could assume a ``regular enough" solution to (\ref{sto. NS}) and apply It\^o's formula for each fixed $x$ and then integrate on $\mathbb T^3$ on both sides. Then a local energy equality follows by integration by part. Our main task here is to legitimize the calculations in the context of martingale weak solutions. The construction of solutions is based on the so-called Leray regularization (c.f. \cite{Ler34}) : For each fixed $\varepsilon>0$ and $N:=N(\varepsilon)=[1/\varepsilon]+1$, we consider the regularised stochastic Navier-Stokes equations
\begin{equation}\label{regularised sto. NS}
	\left\{\begin{aligned}
		&{\rm d}u^\varepsilon(t)=\Big(-{\rm div}\big((\psi_\varepsilon \ast u^\varepsilon)(t)\otimes u^\varepsilon(t)\big)+\nu \Delta u^\varepsilon(t)-\nabla p^\varepsilon(t)\Big){\rm d}t+\sum_{k=1}^N\psi_\varepsilon\ast\varsigma_k(u^\varepsilon(t)){\rm dB}^k_t,\\
		&{\rm div}\ u^\varepsilon(t)=0,\\
		&u^\varepsilon(0)=\psi_\varepsilon \ast u_0=:u_0^\varepsilon,
	\end{aligned}\right.
\end{equation}
where $\psi_\varepsilon=\varepsilon^{-3}\psi(\cdot/\varepsilon)$ is a standard mollifier, with function $\psi\in C^\infty_c\left(\mathbb R^3\right)$ such that $0\leq\psi\leq1$, $supp~\psi\subset B_1(0)$ and $\|\psi\|_{L^1}=1$. Note that in case where $\{\varsigma_k\}$ satisfy Lipschitz continuity one can work on the unique probabilistically strong solution to (\ref{regularised sto. NS}) directly. This is however not a necessary path to our final goal as we are dealing with tightness and compactness in what follows. In present situation where Lipschitz continuity does not necessarily hold for $\{\varsigma_k\}$ we may deal with martingale weak solutions to (\ref{regularised sto. NS}).

\begin{definition}\label{def. weak solu.^varepsilon}
We say that there exists a martingale weak solution $\Big\{\big(\Omega^\varepsilon, \mathscr{F}^\varepsilon,\{\mathscr{F}_t^\varepsilon\}_{t\geq0}, \mathbb{P}^\varepsilon\big),u^\varepsilon\Big\}$ to (\ref{sto. NS}) on $[0,T]$ with initial data $u_0\in L^2_\sigma(\mathbb{T}^3)$, if there exist a probability basis $\big(\Omega^\varepsilon, \mathscr{F}^\varepsilon,\{\mathscr{F}_t^\varepsilon\}_{t\geq0}, \mathbb{P}^\varepsilon\big)$ and a progressively measurable process $u^\varepsilon : [0,T]\times\Omega\longrightarrow L^2_\sigma(\mathbb{T}^3)$ with $\mathbb{P}$-$a.s.$ paths in $C\big([0,T]; H^{-\alpha}\big)\cap L^2\big([0,T]; H^{1}\big)$ for some $\alpha\in\mathbb R$, such that for any test function $\varphi\in C^\infty(\mathbb{T}^3)$, the process
\begin{align}
	t\longmapsto\big\langle u^\varepsilon(t),\varphi\big\rangle - \big\langle u_0,\varphi\big\rangle\ - \int^{t}_{0}\Big\langle-{\rm div}\big(u^\varepsilon(s)\otimes u^\varepsilon(s)\big)+\nu \Delta u^\varepsilon(s)-\nabla p^\varepsilon(s) ~,~ \varphi\Big\rangle~{\rm d}s,
	\quad t\in[0,T]\nonumber
\end{align}
( $p^\varepsilon:=(-\Delta)^{-1}{\rm divdiv}\big(u^\varepsilon\otimes u^\varepsilon\big)$ ) is a real-valued continuous square-integrable $\{\mathscr{F}_t\}$-martingale with quadratic variation $\displaystyle\sum_{k=1}^\infty\int^{\cdot}_{0}\big\langle\psi_\varepsilon\ast\varsigma_k(u^\varepsilon(s)),\varphi\big\rangle^2{\rm d}s$.
\end{definition}

\par One of the advantages of working on Leray regularization is that the system (\ref{regularised sto. NS}) has a solution which is smooth in $x\in\mathbb T^3$ for all $t>0$:

\begin{proposition}\label{Prop. Leray regularization}
	Let $\{\varsigma_k\}$ satisfy (\ref{con. linear growth}-\ref{con. uniform decay of the tails}). Then for any $\varepsilon>0$ and given $u_0\in L^2_\sigma$, there exists a martingale weak solution $\Big\{\big(\Omega^\varepsilon, \mathscr{F}^\varepsilon,\{\mathscr{F}_t^\varepsilon\}_{t\geq0}, \mathbb{P}^\varepsilon\big),u^\varepsilon\Big\}$ on $[0,T]$ to system (\ref{regularised sto. NS}) such that
	
	\noindent{\bf (I)} $u^\varepsilon\in L^2\big(\Omega^\varepsilon;C([0,T];L^2)\big)\cap L^2\big(\Omega^\varepsilon;L^2([0,T];H^1)\big)$;\\
	
	\noindent{\bf (II)} the ``noise term"
	\begin{align}\label{noise term N^varepsilon}
		{\mathfrak N}^\varepsilon:=u^\varepsilon-u^\varepsilon_0-\int_0^\cdot\big[-{\rm div}\big((\psi_\varepsilon \ast u^\varepsilon)(s)\otimes u^\varepsilon(s)\big)+\nu \Delta u^\varepsilon(s)-\nabla p^\varepsilon(s)\big]{\rm d}s
	\end{align}
	satisfies a fraction Sobolev estimate for any $r\geq2$ and $0<\alpha<1/2$ :
	\begin{align}\label{fraction Sobolev estimate for noise term}
		\mathbb E^\varepsilon\left\|{\mathfrak N}^\varepsilon\right\|_{{\mathcal W}^{\alpha,r}([0,T];L^2)}^r
	&\leq C(\alpha,p)\mathbb E^\varepsilon\int_0^T\left(\sum_{k=1}^N\left\|\psi_\varepsilon \ast \varsigma_k(u^\varepsilon(s))\right\|_{L^2}^2\right)^{r/2}{\rm d}s~;
	\end{align}
	
	\noindent{\bf (III)} the following bounds also hold for any $m\in\mathbb N$, $r\geq2$, $\gamma>3/2$ and $T>0$ :
	\begin{align}
	    &\mathbb E^\varepsilon\big\|u^\varepsilon\big\|_{C([0,T];L^2)}^r\leq C_1\big(\|u_0\|_{L^2},T,r\big)~,\ \ \ \ \ \ \mathbb E^\varepsilon\big\|u^\varepsilon\big\|_{C([0,T];H^m)}^r\leq C_2\big(\varepsilon,\|u_0^\varepsilon\|_{H^m},T,m,r\big)~,\label{standard uniform bounds 1}\\
		&\mathbb E^\varepsilon\big\|{\rm div}\big((\psi_\varepsilon\ast u^\varepsilon)\otimes u^\varepsilon\big)\big\|_{L^2([0,T];H^{-\gamma})}\nonumber\\		
		&+\mathbb E^\varepsilon\big\|u^\varepsilon\big\|_{L^2([0,T];H^1)}^2+\mathbb E^\varepsilon\big\|\Delta u^\varepsilon\big\|_{L^2([0,T];H^{-1})}^2
		\leq C_3(\|u_0\|_{L^2},T,\gamma)~,\label{standard uniform bounds 2}
	\end{align}
	with constants $C_1,C_2,C_3>0$.
	And, in particular, the path-wise regularity holds:
	\begin{align}\label{smoothness of regularised solution}
		\mathbb P^\varepsilon\Big(u^\varepsilon,p^\varepsilon\in C\big((0,T];H^m\big)~for~all~m\in\mathbb N\Big)=1~.
	\end{align}
\end{proposition}
With this in hand, we are able to justify application of It\^o's formula to the process $t\mapsto\frac{1}{2}\big|u^\varepsilon(t,x)\big|^2\varphi(t,x)$ for arbitrarily fixed $x\in\mathbb T^3$ (with test function $\varphi$). Then integrating over $\mathbb T^3$ and applying integration by part, this would lead to a path-wise local energy equality for $u^\varepsilon$. 

\begin{proposition}\label{Prop. LEE for Leray regularization}
Let $\Big\{\big(\Omega^\varepsilon, \mathscr{F}^\varepsilon,\{\mathscr{F}_t^\varepsilon\}_{t\geq0}, \mathbb{P}^\varepsilon\big),u^\varepsilon\Big\}$ be a martingale weak solution to system (\ref{regularised sto. NS}), as stated in Proposition \ref{Prop. Leray regularization}. Then for any test function $\varphi\in C^\infty(\mathbb T^3)$, the following equality holds $\mathbb{P}^\varepsilon$-a.s.:
	\begin{align}\label{LEE-a for u^varepsilon}
	&\ \ \ \ \ \ \int_{\mathbb{T}^3}\big|u^\varepsilon(t)\big|^2\varphi(t){\rm d}x+2\nu\int_0^t\int_{\mathbb{T}^3}\big|\nabla u^\varepsilon(s)\big|^2\varphi(s){\rm d}x{\rm d}s\nonumber\\
	&=\int_0^t\int_{\mathbb{T}^3}\big|u^\varepsilon(s)\big|^2\Big(\partial_\tau\varphi(s)+\nu\Delta\varphi(s)\Big){\rm d}x{\rm d}s+\sum_{k=1}^N\int_0^t\int_{\mathbb{T}^3}|\psi_\varepsilon\ast\varsigma_k(u^\varepsilon(s))|^2\varphi(s){\rm d}x{\rm d}s\nonumber\\
	&\ \ \ \ \ \int_0^t\int_{\mathbb{T}^3}\left(\big|u^\varepsilon(s)\big|^2\big(\psi_\varepsilon\ast u^\varepsilon(s)\big)+2 p^\varepsilon(s) u^\varepsilon(s)\right)\cdot\nabla\varphi(s){\rm d}x{\rm d}s+\mathbf N^\varepsilon_t(\varphi)
	~,\ \ \ \forall t\in[0,T]
\end{align}
where $\mathbf N^\varepsilon(\varphi)$ is a continuous square-integrable $\mathscr{F}_t^\varepsilon$-martingale with quadratic variation
\begin{align}\label{quadratic variation of N^varepsilon}
	\big\langle \mathbf N^\varepsilon(\varphi)\big\rangle_t=4\sum_{k=1}^N\int_0^t\left(\int_{\mathbb{T}^3}\psi_\varepsilon\ast\varsigma_k(\tilde u^\varepsilon(s))\cdot\tilde u^\varepsilon(s)\varphi(s){\rm d}x\right)^2{\rm d}s.
\end{align}
\end{proposition}
To show (\ref{LEI-I}), the remaining task then is to investigate suitable tightness and compactness which allows us to ``send $\varepsilon\rightarrow0$". Note that how $N$ depends on $\varepsilon$ is not important here as they can also be sent to infinity independently. We introduce the relation ``$N(\varepsilon)=[1/\varepsilon]+1$" here only to simplify notations.
\par Also, we are able to use the method of change of dependent variables to get further uniform bounds for the vorticity $\omega^\varepsilon=\nabla\times u^\varepsilon$ and gradient of the vorticity. Differently from \cite{Con90}, the method should be adapted for stochastic calculus in our settings. This is achieved by applying It\^o's formula to deduce a stochastic transport equality for the variable $q\big(\omega^\varepsilon\big)$ with specially chosen smooth function $q$, and carefully bounding each term by the uniform bounds (\ref{standard uniform bounds 1}). The result reads as follows.

\begin{proposition}\label{Prop. vorticity bounds for Leray regularised solution}
	Consider the same setting as in Proposition 2.1, and additionally assume that $\{\varsigma_k\}$ satisfy (\ref{con. control on vorticity}). Let also  $u_0\in L^2_\sigma$ and $\omega_0=\nabla\times u_0\in\mathcal M_r\big(\mathbb T^3,\mathscr B(\mathbb T^3)\big)$, then a martingale weak solution to (\ref{regularised sto. NS}) can be constructed such that its vorticity $\omega^\varepsilon=\nabla\times u^\varepsilon$  satisfies
	\begin{align}\label{vorticity bounds for Leray regularised solution}
		\mathbb E^\varepsilon\sup_{0\leq t\leq T}\|\omega^\varepsilon(t)\|_{L^1}+\mathbb E^\varepsilon\int_0^T\int_{\mathbb T^3}|\nabla\omega^\varepsilon(t,x)|^{\frac{4}{3+\delta}}\mathd x\mathd t\leq C\left(\|u_0\|_{L^2},\|\omega_0\|_{var},\delta,T\right),\quad\forall0<\delta\leq1/2
	\end{align}
	with constant $C>0$.
\end{proposition}

\par We give the proofs of Proposition \ref{Prop. Leray regularization}, \ref{Prop. LEE for Leray regularization} and \ref{Prop. vorticity bounds for Leray regularised solution} in Section \ref{Sec. proof of Leray regularised solution}. \\


\subsection{Existence of martingale suitable weak solutions} We apply a compactness and tightness argument based on the framework of fractional Sobolev space proposed in \cite{FG95}. The following compactness lemmas will be used to deduce tightness of the laws $\mathcal{L}^{u^\varepsilon}$ ($\varepsilon>0$).

\begin{lemma}[\cite{FG95}, Theorem 2.1]\label{Lem. fractional compact embeddings 1}
	Let $B_0\subset B\subset B_1$ be Banach spaces with $B_0$ and $B_1$ reflexive and suppose the embedding $B_0\subset B$ is compact. Let $1<r<\infty$ and $0<\alpha<1$ be given. Then the embedding
	\begin{align}
		L^r\big([0,T];B_0\big)\cap{\mathcal W}^{\alpha,r}\big([0,T];B_1\big)\subset L^r\big([0,T];B\big)\nonumber
	\end{align}
	is also compact.
\end{lemma}

\begin{lemma}[\cite{FG95}, Theorem 2.2]\label{Lem. fractional compact embeddings 2}
	Given finitely many Banach spaces $X,B_1,B_2,...,B_n$ ($n\in\mathbb N$) s.t.
	\begin{align}
		B_i\subset X~compactly,~\forall i=1,...,n\nonumber
	\end{align}
	and any indeces $\alpha_1,...,\alpha_n\in(0,1)$ and $p_1,...,p_n\in(1,\infty)$ s.t.
	\begin{align}
		\alpha_i p_i>1,~\forall i=1,...,n\nonumber
	\end{align}
	the embedding
	\begin{align}
		{\mathcal W}^{\alpha_1,p_1}\big([0,T];B_1\big)+...+{\mathcal W}^{\alpha_n,p_n}\big([0,T];B_1\big)\subset C\big([0,T];X\big)\nonumber
	\end{align}
	is compact.
\end{lemma}

The following tightness lemma will be used to deal with the martingale term in local energy balance. It can be found in the text book \cite{KS91} as an exercise problem (with answer).
\begin{lemma}[\cite{KS91}, Problem 2.4.11]\label{Lem. tightness on C-space}
	Let $\left\{X^{(n)}\right\}_{n\in\mathbb N}$ be a sequence of continuous stochastic processes $X^{(n)}=\left\{X^{(n)}_t;0\leq t<\infty\right\}$ on probability space $\left(\Omega, \mathscr{F},\mathbb{P}\right)$, satisfying the following conditions:
	\begin{align}
		&\sup_{n\in\mathbb N}\mathbb E\left|X^{(n)}_0\right|^\gamma<\infty~;\label{condition 1, tightness on C-space}\\
		&\sup_{n\in\mathbb N}\mathbb E\left|X^{(n)}_t-X^{(n)}_s\right|^\alpha\leq C_T|t-s|^{1+\beta}~,\quad\forall~T>0~and~\forall~s,t\in[0,T],\label{condition 2, tightness on C-space}
	\end{align}
	for some positive constants $\alpha,\beta,\gamma$ (universal) and $C_T$ (depending on $T$), then the probability measures $\mathbb P\left(X^{(n)}\right)^{-1}$ ($n\in\mathbb N$) induced by these processes on $\Big(C[0,\infty),\mathscr{B}\big(C[0,\infty)\big)\Big)$ form a tight sequence.
\end{lemma}
\begin{remark}
	It is not hard to see (following the proof in \cite{KS91}) that Lemma \ref{Lem. tightness on C-space} also applies to the case where each process $X^{(n)}$ is given on a perhaps distinct probability space $\left(\Omega^{(n)}, \mathscr{F}^{(n)},\mathbb{P}^{(n)}\right)$.
\end{remark}

\begin{proof}[\bf{Proof of Theorem \ref{Thm. existence of martingale suitable solution}}]
In what follow, we set $N=[1/\varepsilon]+1$. By ``$\varepsilon\rightarrow0$", or sometimes ``$\varepsilon>0$", we always mean that a sequence $\{\varepsilon_n\}$ of positive numbers is taken such that $\varepsilon_n\rightarrow0$ as $n\rightarrow\infty$. And when it's necessary,  some subsequence $\{\varepsilon_{n_k}\}$ will be subtracted from the original one such that $\varepsilon_{n_k}\rightarrow0$ as $k\rightarrow\infty$.
\par Let now $u_0\in L^2_\sigma$ be arbitrarily given. Then for each $\varepsilon>0$, Proposition \ref{Prop. Leray regularization} yields a martingale weak solution $\big\{\big(\Omega^\varepsilon, \mathscr{F}^\varepsilon,\{\mathscr{F}^\varepsilon_t\}_{t\geq0}, \mathbb{P}^\varepsilon\big),u^\varepsilon\big\}$ to system (\ref{regularised sto. NS}) with initial value $\psi_\varepsilon\ast u_0=:u_0^\varepsilon$, satisfying all properties stated in the proposition. And, of course, the pressure term is given by $p^\varepsilon=\Delta^{-1}{\rm divdiv}\big((\psi_\varepsilon\ast u^\varepsilon)\otimes u^\varepsilon\big)$.\\

\noindent$\bf\bullet~I.~Tightness~of~laws~\mathcal{L}^{(u^\varepsilon,\mathbf N^\varepsilon(\varphi))}\quad (\varepsilon>0)$
\par We first show certain tightness of the laws $\left\{\mathcal{L}^{u^\varepsilon}:=\mathbb{P}^\varepsilon(u^\varepsilon)^{-1}\right\}_{\varepsilon>0}$ . For this we write (\ref{noise term N^varepsilon}) as
\begin{align}\label{decomposition of u^varepsilon}
	u^\varepsilon(t)=u_0^\varepsilon+J_1^\varepsilon(t)+J_2^\varepsilon(t)+{\mathfrak N}^\varepsilon(t)
\end{align}
with
\begin{align}
	J_1^\varepsilon(t)&=-\int_0^t{\rm div}\big[(\psi_\varepsilon \ast u^\varepsilon)(s)\otimes u^\varepsilon(s)\big]{\rm d}s,\nonumber\\
	J_2^\varepsilon(t)&=\int_0^t\Big(\nu \Delta u^\varepsilon(s)-\nabla p^\varepsilon(s)\Big){\rm d}s.\nonumber
\end{align}
Then by (\ref{standard uniform bounds 2}), we see that
\begin{align}\label{uniform bounds for J1 and J2}
	\mathbb E^\varepsilon\left\|J_1^\varepsilon\right\|_{{\mathcal W}^{1,2}([0,T];H^{-\gamma})}+\mathbb E^\varepsilon\left\|J_2^\varepsilon\right\|_{{\mathcal W}^{1,2}([0,T];H^{-1})}\leq A(\nu,T,\gamma)
\end{align}  	
for any $\gamma>3/2$, $T>0$ and some constant $A>0$. 
And by (\ref{fraction Sobolev estimate for noise term}), (\ref{con. linear growth}) and (\ref{standard uniform bounds 1}), we have the estimate 
\begin{align}\label{uniform bounds for J3}
	\mathbb E^\varepsilon\left\|{\mathfrak N}^\varepsilon\right\|_{{\mathcal W}^{\alpha,r}([0,T];L^2)}^r
	&\leq C(\alpha,r)\mathbb E^\varepsilon\int_0^T\left(\sum_{k=1}^N\left\|\psi_\varepsilon \ast \varsigma_k(u^\varepsilon(s))\right\|_{L^2}^2\right)^{r/2}{\rm d}s\nonumber\\
	&\lesssim C(\alpha,r)\mathbb E^\varepsilon\int_0^T\big(1+\|u^\varepsilon(s)\|_{L^2}^2\big)^{r/2}{\rm d}s\nonumber\\	
	&\leq C(\alpha,r)B\big(\|u_0\|_{L^2},\nu,T,r\big)
\end{align}
for any $r\geq2$, $0<\alpha<1/2$ and some constant $B>0$. Hence, (\ref{decomposition of u^varepsilon}), (\ref{uniform bounds for J1 and J2}) and (\ref{uniform bounds for J3}) show that
\begin{align}
	\mathbb E^\varepsilon\left\|u^\varepsilon\right\|_X+\mathbb E^\varepsilon\left\|u^\varepsilon\right\|_Y\nonumber
\end{align}
is bounded uniformly with respect to $\varepsilon>0$, where
\begin{align}
	X&=L^2\big([0,T];H^1\big)\cap {\mathcal W}^{\alpha,2}\big([0,T];H^{-\gamma}\big),\nonumber\\
	Y&={\mathcal W}^{1,2}\big([0,T];H^{-\gamma}\big)+{\mathcal W}^{1,2}\big([0,T];H^{-1}\big)+{\mathcal W}^{\alpha,r}\big([0,T];L^2\big).\nonumber
\end{align}
By Lemma \ref{Lem. fractional compact embeddings 1} and Lemma \ref{Lem. fractional compact embeddings 2}, the embeddings $X\subset L^2\big([0,T];L^2\big)$ and $Y\subset C\big([0,T];H^{-s}\big)$ ($s>3/2$) are both compact. Hence, the family $\left\{\mathcal{L}^{u^\varepsilon}\right\}_{\varepsilon>0}$ of laws is tight on $C\big([0,T];H^{-s}\big)\cap L^2\big([0,T];L^2\big)$. Further more, since any bounded set in $L^\infty\big([0,T];L^2\big)\cap X$ is compact in $L^2\big([0,T];L^2\big)$, by the interpolation
\begin{align}
	\|f\|_{L^3([0,T];L^3)}
	&\leq\|f\|_{L^2([0,T];L^2)}^{1/6}\|f\|_{L^{10/3}([0,T];L^{10/3})}^{5/6}\nonumber\\
	&\leq\|f\|_{L^2([0,T];L^2)}^{1/6}\Big(\|f\|_{L^\infty([0,T];L^2)}^{2/5}\|f\|_{L^2([0,T];L^6)}^{3/5}\Big)^{5/6}\nonumber\\
	&\lesssim\|f\|_{L^2([0,T];L^2)}^{1/6}\|f\|_{L^\infty([0,T];L^2)}^{1/3}\|f\|_{L^2([0,T];H^1)}^{1/2},\nonumber
\end{align}
we see that the embedding $L^\infty\big([0,T];L^2\big)\cap X\subset L^3([0,T];L^3)$ is also compact. Thus, $\left\{\mathcal{L}^{u^\varepsilon}\right\}_{\varepsilon>0}$ is also tight on $L^3([0,T];L^3)$ and hence on $Z:=C\big([0,T];H^{-s}\big)\cap L^3\big([0,T];L^3\big)$. Here, $Z$ is equipped with the norm $\|\cdot\|_Z=\|\cdot\|_{C([0,T];H^{-s})}+\|\cdot\|_{L^3([0,T];L^3)}$.

\par Next, we show tightness of the laws $\left\{\mathcal{L}^{\mathbf N^\varepsilon(\varphi)}:=\mathbb{P}^\varepsilon\left(\mathbf N^\varepsilon(\varphi)\right)^{-1}\right\}_{\varepsilon>0}$ , where $\mathbf N^\varepsilon(\varphi)$ is the martingale term given in (\ref{LEE-a for u^varepsilon}). This is straightforward by Lemma \ref{Lem. tightness on C-space}, since from (\ref{quadratic variation of N^varepsilon}), one deduce by BDG-inequality H\"older's inequality, standard mollification estimate and (\ref{con. linear growth}) that
\begin{align}
	\mathbb{E}^\varepsilon\big|\mathbf N^\varepsilon_t(\varphi)-\mathbf N^\varepsilon_s(\varphi)\big|^m
	\leq & G_m2^m\mathbb{E}^\varepsilon\left[\sum_{k=1}^N\int_s^t\left(\int_{\mathbb{T}^3}\psi_\varepsilon\ast\varsigma_k( u^\varepsilon(\tau))\cdot u^\varepsilon(\tau)\varphi(\tau){\rm d}x\right)^2{\rm d}\tau\right]^{m/2}\nonumber\\
	\leq & G_m2^m\|\varphi\|_{L^\infty}^m\mathbb{E}^\varepsilon\left(\sum_{k=1}^N\int_s^t\big\|\psi_\varepsilon\ast\varsigma_k( u^\varepsilon(\tau))\big\|_{L^2}^2\big\|u^\varepsilon(\tau)\big\|_{L^2}^2{\rm d}\tau\right)^{m/2}\nonumber\\
	\leq & G_m2^m\|\varphi\|_{L^\infty}^m\mathbb{E}^\varepsilon\left[\int_s^t\left(1+\big\|u^\varepsilon(\tau)\big\|_{L^2}^2\right)^2{\rm d}\tau\right]^{m/2}\nonumber\\
	\leq & G_m2^m\|\varphi\|_{L^\infty}^m\mathbb{E}^\varepsilon\left(1+\|u^\varepsilon\|_{C([0,T];L^2)}^2\right)^m\cdot|t-s|^{m/2},~\forall~0\leq s\leq t\leq T\nonumber
\end{align}
where $G_m$ is nothing but the constant from BDG-inequality. Taking, say $m=4$, then by (\ref{standard uniform bounds 1}) one have
\begin{align}
	\sup_{\varepsilon>0}\mathbb{E}^\varepsilon\big|\mathbf N^\varepsilon_t(\varphi)-\mathbf N^\varepsilon_s(\varphi)\big|^4\lesssim\|\varphi\|_{L^\infty}^4\tilde C_1\cdot|t-s|^2,~\forall~0\leq s\leq t\leq T\nonumber
\end{align}
with constant $\tilde C_1=\tilde C_1\big(\|u_0\|_{L^2},T\big)>0$ that depends only on $\|u_0\|_{L^2},T$. Then by Lemma \ref{Lem. tightness on C-space} (clearly $\mathbf N^\varepsilon_0(\varphi)=0$), it is straightforward that the laws $\mathcal{L}^{\mathbf N^\varepsilon(\varphi)}$ ($\varepsilon>0$) forms a tight family of measures on $C[0,T]$.
\par Finally, we introduce addition and scalar multiplication on the Cartesian product $\mathscr{V}=Z\times C[0,T]$ :
\begin{align}
	(u,\omega)+(u',\omega')&=(u+u',\omega+\omega')\nonumber\\
	\lambda(u,\omega)&=(\lambda u,\lambda\omega)\nonumber
\end{align}
where $(u,\omega),(u',\omega')\in\mathscr{V}$ and $\lambda\in\mathbb R$. This makes $\mathscr{V}$ a vector space. Also, we equip $\mathscr{V}$ with the norm
\begin{align}
	\big\|(u,\omega)\big\|_{\mathscr{V}}:=\|u\|_Z+\sup_{t\in[0,T]}|\omega(t)|~,\quad (u,\omega)\in\mathscr{V}.\nonumber
\end{align}
There is no doubt that $\left(\mathscr{V},\|\cdot\|_{\mathscr{V}}\right)$ is a separable Banach space as are $Z$ and $C[0,T]$. Now, we view $\big(u^\varepsilon,\mathbf N^\varepsilon(\varphi)\big)$ as random variables taking values in the space $\mathscr{V}$. By what we prove above, for every $0<\delta<1$, there exist compact subsets $U_{\delta}$ of $Z$ and $K_{\delta}$ of $C[0,T]$ (both independent of $\varepsilon$) such that
\begin{align}
	&\mathcal{L}^{u^\varepsilon}\big(U_{\delta}^c\big)=\mathbb P^\varepsilon\big(u^\varepsilon\in U_{\delta}^c\big)\leq\delta~;\nonumber\\
	&\mathcal{L}^{\mathbf N^\varepsilon(\varphi)}\big(K_{\delta}^c\big)=\mathbb P^\varepsilon\big(\mathbf N^\varepsilon(\varphi)\in K_{\delta}^c\big)\leq\delta~,\quad \forall\varepsilon>0.\nonumber
\end{align}
Hence, for every $0<\delta<1$, it is clear that $U_{\delta/2}\times K_{\delta/2}$ is compact subset of $\mathscr{V}=Z\times C[0,T]$ and we have
\begin{align}
	\mathcal{L}^{(u^\varepsilon,\mathbf N^\varepsilon(\varphi))}\Big(\big(U_{\delta/2}\times K_{\delta/2}\big)^c\Big)
	&=\mathbb P^\varepsilon\Big(\big(u^\varepsilon,\mathbf N^\varepsilon(\varphi)\big)\in\big(U_{\delta/2}\times K_{\delta/2}\big)^c\Big)\nonumber\\
	&=\mathbb P^\varepsilon\Big(\big\{ u^\varepsilon\in U_{\delta/2}^c\big\}\cup\big\{\mathbf N^\varepsilon(\varphi)\in K_{\delta/2}^c\big\}\Big)\nonumber\\
	&\leq\mathbb P^\varepsilon\Big(u^\varepsilon\in U_{\delta/2}^c\Big)+\mathbb P^\varepsilon\Big(\mathbf N^\varepsilon(\varphi)\in K_{\delta/2}^c\Big)\nonumber\\
	&\leq\delta/2+\delta/2=\delta.\nonumber
\end{align}
This shows that the laws $\mathcal{L}^{(u^\varepsilon,\mathbf N^\varepsilon(\varphi))}=\mathbb P^\varepsilon\big(u^\varepsilon,\mathbf N^\varepsilon(\varphi)\big)^{-1}$ ($\varepsilon>0$) form a tight family of measures on $\mathscr{V}=Z\times C[0,T]$.\\

\noindent$\bf\bullet~II.~Existence~of~martingale~weak~solutions$ 
\par Now, we show the existence of martingale weak solutions. With the tightness results above, Skorohod Theorem yields a probability basis which we denote by $\big(\tilde \Omega, \mathscr{\tilde F}, \mathbb{\tilde P}\big)$ and a sequence $\Big\{\big(u,\mathbf N(\varphi)\big),\big(\tilde u^\varepsilon,\tilde{\mathbf N}^\varepsilon(\varphi)\big)\Big\}$ of $\mathscr{V}$-valued random variables on the basis such that
\begin{align}
	\mathcal{L}^{(\tilde u^\varepsilon,\tilde{\mathbf N}^\varepsilon(\varphi))}=\mathcal{L}^{(u^\varepsilon,\mathbf N^\varepsilon(\varphi))}\ \ \  for\ all\ \ \varepsilon,\nonumber
\end{align}
and
\begin{align}
	\big(\tilde u^\varepsilon,\tilde{\mathbf N}^\varepsilon(\varphi)\big)\longrightarrow \big(u,\mathbf N(\varphi)\big) \ \ \ strongly\ in\ \ \ \mathscr{V},\ \ \mathbb{\tilde P}-a.s..\nonumber
\end{align}
In particular, we have for all $\varepsilon$ that
\begin{align}\label{equal Laws}
	\mathcal{L}^{\tilde u^\varepsilon}=\mathcal{L}^{u^\varepsilon};\quad
	\mathcal{L}^{\tilde{\mathbf N}^\varepsilon(\varphi)}=\mathcal{L}^{\mathbf N^\varepsilon(\varphi)},
\end{align}
and $\mathbb{\tilde P}$-$a.s.$ that
\begin{align}
	&\tilde u^\varepsilon \longrightarrow u \ \ \ strongly\ in\ \ \ Z=C\big([0,T];H^{-s}\big)\cap L^3\big([0,T];L^3\big);\label{strong convergences of u^varepsilon}\\
	&\tilde{\mathbf N}^\varepsilon(\varphi) \longrightarrow \mathbf N(\varphi) \ \ \ strongly\ in\ \ \ C[0,T].\label{strong convergence of N^varepsilon}
\end{align}
Define
\begin{align}
	p&:=\Delta^{-1}{\rm divdiv}(u\otimes u);\quad \tilde p^\varepsilon :=\Delta^{-1}{\rm divdiv}\big((\psi_\varepsilon\ast\tilde u^\varepsilon)\otimes\tilde u^\varepsilon\big),\nonumber
\end{align}
then we also have
\begin{align}\label{equal Laws of p}
	\mathcal L(\tilde p^\varepsilon)=\mathcal L(p^\varepsilon)\ for\  all\  \varepsilon
\end{align}
and
\begin{align}\label{strong convergences of p^varepsilon}
	\tilde p^\varepsilon\longrightarrow p \ \ \ strongly\ in\ \ \ L^{3/2}\big([0,T];L^{3/2}\big),\ \ \mathbb{\tilde P}-a.s..
\end{align}
In particular, by (\ref{equal Laws}), (\ref{standard uniform bounds 1}) and (\ref{standard uniform bounds 2}), the uniform bounds  also holds for $\big\{\tilde u^\varepsilon\big\}$:
\begin{align}
		\mathbb{\tilde  E}\big\|\tilde u^\varepsilon\big\|_{C([0,T];L^2)}^r+\mathbb{\tilde  E}\big\|\tilde u^\varepsilon\big\|_{L^2([0,T];H^1)}^2
		&\leq C_1\big(\|u_0\|_{L^2},\nu,T,r\big)~.\label{standard uniform bounds 1b}
\end{align}
Hence, after subtracting a subsequence, we have
\begin{equation}\label{weak convergences of u^varepsilon}
	\left\{\begin{aligned}
	&\tilde u^\varepsilon\stackrel{\ast}{\rightharpoonup}u\ \ \ in\ \ \ L^r\big(\tilde \Omega;L^\infty([0,T];L^2)\big),\\
	&\tilde u^\varepsilon\rightharpoonup u\ \ \ in\ \ \ L^2\big(\tilde \Omega;L^2([0,T];H^1)\big),\\
	&\nabla\times\tilde u^\varepsilon\rightharpoonup\nabla\times u\ \ \ in\ \ \ L^2\big([0,T]\times\tilde \Omega\times\mathbb T^3\big).
\end{aligned}\right.
\end{equation}

In case that $\nabla\times u_0\in\mathcal M_r\big(\mathbb T^3,\mathscr B(\mathbb T^3)\big)$  and condition (\ref{con. control on vorticity}) holds, we would also have the uniform bounds (\ref{vorticity bounds for Leray regularised solution}) for $\tilde \omega^\varepsilon=\nabla\times\tilde u^\varepsilon$ with pre-given $0<\delta<1/2$. Then after subtracting a subsequence, we have
\begin{align}\label{weak convergences for vorticity}
	\nabla\tilde \omega^\varepsilon\rightharpoonup\nabla(\nabla\times u)\ \ \ in\ \ \ L^{\frac{4}{3+\delta}}\big(\tilde\Omega\times[0,T]\times\mathbb T^3\big).
\end{align}
Next, for every $R>0$, it is clear that the subset
\begin{align}
	V_R:=\left\{f\in L^2\big([0,T]\times\tilde \Omega\times\mathbb T^3\big):~\mathop{\rm ess~sup}_{t\in[0,T]}\mathbb{\tilde  E}\|f(t)\|_{L^1}\leq R\right\}\nonumber
\end{align}
is convex and strongly closed in $L^2\big([0,T]\times\tilde \Omega\times\mathbb T^3\big)$ and hence also weakly closed. By (\ref{vorticity bounds for Leray regularised solution}) and (\ref{equal Laws}) we are allowed to set
\begin{align}
	R(\varepsilon'):=\sup_{\varepsilon\in(0,\varepsilon']}\mathop{\rm ess~sup}_{t\in[0,T]}\mathbb{\tilde  E}\|\tilde \omega^\varepsilon(t)\|_{L^1}<\infty\quad for\quad\varepsilon'>0.\nonumber
\end{align}
Then by (\ref{weak convergences of u^varepsilon}) and the fact that $\left\{\tilde \omega^\varepsilon\right\}_{\varepsilon\in(0,\varepsilon']}\subset V_{R(\varepsilon')}$~, we have $\nabla\times u\in V_{R(\varepsilon')}$ as well, which means that
\begin{align}
	\mathop{\rm ess~sup}_{t\in[0,T]}\mathbb{\tilde  E}\|\nabla\times u(t)\|_{L^1}\leq R(\varepsilon')~,\quad\forall\varepsilon'>0.\nonumber
\end{align}
Sending $\varepsilon'\rightarrow0$, we see that
\begin{align}
	\mathop{\rm ess~sup}_{t\in[0,T]}\mathbb{\tilde  E}\|\nabla\times u(t)\|_{L^1}\leq\mathop{\rm limsup}_{\varepsilon\rightarrow0}\left(\mathop{\rm ess~sup}_{t\in[0,T]}\mathbb{\tilde  E}\|\tilde \omega^\varepsilon(t)\|_{L^1}\right).
\end{align}
This proves (\ref{further bounds for vorticity}).

Now, taking any test function $\varphi\in C^\infty\big(\mathbb T^3\big)$, and by (\ref{strong convergences of u^varepsilon}), we have $\mathbb{\tilde P}$-$a.s.$ that
\begin{align}
    &\ \ \ \ \ \sup_{t\in[0,T]}\big|\langle \tilde u^\varepsilon(t)-u(t),\varphi\rangle\big|\leq\|\varphi\|_{H^s}\big\|\tilde u^\varepsilon-u\big\|_{C([0,t];H^{-s})}\longrightarrow0,\nonumber\\
	&\ \ \ \ \ \sup_{t\in[0,T]}\left|\int_0^t\big\langle{\rm div}\big((\psi_\varepsilon\ast\tilde u^\varepsilon)\otimes\tilde u^\varepsilon\big)+\nabla p^\varepsilon-{\rm div}(u\otimes u)-\nabla p,\varphi \big\rangle{\rm d}s\right|\nonumber\\
	&\lesssim\|\nabla\varphi\|_{L^\infty}\Big(\big\|\psi_\varepsilon\ast(\tilde u^\varepsilon-u)\big\|_{L^2}\big\|\tilde u^\varepsilon\big\|_{L^2}+\big\|\psi_\varepsilon\ast u-u\big\|_{L^2}\big\|\tilde u\big\|_{L^2}\Big)\longrightarrow0,\nonumber\\
	&\ \ \ \ \ \sup_{t\in[0,T]}\left|\int_0^t\big\langle\Delta\tilde u^\varepsilon-\Delta u,\varphi \big\rangle{\rm d}s\right|\leq\|\Delta\varphi\|_{L^2}\big\|\tilde u^\varepsilon-u\big\|_{L^2}\longrightarrow0.\nonumber
\end{align}
Hence, let
\begin{align}\label{tilde M^varepsilon}
	\tilde{\mathbf M}^\varepsilon_t(\varphi):=\langle\tilde u^\varepsilon(t),\varphi\rangle-\langle u_0^\varepsilon,\varphi\rangle+\int_0^t\Big\langle{\rm div}\big((\psi_\varepsilon\ast\tilde u^\varepsilon)(s)\otimes\tilde u^\varepsilon(s)\big)+\nabla\tilde p^\varepsilon(s)-\nu\Delta\tilde u^\varepsilon(s),\varphi\Big\rangle{\rm d}s,
\end{align}
we have $\mathbb{\tilde P}$-$a.s.$ that
\begin{align}\label{a.s.-convergence of M^varepsilon}
	\tilde{\mathbf M}^\varepsilon(\varphi)\stackrel{C[0,T]}{\longrightarrow}\langle u(\cdot),\varphi\rangle-\langle u_0,\varphi\rangle+\int_0^\cdot\Big\langle{\rm div}\big(u(s)\otimes u(s)\big)+\nabla p(s)-\nu\Delta u(s),\varphi\Big\rangle{\rm d}s=:\mathbf M(\varphi).
\end{align}
By (\ref{equal Laws}), (\ref{equal Laws of p}) and the fact that $u^\varepsilon$ is a Martingale weak solution to (\ref{regularised sto. NS}) with initial data $u_0^\varepsilon$, it's easily seen that for each $\varepsilon$, $\tilde{\mathbf M}^\varepsilon(\varphi)$ is a continuous square-integrable martingale (w.r.t. the filtration $\{\mathscr{\tilde F}_t^\varepsilon:=\sigma(\tilde u^\varepsilon(s);0\leq s\leq t)\}_{t\geq0}$ ) with quadratic variation
\begin{align}\label{quadratic variation of tilde M^varepsilon}
	\big\langle \tilde{\mathbf M}^\varepsilon(\varphi)\big\rangle_t=\sum_{k=1}^N\int_0^t\big\langle\psi_\varepsilon\ast\varsigma_k(\tilde u^\varepsilon(s)),\varphi\big\rangle^2{\rm d}s.
\end{align}
On the other hand, by BDG-inequality, H\"older's inequality, (\ref{con. linear growth}) and (\ref{standard uniform bounds 1b}), it follows that 
\begin{align}
	\mathbb{\tilde E}\sup_{t\in[0,T]}\Big|\tilde{\mathbf M}^{\varepsilon}_{t}(\varphi)\Big|^4
	&\lesssim\mathbb{\tilde E}\left(\sum_{k=1}^N\int_0^T\big\langle\psi_\varepsilon\ast\varsigma_k(\tilde u^\varepsilon(s)),\varphi\big\rangle^2{\rm d}s\right)^2\nonumber\\
	&\lesssim\|\varphi\|_{L^2}^4\mathbb{\tilde E}\left(\sum_{k=1}^N\int_0^T\big\|\varsigma_k(\tilde u^\varepsilon(s))\big\|_{L^2}^2{\rm d}s\right)^2\nonumber\\
	&\lesssim\|\varphi\|_{L^2}^4\mathbb{\tilde E}\left(T+\big\|\tilde u^\varepsilon\big\|_{L^2}^2\right)^2
\end{align}
where the far right-hand-side is bounded uniformly in $\varepsilon$. 
Therefore together by (\ref{a.s.-convergence of M^varepsilon}) it is not hard to show the convergence (after possibly subtracting a subsequence):
\begin{align}\label{L^2-convergence of M^varepsilon}
	\lim_{\varepsilon\rightarrow0}~\mathbb{\tilde E}\sup_{t\in[0,T]}\Big|\tilde{\mathbf M}^{\varepsilon}_{t}(\varphi)-\mathbf M_{t}(\varphi)\Big|^2=0~,
\end{align}
as for any $\delta>0$ one could write
\begin{align}
	\mathbb{\tilde E}\sup_{t\in[0,T]}\Big|\tilde{\mathbf M}^{\varepsilon}_{t}(\varphi)-\mathbf M_{t}(\varphi)\Big|^2
	&\leq\int_{\{\|\tilde{\mathbf M}^{\varepsilon}(\varphi)-\mathbf M(\varphi)\|_{C[0,T]}>\delta\}}\big\|\tilde{\mathbf M}^{\varepsilon}(\varphi)-\mathbf M(\varphi)\big\|_{C[0,T]}^2{\rm d}\mathbb{\tilde P}+\delta^2\nonumber\\
	&\leq\mathbb{\tilde P}\left(\|\tilde{\mathbf M}^{\varepsilon}(\varphi)-\mathbf M(\varphi)\|_{C[0,T]}>\delta\right)^{\frac{1}{2}}\cdot\left(\mathbb{\tilde E}\big\|\tilde{\mathbf M}^{\varepsilon}(\varphi)-\mathbf M(\varphi)\big\|_{C[0,T]}^4\right)^{\frac{1}{2}}+\delta^2.\nonumber
\end{align} 
Next, we show
\begin{align}
	\sum_{k=1}^N\int_{0}^{t}\big\langle\psi_\varepsilon\ast\varsigma_k(\tilde u^\varepsilon(s)),\varphi\big\rangle^2{\rm d}s\stackrel{L^1(\tilde\Omega)}{\longrightarrow}\sum_{k=1}^\infty\int_0^t\big\langle\varsigma_k(u(s)),\varphi\big\rangle^2{\rm d}s.\label{convergence of quadratic variation of M^varepsilon}
\end{align}
By H\"older's inequality, (\ref{con. linear growth}), (\ref{standard uniform bounds 1b}), we have
\begin{align}
	&\mathbb{\tilde E}\left|\sum_{k=1}^N\int_0^t\big\langle\psi_\varepsilon\ast\varsigma_k(\tilde u^\varepsilon(s)),\varphi\big\rangle^2{\rm d}s-\sum_{k=1}^\infty\int_0^t\big\langle\varsigma_k(u(s)),\varphi\big\rangle^2{\rm d}s\right|\nonumber\\
	\leq &\mathbb{\tilde E}\left|\sum_{k=1}^N\int_0^t\big\langle\psi_\varepsilon\ast\varsigma_k(\tilde u^\varepsilon(s)),\varphi\big\rangle^2{\rm d}s-\sum_{k=1}^N\int_0^t\big\langle\varsigma_k(u(s)),\varphi\big\rangle^2{\rm d}s\right|+\mathbb{\tilde E}\sum_{k=N+1}^\infty\int_0^t\big\langle\varsigma_k(u(s)),\varphi\big\rangle^2{\rm d}s\nonumber\\
	\leq &\|\varphi\|_{L^2}^2\mathbb{\tilde E}\sum_{k=1}^N\int_0^T\big\|\psi_\varepsilon\ast\varsigma_k(\tilde u^\varepsilon(s))-\varsigma_k(u(s))\big\|_{L^2}\Big(\big\|\varsigma_k(\tilde u^\varepsilon(s))\big\|_{L^2}+\big\|\varsigma_k(u(s))\big\|_{L^2}\Big){\rm d}s\nonumber\\
	&\quad +\|\varphi\|_{L^2}^2\mathbb{\tilde E}\int_0^T\sum_{k=N+1}^\infty\big\|\varsigma_k(u(s))\big\|_{L^2}^2{\rm d}s\nonumber\\
	\lesssim &\|\varphi\|_{L^2}^2\left(\mathbb{\tilde E}\sum_{k=1}^N\int_0^T\big\|\psi_\varepsilon\ast\varsigma_k(\tilde u^\varepsilon(s))-\varsigma_k(u(s))\big\|_{L^2}^2{\rm d}s\right)^{\frac{1}{2}}\left[\mathbb{\tilde E}\int_0^T\left(2+\big\|\tilde u^\varepsilon(s)\big\|_{L^2}^2+\|u(s)\|_{L^2}^2\right){\rm d}s\right]^{\frac{1}{2}}\nonumber\\
	&\quad +\left(\sup_{u\in L^2}\frac{\sum_{k>N}\left\|\varsigma_k(u)\right\|_{L^2}^2}{1+\|u\|_{L^2}^2}\right)\|\varphi\|_{L^2}^2\mathbb{\tilde E}\int_0^T\left(1+\|u(s)\|_{L^2}^2\right){\rm d}s\nonumber\\
	\lesssim &\|\varphi\|_{L^2}^2\left[T(2+C_1)+\mathbb{\tilde E}\|u\|_{L^2}^2\right]^{\frac{1}{2}}\left(\mathbb{\tilde E}\sum_{k=1}^N\int_0^T\big\|\psi_\varepsilon\ast\varsigma_k(\tilde u^\varepsilon(s))-\varsigma_k(u(s))\big\|_{L^2}^2{\rm d}s\right)^{\frac{1}{2}}\nonumber\\
	&\quad +\|\varphi\|_{L^2}^2\left(T+\mathbb{\tilde E}\|u\|_{L^2}^2\right)\left(\sup_{u\in L^2}\frac{\sum_{k>N}\left\|\varsigma_k(u)\right\|_{L^2}^2}{1+\|u\|_{L^2}^2}\right).
\end{align}

{\noindent Therefore by (\ref{con. uniform decay of the tails}) it suffices to show that, as $\varepsilon\rightarrow0$,}
\begin{align}
	\mathbb{\tilde E}\int_0^T\sum_{k=1}^\infty\big\|\psi_{\varepsilon}\ast\varsigma_k(u(s))-\varsigma_k(u(s))\big\|_{L^2}^2{\rm d}s\longrightarrow0~;\label{control of convergence of M(varphi) 1}\\
	\mathbb{\tilde E}\int_0^T\sum_{k=1}^\infty\big\|\varsigma_k(\tilde u^{\varepsilon}(s))-\varsigma_k(u(s))\big\|_{L^2}^2{\rm d}s\longrightarrow0~.\label{control of convergence of M(varphi) 2}
\end{align}
(\ref{control of convergence of M(varphi) 1}) follows by Dominated Convergence Theorem since\begin{align}
	&\sum_{k=1}^\infty\big\|\psi_{\varepsilon}\ast\varsigma_k(u(\cdot))-\varsigma_k(u(\cdot))\big\|_{L^2}^2\longrightarrow0,~~{\rm d}\mathbb{\tilde P}\times{\rm d}t-a.e.\quad \big(again~by~Dominated~Convergence\big);\nonumber\\
	&\sum_{k=1}^\infty\big\|\psi_{\varepsilon}\ast\varsigma_k(u(\cdot))-\varsigma_k(u(\cdot))\big\|_{L^2}^2\lesssim1+\|u(\cdot)\|_{L^2}^2\ \big(\in L^1(\tilde\Omega\times[0,T])\big)\quad \big(by~(\ref{con. linear growth})~and~(\ref{standard uniform bounds 1b})\big).\nonumber
\end{align}
For (\ref{control of convergence of M(varphi) 2}), we split the infinite sum by an integer $N'\in\mathbb N$ to be chosen later:
\begin{align}
	&\mathbb{\tilde E}\int_0^T\sum_{k=1}^\infty\big\|\varsigma_k(\tilde u^{\varepsilon}(s))-\varsigma_k(u(s))\big\|_{L^2}^2{\rm d}s\nonumber\\
	=&\mathbb{\tilde E}\int_0^T\sum_{k=1}^{N'}\big\|\varsigma_k(\tilde u^{\varepsilon}(s))-\varsigma_k(u(s))\big\|_{L^2}^2{\rm d}s+\mathbb{\tilde E}\int_0^T\sum_{k=N'+1}^\infty\big\|\varsigma_k(\tilde u^{\varepsilon}(s))-\varsigma_k(u(s))\big\|_{L^2}^2{\rm d}s\nonumber\\
	\leq &\sum_{k=1}^{N'}\mathbb{\tilde E}\int_0^T\big\|\varsigma_k(\tilde u^{\varepsilon}(s))-\varsigma_k(u(s))\big\|_{L^2}^2{\rm d}s\nonumber\\
	&+2\left(\sup_{u\in L^2}\frac{\sum_{k>N'}\left\|\varsigma_k(u)\right\|_{L^2}^2}{1+\|u\|_{L^2}^2}\right)\mathbb{\tilde E}\int_0^T\left(2+\|u(s)\|_{L^2}^2+\big\|\tilde u^{\varepsilon}(s)\big\|_{L^2}^2\right){\rm d}s\nonumber\\
	\leq &\sum_{k=1}^{N'}\mathbb{\tilde E}\int_0^T\big\|\varsigma_k(\tilde u^{\varepsilon}(s))-\varsigma_k(u(s))\big\|_{L^2}^2{\rm d}s+2\left[T(2+C_1)+\mathbb{\tilde E}\|u\|_{L^2}^2\right]\left(\sup_{u\in L^2}\frac{\sum_{k>N'}\left\|\varsigma_k(u)\right\|_{L^2}^2}{1+\|u\|_{L^2}^2}\right),\nonumber
\end{align}
where we have used (\ref{standard uniform bounds 1b}) to get the last inequality. Then sending $\varepsilon\rightarrow0$, the first term on the right-hand-side tends to zero by Dominated Convergence Theorem since
\begin{align}
	&\big\|\varsigma_k(\tilde u^{\varepsilon}(\cdot))-\varsigma_k(u(\cdot))\big\|_{L^2}^2\longrightarrow0,~~{\rm d}\mathbb{\tilde P}\times{\rm d}t-a.e.\quad \big(by~(\ref{strong convergences of u^varepsilon})~and~continuity~of~\varsigma_k~'s\big);\nonumber\\
	&\big\|\varsigma_k(\tilde u^{\varepsilon}(\cdot))-\varsigma_k(u(\cdot))\big\|_{L^2}^2\lesssim1+\|u(\cdot)\|_{L^2}^2+\big\|\tilde u^{\varepsilon}(\cdot)\big\|_{L^2}^2\ \big(\in L^1(\tilde\Omega\times[0,T])\big)\quad \big(by~(\ref{con. linear growth})~and~(\ref{standard uniform bounds 1b})\big),\nonumber
\end{align}
and we deduce
\begin{align}
	\lim_{\varepsilon\rightarrow0}\mathbb{\tilde E}\int_0^T\sum_{k=1}^\infty\big\|\varsigma_k(\tilde u^{\varepsilon}(s))-\varsigma_k(u(s))\big\|_{L^2}^2{\rm d}s\leq 2\left[T(2+C_1)+\mathbb{\tilde E}\|u\|_{L^2}^2\right]\left(\sup_{u\in L^2}\frac{\sum_{k>N'}\left\|\varsigma_k(u)\right\|_{L^2}^2}{1+\|u\|_{L^2}^2}\right).\nonumber
\end{align}
As this is valid for arbitrary $N'\in\mathbb N$ and the left hand side does not depend on $N'$, we see that (\ref{control of convergence of M(varphi) 2}) follows by (\ref{con. uniform decay of the tails}). And finally (\ref{convergence of quadratic variation of M^varepsilon}) is proved. 

\par Now take any bounded continuous function $\phi$ on $Z$ and by (\ref{quadratic variation of M^varepsilon}) we can write for $0\leq t_1<t_2\leq T$ that
\begin{align}
	&\mathbb{\tilde E}\Big[\left(\mathbf M^{\varepsilon}_{t_2}(\varphi)-\mathbf M^{\varepsilon}_{t_1}(\varphi)\right)\phi\left(\tilde u^\varepsilon|_{[0,t_1]}\right)\Big]=0,\label{martingale property of M^variepsilon}\\
	&\mathbb{\tilde E}\bigg[\Big(\mathbf M^{\varepsilon}_{t_2}(\varphi)^2-\mathbf M^{\varepsilon}_{t_1}(\varphi)^2\Big)\phi\left(\tilde u^\varepsilon|_{[0,t_1]}\right)\bigg]\nonumber\\
	=&\mathbb{\tilde E}\Bigg[\left(\sum_{k=1}^N\int_{t_1}^{t_2}\big\langle\psi_\varepsilon\ast\varsigma_k(\tilde u^\varepsilon(s)),\varphi\big\rangle^2{\rm d}s\right)\phi\left(\tilde u^\varepsilon|_{[0,t_1]}\right)\Bigg].\label{quadratic variation equality of M^varepsilon}
\end{align}
Then by using (\ref{L^2-convergence of M^varepsilon}), (\ref{convergence of quadratic variation of M^varepsilon}) and (\ref{strong convergences of u^varepsilon}), one can take limits in (\ref{martingale property of M^variepsilon}) and (\ref{quadratic variation equality of M^varepsilon}) to see that the process
\begin{align}
	\mathbf M(\varphi):=\langle u(\cdot),\varphi\rangle-\langle u_0,\varphi\rangle+\int_0^\cdot\Big\langle{\rm div}\big(u(s)\otimes u(s)\big)+\nabla p(s)-\nu\Delta u(s),\varphi\Big\rangle{\rm d}s\nonumber
\end{align}
is a continuous square-integrable martingale (w.r.t. the filtration $\big\{\mathscr{\tilde F}_t:=\sigma( u(s);0\leq s\leq t)\big\}_{t\geq0}$ ) with quadratic variation
\begin{align}
	\big\langle \mathbf M(\varphi)\big\rangle=\sum_{k=1}^\infty\int_0^t\big\langle\varsigma_k(u(s)),\varphi\big\rangle^2{\rm d}s.\nonumber
\end{align}
\\

\noindent$\bf\bullet~III.~Derivation~of~local~energy~inequality$
\par It remains to show (\ref{LEI-I}). The local energy balance of $u^\varepsilon$ is given by (\ref{LEE-a for u^varepsilon}) in Proposition \ref{Prop. LEE for Leray regularization}. Together with (\ref{equal Laws}), one should see that the pair $\big(\tilde u^\varepsilon,\tilde{\mathbf N}^\varepsilon(\varphi)\big)$ should satisfies the equality $\tilde{\mathbb P}$-$a.s.$ with (scalar) test function $\varphi\in C^\infty_c([0,T]\times\mathbb T^3)$ :
\begin{align}\label{LEE for u^varepsilon}
	&\ \ \ \ \ \ \int_{\mathbb{T}^3}\big|\tilde u^\varepsilon(t)\big|^2\varphi(t){\rm d}x+2\nu\int_0^t\int_{\mathbb{T}^3}\big|\nabla \tilde u^\varepsilon(s)\big|^2\varphi(s){\rm d}x{\rm d}s\nonumber\\
	&=\underbrace{\int_0^t\int_{\mathbb{T}^3}\big|\tilde u^\varepsilon(s)\big|^2\Big(\partial_\tau\varphi(s)+\nu\Delta\varphi(s)\Big){\rm d}x{\rm d}s+\sum_{k=1}^N\int_0^t\int_{\mathbb{T}^3}|\psi_\varepsilon\ast\varsigma_k(\tilde u^\varepsilon(s))|^2\varphi(s){\rm d}x{\rm d}s}_{I^\varepsilon_1(t)}\nonumber\\
	&\ \ \ \ \ +\underbrace{\int_0^t\int_{\mathbb{T}^3}\big|\tilde u^\varepsilon(s)\big|^2\big(\psi_\varepsilon\ast\tilde u^\varepsilon(s)\big)\cdot\nabla \varphi(s){\rm d}x{\rm d}s}_{I^\varepsilon_2(t)}
	+\underbrace{\int_0^t\int_{\mathbb{T}^3}2\tilde p^\varepsilon(s)\tilde u^\varepsilon(s)\cdot\nabla\varphi(s){\rm d}x{\rm d}s}_{I^\varepsilon_3(t)}+\tilde{\mathbf N}^\varepsilon_t(\varphi),\quad \forall t\in[0,T].
\end{align}
Here, again by (\ref{equal Laws}), $\tilde{\mathbf N}^\varepsilon(\varphi)$ is a continuous square-integrable martingale (w.r.t. the filtration $\big\{\mathscr{\tilde F}_t^\varepsilon:=\sigma(\tilde u^\varepsilon(s);0\leq s\leq t)\big\}_{t\geq0}$ ) with quadratic variation
\begin{align}\label{quadratic variation of tilde N^varepsilon}
	\big\langle \tilde{\mathbf N}^\varepsilon(\varphi)\big\rangle_t=4\sum_{k=1}^N\int_0^t\left(\int_{\mathbb{T}^3}\psi_\varepsilon\ast\varsigma_k(\tilde u^\varepsilon(s))\cdot\tilde u^\varepsilon(s)\varphi(s){\rm d}x\right)^2{\rm d}s.
\end{align}

\par To this end, we deduce compactness for each term on the righthand side of (\ref{LEE for u^varepsilon}). We first deal with $I^\varepsilon_1$. By using Fatou's lemma, H\"older's inequality and (\ref{con. linear growth}), we have
\begin{align}
	&\sup_{t\in[0,T]}\left|\sum_{k=1}^N\int_0^t\int_{\mathbb{T}^3}|\psi_\varepsilon\ast\varsigma_k(\tilde u^\varepsilon(s))|^2\varphi(s){\rm d}x{\rm d}s-\sum_{k=1}^\infty\int_0^t\int_{\mathbb{T}^3}|\varsigma_k(u(s))|^2\varphi(s){\rm d}x{\rm d}s\right|\nonumber\\
	\leq &\sum_{k>N}\int_0^T\int_{\mathbb{T}^3}|\psi_\varepsilon\ast\varsigma_k(\tilde u^\varepsilon(s))|^2\varphi(s){\rm d}x{\rm d}s\nonumber\\
	&+\sum_{k=1}^\infty\int_0^T\int_{\mathbb{T}^3}\big|\psi_\varepsilon\ast\varsigma_k(\tilde u^\varepsilon(s))+\varsigma_k(u(s))\big|\big|\psi_\varepsilon\ast\varsigma_k(\tilde u^\varepsilon(s))-\varsigma_k(u(s))\big|\varphi(s){\rm d}x{\rm d}s\nonumber\\
	\leq &\|\varphi\|_{L^\infty}\int_0^T\sum_{k>N}\big\|\psi_\varepsilon\ast\varsigma_k(\tilde u^\varepsilon(s))\big\|_{L^2}^2{\rm d}s\nonumber\\
	&+\|\varphi\|_{L^\infty}\left(\int_0^T\sum_{k=1}^\infty\big\|\psi_\varepsilon\ast\varsigma_k(\tilde u^\varepsilon(s))+\varsigma_k(u(s))\big\|_{L^2}^2{\rm d}s \right)^{\frac{1}{2}}\left(\int_0^T\sum_{k=1}^\infty\big\|\psi_\varepsilon\ast\varsigma_k(\tilde u^\varepsilon(s))-\varsigma_k(u(s))\big\|_{L^2}^2{\rm d}s\right)^{\frac{1}{2}}\nonumber\\
	\lesssim &\left(\sup_{u\in L^2}\frac{\sum_{k>N}\big\|\varsigma_k(u)\big\|_{L^2}^2}{1+\|u\|_{L^2}^2}\right)\int_0^T\big(1+\|\tilde u^\varepsilon(s)\|_{L^2}^2\big){\rm d}s\nonumber\\
	&+\left[\int_0^T\big(1+\|\tilde u^\varepsilon(s)\|_{L^2}^2\big){\rm d}s+\int_0^T\big(1+\|u(s)\|_{L^2}^2\big){\rm d}s\right]^{\frac{1}{2}}\left(\int_0^T\sum_{k=1}^\infty\big\|\psi_\varepsilon\ast\varsigma_k(\tilde u^\varepsilon(s))-\varsigma_k(u(s))\big\|_{L^2}^2{\rm d}s\right)^{\frac{1}{2}}.\nonumber
\end{align}
By (\ref{strong convergences of u^varepsilon}), we see that $\mathbb{\tilde P}$-$a.s.$, $\displaystyle\int_0^T\big(1+\|\tilde u^\varepsilon(s)\|_{L^2}^2\big){\rm d}s$ is bounded uniformly in $\varepsilon$. Hence we arrive at
\begin{align}\label{a.s.-convergence of noise term in I1}
	&\sup_{t\in[0,T]}\left|\sum_{k=1}^N\int_0^t\int_{\mathbb{T}^3}|\psi_\varepsilon\ast\varsigma_k(\tilde u^\varepsilon(s))|^2\varphi(s){\rm d}x{\rm d}s-\sum_{k=1}^\infty\int_0^t\int_{\mathbb{T}^3}|\varsigma_k(u(s))|^2\varphi(s){\rm d}x{\rm d}s\right|\nonumber\\
	\lesssim &\left(\sup_{u\in L^2}\frac{\sum_{k>N}\big\|\varsigma_k(u)\big\|_{L^2}^2}{1+\|u\|_{L^2}^2}\right)+\left(\int_0^T\sum_{k=1}^\infty\big\|\psi_\varepsilon\ast\varsigma_k(\tilde u^\varepsilon(s))-\varsigma_k(u(s))\big\|_{L^2}^2{\rm d}s\right)^{\frac{1}{2}}\nonumber\\
	\lesssim &\left(\sup_{u\in L^2}\frac{\sum_{k>N}\big\|\varsigma_k(u)\big\|_{L^2}^2}{1+\|u\|_{L^2}^2}\right)+\left(\int_0^T\sum_{k=1}^\infty\big\|\varsigma_k(\tilde u^\varepsilon(s))-\varsigma_k(u(s))\big\|_{L^2}^2{\rm d}s\right)^{\frac{1}{2}}\nonumber\\
	&+\left(\int_0^T\sum_{k=1}^\infty\big\|\psi_\varepsilon\ast\varsigma_k(u(s))-\varsigma_k(u(s))\big\|_{L^2}^2{\rm d}s\right)^{\frac{1}{2}}
\end{align}
where, by (\ref{con. uniform decay of the tails}), (\ref{control of convergence of M(varphi) 1}) and (\ref{control of convergence of M(varphi) 2}), the far right-hand-side converges to zero $\mathbb{\tilde P}$-$a.s.$ as $\varepsilon\rightarrow0$ (after possibly subtracting a subsequence). Now by (\ref{strong convergences of u^varepsilon}) and (\ref{a.s.-convergence of noise term in I1}), we see that $\mathbb{\tilde P}-a.s.$,
\begin{align}\label{a.s.-convergence of I1}
	I^\varepsilon_1\stackrel{C[0,T]}{\longrightarrow}\int_0^\cdot\int_{\mathbb{T}^3}\big|u(s)\big|^2\Big(\partial_\tau\varphi(s)+\nu\Delta\varphi(s)\Big){\rm d}x{\rm d}s+\sum_{k=1}^\infty\int_0^\cdot\int_{\mathbb{T}^3}|\varsigma_k(u(s))|^2\varphi(s){\rm d}x{\rm d}s.
\end{align}
For $I^\varepsilon_2$, we have $\mathbb{\tilde P}$-$a.s.$ that
\begin{align}\label{a.s.-convergence of I2}
	&\sup_{t\in[0,T]}\left|I^\varepsilon_2(t)-\int_0^t\int_{\mathbb{T}^3}\big|u(s)\big|^2u\cdot\nabla \varphi(s){\rm d}x{\rm d}s\right|\nonumber\\
	&\leq\|\nabla\varphi\|_{L^\infty}\int_0^T\int_{\mathbb{T}^3}\Big|\underbrace{\big|\tilde u^\varepsilon(s)\big|^2\psi_\varepsilon\ast\tilde u^\varepsilon(s)-|u(s)|^2u(s)}_{=|\tilde u^\varepsilon-u|^2\psi_\varepsilon\ast\tilde u^\varepsilon+2u\cdot(\tilde u^\varepsilon-u)\psi_\varepsilon\ast\tilde u^\varepsilon+|u|^2(\psi_\varepsilon\ast\tilde u^\varepsilon-u)}\Big|{\rm d}x{\rm d}s\nonumber\\
	&\leq\|\nabla\varphi\|_{L^\infty}\Big(\big\|\tilde u^\varepsilon-u\big\|_{L^3}^2\big\|\tilde u^\varepsilon\big\|_{L^3}+2\|u\|_{L^3}\big\|\tilde u^\varepsilon-u\big\|_{L^3}\big\|\tilde u^\varepsilon\big\|_{L^3}+\|u\|_{L^3}^2\big\|\psi_\varepsilon\ast\tilde u^\varepsilon-u\big\|_{L^3}\Big)\rightarrow0
\end{align}
where we have used $\big\|\psi_\varepsilon\ast\tilde u^\varepsilon-u\big\|_{L^3}\leq\big\|\psi_\varepsilon\ast(\tilde u^\varepsilon-u)\big\|_{L^3}+\big\|\psi_\varepsilon\ast u-u\big\|_{L^3}$ and (\ref{strong convergences of u^varepsilon}). 
And by (\ref{strong convergences of u^varepsilon}) and (\ref{strong convergences of p^varepsilon}), we have for $I^\varepsilon_3$ that
\begin{align}\label{a.s.-convergence of I3}
	&\sup_{t\in[0,T]}\left|I^\varepsilon_3(t)-\int_0^t\int_{\mathbb{T}^3}2p(s)u(s)\cdot\nabla\varphi(s){\rm d}x{\rm d}s\right|\nonumber\\
	&\leq2\|\nabla\varphi\|_{L^\infty}\int_0^T\int_{\mathbb{T}^3}\Big|\tilde p^\varepsilon(s)\tilde u^\varepsilon(s)-p(s)u(s)\Big|{\rm d}x{\rm d}s\nonumber\\
	&\leq2\|\nabla\varphi\|_{L^\infty}\Big(\big\|\tilde p^\varepsilon\big\|_{L^{3/2}}\big\|\tilde u^\varepsilon-u\big\|_{L^3}+\|u\|_{L^3}\big\|\tilde p^\varepsilon-p\big\|_{L^{3/2}}\Big)\rightarrow0,\ \ \mathbb{\tilde P}-a.s..
\end{align}
Next we show that $\big\{I^\varepsilon_1\big\}$, $\big\{I^\varepsilon_2\big\}$ and $\big\{I^\varepsilon_3\big\}$ are uniformly integrable families of $C[0,T]$-valued random variables. For $I^\varepsilon_1$, by (\ref{con. linear growth}), we obviously have
\begin{align}
	\mathbb{\tilde E}\sup_{t\in[0,T]}\big|I^\varepsilon_1(t)\big|^p\leq T^p\|\varphi\|_{W^{2,\infty}}^p\mathbb{\tilde E}\left(1+\big\|\tilde u^\varepsilon\big\|_{L^\infty([0,T];L^2)}^2\right)^{p}\nonumber
\end{align}
which is uniformly bounded in $\varepsilon$ by (\ref{standard uniform bounds 1b}).
For $I^\varepsilon_2$, by H\"older's inequality, Sobolev imbedding $H^{3/4}\hookrightarrow L^4$ and Sobolev interpolation $H^{3/4}=\big[L^2,H^1\big]_{1/4}$, we have
\begin{align}
	\sup_{t\in[0,T]}\big|I^\varepsilon_2(t)\big|
	&\leq\|\nabla\varphi\|_{L^{\infty}}\int_0^T\int_{\mathbb{T}^3}\Big|\big|\tilde u^\varepsilon(s)\big|^2\psi_\varepsilon\ast\tilde u^\varepsilon(s)\Big|{\rm d}x{\rm d}s\nonumber\\
	&\leq\|\nabla\varphi\|_{L^{\infty}}\int_0^T\big\|\psi_\varepsilon\ast\tilde u^\varepsilon(s)\big\|_{L^2}\big\|\tilde u^\varepsilon(s)\big\|_{L^4}^2{\rm d}s\nonumber\\
	&\lesssim\|\nabla\varphi\|_{L^{\infty}}\big\|\tilde u^\varepsilon\big\|_{L^\infty([0,T];L^2)}\int_0^T\big\|\tilde u^\varepsilon(s)\big\|_{H^{3/4}}^2{\rm d}s\nonumber\\
	&\lesssim\|\nabla\varphi\|_{L^{\infty}}\big\|\tilde u^\varepsilon\big\|_{L^\infty([0,T];L^2)}\int_0^T\big\|\tilde u^\varepsilon(s)\big\|_{L^2}^{1/2}\big\|\nabla\tilde u^\varepsilon(s)\big\|_{L^2}^{3/2}{\rm d}s\nonumber\\
	&\lesssim\|\nabla\varphi\|_{L^{\infty}}\big\|\tilde u^\varepsilon\big\|_{L^\infty([0,T];L^2)}^{3/2}\big\|\tilde u^\varepsilon\big\|_{L^2([0,T];H^1)}^{3/2}T^{1/4}.\nonumber
\end{align}
Hence, we may have, say,
\begin{align}
	\mathbb{\tilde E}\sup_{t\in[0,T]}\big|I^\varepsilon_2(t)\big|^{10/9}
	&\lesssim T^{5/18}\|\nabla\varphi\|_{L^{\infty}}^{10/9}\mathbb{\tilde E}\Big(\big\|\tilde u^\varepsilon\big\|_{L^\infty([0,T];L^2)}^{5/3}\big\|\tilde u^\varepsilon\big\|_{L^2([0,T];H^1)}^{5/3}\Big)\nonumber\\
	&\lesssim T^{5/18}\|\nabla\varphi\|_{L^{\infty}}^{10/9}\Big(\mathbb{\tilde E}\big\|\tilde u^\varepsilon\big\|_{L^\infty([0,T];L^2)}^{10}\Big)^{1/6}\Big(\mathbb{\tilde E}\big\|\tilde u^\varepsilon\big\|_{L^2([0,T];H^1)}^{2}\Big)^{5/6}.\nonumber
\end{align}
Here, the far right-hand-side is again uniformly bounded in $\varepsilon$ by (\ref{standard uniform bounds 1b}).
$I^\varepsilon_3$ can be estimated similarly as $I^\varepsilon_2$, noting the $L^2$-boundedness of the operator $\Delta^{-1}{\rm divdiv}$. Hence, we see that $\big\{I^\varepsilon_1\big\}$, $\big\{I^\varepsilon_2\big\}$ and $\big\{I^\varepsilon_3\big\}$ are uniformly integrable (here, we have to choose $r\geq10$ in (\ref{standard uniform bounds 1b}) ). And now, we have the convergences
\begin{equation}\label{convergences of I1 I2 I3}
	\mathbb{\tilde P}-a.s.\ and\ strongly\ in\ \ \ L^1(\tilde \Omega)\ \left\{\begin{aligned}
		&I^\varepsilon_1\stackrel{C[0,T]}{\longrightarrow}\int_0^\cdot\int_{\mathbb{T}^3}\big|u(s)\big|^2\Big(\partial_\tau\varphi(s)+\nu\Delta\varphi(s)\Big){\rm d}x{\rm d}s\\
		&\ \ \ \ \ \ \ \ \ \ \ \ \ \ +\sum_{k=1}^\infty\int_0^\cdot\int_{\mathbb{T}^3}|\varsigma_k(u(s))|^2\varphi(s){\rm d}x{\rm d}s,\\
	&I^\varepsilon_2\stackrel{C[0,T]}{\longrightarrow}\int_0^\cdot\int_{\mathbb{T}^3}\big|u(s)\big|^2u\cdot\nabla \varphi(s){\rm d}x{\rm d}s,\\
	&I^\varepsilon_3\stackrel{C[0,T]}{\longrightarrow}\int_0^\cdot\int_{\mathbb{T}^3}2p(s)u(s)\cdot\nabla\varphi(s){\rm d}x{\rm d}s.
	\end{aligned}\right.
\end{equation}

\par Next, we deal with the martingale term $\tilde{\mathbf N}^\varepsilon(\varphi)$. Recalling (\ref{strong convergence of N^varepsilon}), we already have $\tilde{\mathbb P}$-$a.s.$ convergence in space $C[0,T]$. Now we apply a similar uniform integrability argument as when we deal with $\tilde{\mathbf M}^\varepsilon(\varphi)$. By (\ref{quadratic variation of tilde N^varepsilon}), BDG-inequality, H\"older's inequality and (\ref{con. linear growth}), we have the estimate
\begin{align}
	\mathbb{\tilde E}\sup_{t\in[0,T]}\big|\tilde{\mathbf N}^\varepsilon_t(\varphi)\big|^{4}
	&\lesssim\mathbb{\tilde E}\left[\sum_{k=1}^N\int_0^T\left(\int_{\mathbb{T}^3}\psi_\varepsilon\ast\varsigma_k(\tilde u^\varepsilon(s))\cdot\tilde u^\varepsilon(s)\varphi(s){\rm d}x\right)^2{\rm d}s\right]^2\nonumber\\
	&\lesssim\|\varphi\|_{L^\infty}^4\mathbb{\tilde E}\left(\sum_{k=1}^N\int_0^T\big\|\varsigma_k(\tilde u^\varepsilon(s))\big\|_{L^2}^2\big\|\tilde u^\varepsilon(s)\big\|_{L^2}^2{\rm d}s\right)^2\nonumber\\
	&\lesssim\|\varphi\|_{L^\infty}^4\mathbb{\tilde E}\left[\int_0^T\left(1+\big\|\tilde u^\varepsilon(s)\big\|_{L^2}^2\right)\big\|\tilde u^\varepsilon(s)\big\|_{L^2}^2{\rm d}s\right]^2\nonumber
\end{align}
where by (\ref{standard uniform bounds 1b}) the right-hand-side is uniformly bounded in $\varepsilon$. 
Therefore together by (\ref{a.s.-convergence of M^varepsilon}) it is not hard to show the convergence (after possibly subtracting a subsequence):
\begin{align}\label{L^2-convergence of N^varepsilon}
	\lim_{\varepsilon\rightarrow0}~\mathbb{\tilde E}\sup_{t\in[0,T]}\Big|\tilde{\mathbf N}^{\varepsilon}_{t}(\varphi)-\mathbf N_{t}(\varphi)\Big|^2=0~.
\end{align} 
Next, we show
\begin{align}\label{convergence of quadratic variation of N^varepsilon}
    \sum_{k=1}^N\int_0^t\left(\int_{\mathbb{T}^3}\psi_\varepsilon\ast\varsigma_k(\tilde u^\varepsilon(s))\cdot\tilde u^\varepsilon(s)\varphi(s){\rm d}x\right)^2{\rm d}s\stackrel{L^1(\tilde\Omega)}{\longrightarrow}\sum_{k=1}^\infty\int_0^t\left(\int_{\mathbb{T}^3}\varsigma_k(u(s))\cdot u(s)\varphi(s){\rm d}x\right)^2{\rm d}s.
\end{align}
Firstly, by H\"older's inequality and (\ref{con. uniform decay of the tails}) we see that
\begin{align}
	&\mathbb{\tilde E}\sum_{k=N+1}^\infty\int_0^T\left(\int_{\mathbb{T}^3}\varsigma_k(u(s))\cdot u(s)\varphi(s){\rm d}x\right)^2{\rm d}s\nonumber\\
	\leq &\|\varphi\|_{L^\infty}^2\mathbb{\tilde E}\sum_{k=N+1}^\infty\int_0^T\big\|\varsigma_k(u(s))\big\|_{L^2}^2\| u(s)\|_{L^2}^2{\rm d}s\nonumber\\
	\leq &\|\varphi\|_{L^\infty}^2\left(\sup_{u\in L^2}\frac{\sum_{k>N}\big\|\varsigma_k(u)\big\|_{L^2}^2}{1+\|u\|_{L^2}^2}\right)\mathbb{\tilde E}\int_0^T\left(1+\|u(s)\|_{L^2}^2\right)\| u(s)\|_{L^2}^2{\rm d}s\longrightarrow0~,\quad as~\varepsilon\rightarrow0.
\end{align}
Hence it remains to show
\begin{align}
	\lim_{\varepsilon\rightarrow0}\mathbb{\tilde E}\left|\sum_{k=1}^N\int_0^t\left(\int_{\mathbb{T}^3}\psi_\varepsilon\ast\varsigma_k(\tilde u^\varepsilon(s))\cdot\tilde u^\varepsilon(s)\varphi(s){\rm d}x\right)^2{\rm d}s-\sum_{k=1}^N\int_0^t\left(\int_{\mathbb{T}^3}\varsigma_k(u(s))\cdot u(s)\varphi(s){\rm d}x\right)^2{\rm d}s\right|=0.\nonumber
\end{align}
By H\"older's inequality and (\ref{con. linear growth}) we have
\begin{align}
	&\mathbb{\tilde E}\left|\sum_{k=1}^N\int_0^t\left(\int_{\mathbb{T}^3}\psi_\varepsilon\ast\varsigma_k(\tilde u^\varepsilon(s))\cdot\tilde u^\varepsilon(s)\varphi(s){\rm d}x\right)^2{\rm d}s-\sum_{k=1}^N\int_0^t\left(\int_{\mathbb{T}^3}\varsigma_k(u(s))\cdot u(s)\varphi(s){\rm d}x\right)^2{\rm d}s\right|\nonumber\\
	\leq &\|\varphi\|_{L^\infty}^2\mathbb{\tilde E}\sum_{k=1}^N\int_0^T\big\|\psi_\varepsilon\ast\varsigma_k(\tilde u^\varepsilon(s))\cdot\tilde u^\varepsilon(s)-\varsigma_k(u(s))\cdot u(s)\big\|_{L^1}\big\|\psi_\varepsilon\ast\varsigma_k(\tilde u^\varepsilon(s))\cdot\tilde u^\varepsilon(s)+\varsigma_k(u(s))\cdot u(s)\big\|_{L^1}{\rm d}s\nonumber\\
	\lesssim &\|\varphi\|_{L^\infty}^2\left(\mathbb{\tilde E}\sum_{k=1}^N\int_0^T\big\|\psi_\varepsilon\ast\varsigma_k(\tilde u^\varepsilon(s))\cdot\tilde u^\varepsilon(s)-\varsigma_k(u(s))\cdot u(s)\big\|_{L^1}^2{\rm d}s\right)^{\frac{1}{2}}\nonumber\\
	&\quad\quad\quad\cdot\left[\mathbb{\tilde E}\sum_{k=1}^N\int_0^T\left(\big\|\varsigma_k(\tilde u^\varepsilon(s))\big\|_{L^2}^2\|\tilde u^\varepsilon(s)\|_{L^2}^2+\big\|\varsigma_k(u(s))\big\|_{L^2}^2\| u(s)\|_{L^2}^2\right){\rm d}s\right]^{\frac{1}{2}}\nonumber\\
	\lesssim &\|\varphi\|_{L^\infty}^2\left(\mathbb{\tilde E}\sum_{k=1}^N\int_0^T\big\|\psi_\varepsilon\ast\varsigma_k(\tilde u^\varepsilon(s))\cdot\tilde u^\varepsilon(s)-\varsigma_k(u(s))\cdot u(s)\big\|_{L^1}^2{\rm d}s\right)^{\frac{1}{2}}\nonumber\\
	&\quad\quad\quad\cdot\left[\mathbb{\tilde E}\int_0^T\left(1+\|\tilde u^\varepsilon(s)\|_{L^2}^2\right)^2{\rm d}s+\mathbb{\tilde E}\int_0^T\left(1+\| u(s)\|_{L^2}^2\right)^2{\rm d}s\right]^{\frac{1}{2}},\nonumber
\end{align}
where the term $\displaystyle\mathbb{\tilde E}\int_0^T\left(1+\|\tilde u^\varepsilon(s)\|_{L^2}^2\right)^2{\rm d}s$ is bounded uniformly in $\varepsilon$ by (\ref{standard uniform bounds 1b}). Therefore we only need to show
\begin{align}
	\lim_{\varepsilon\rightarrow0}\mathbb{\tilde E}\sum_{k=1}^N\int_0^T\big\|\psi_\varepsilon\ast\varsigma_k(\tilde u^\varepsilon(s))\cdot\tilde u^\varepsilon(s)-\varsigma_k(u(s))\cdot u(s)\big\|_{L^1}^2{\rm d}s=0.\nonumber
\end{align}
By H\"older's inequality and (\ref{con. linear growth}) we see that
\begin{align}
	&\mathbb{\tilde E}\sum_{k=1}^N\int_0^T\Big\|\underbrace{\psi_\varepsilon\ast\varsigma_k(\tilde u^\varepsilon(s))\cdot\tilde u^\varepsilon(s)-\varsigma_k(u(s))\cdot u(s)}_{=\psi_{\varepsilon}\ast\varsigma_k(\tilde u^{\varepsilon})\cdot(\tilde u^{\varepsilon}-u)+\big[\psi_{\varepsilon}\ast\varsigma_k(\tilde u^{\varepsilon})-\psi_{\varepsilon}\ast\varsigma_k(u)\big]\cdot u +\big[\psi_{\varepsilon}\ast\varsigma_k(u)-\varsigma_k(u)\big]\cdot u}\Big\|_{L^1}^2{\rm d}s\nonumber\\
	\lesssim &\mathbb{\tilde E}\int_0^T\left(1+\big\|\tilde u^{\varepsilon}(s)\big\|_{L^2}^2\right)\big\|\tilde u^{\varepsilon}(s)-u(s)\big\|_{L^2}^2{\rm d}s+\mathbb{\tilde E}\sum_{k=1}^N\int_0^T\big\|\varsigma_k(\tilde u^{\varepsilon}(s))-\varsigma_k(u(s))\big\|_{L^2}^2\|u(s)\|_{L^2}^2{\rm d}s\nonumber\\
	&+\mathbb{\tilde E}\sum_{k=1}^N\int_0^T\big\|\psi_{\varepsilon}\ast\varsigma_k(u(s))-\varsigma_k(u(s))\big\|_{L^2}^2\|u(s)\|_{L^2}^2{\rm d}s.		
\end{align}

{\noindent In the following we show the limits:}
\begin{align}
	&\lim_{\varepsilon\rightarrow0}\mathbb{\tilde E}\int_0^T\left(1+\big\|\tilde u^{\varepsilon}(s)\big\|_{L^2}^2\right)\big\|\tilde u^{\varepsilon}(s)-u(s)\big\|_{L^2}^2{\rm d}s=0~,\label{control of convergence of N(varphi) 2}\\
	&\lim_{\varepsilon\rightarrow0}\mathbb{\tilde E}\int_0^T\sum_{k=1}^\infty\big\|\varsigma_k(\tilde u^{\varepsilon}(s))-\varsigma_k(u(s))\big\|_{L^2}^2\big\|u(s)\big\|_{L^2}^2{\rm d}s=0~,\label{control of convergence of N(varphi)3}\\
	&\lim_{\varepsilon\rightarrow0}\mathbb{\tilde E}\int_0^T\sum_{k=1}^\infty\big\|\psi_{\varepsilon}\ast\varsigma_k(u(s))-\varsigma_k(u(s))\big\|_{L^2}^2\|u(s)\|_{L^2}^2{\rm d}s=0~.\label{control of convergence of N(varphi) 1}
\end{align}
(\ref{control of convergence of N(varphi) 1}) follows by Dominated Convergence Theorem, since:
\begin{align}
    &\sum_{k=1}^\infty\big\|\psi_{\varepsilon}\ast\varsigma_k(u(s))-\varsigma_k(u(s))\big\|_{L^2}^2\|u(s)\|_{L^2}^2\longrightarrow0~,~~{\rm d}\mathbb{\tilde P}\times{\rm d}t-a.e.;~(again~by~Dominated~Convergence)\nonumber\\
	&\sum_{k=1}^\infty\big\|\psi_{\varepsilon}\ast\varsigma_k(u(s))-\varsigma_k(u(s))\big\|_{L^2}^2\|u(s)\|_{L^2}^2\lesssim\big(1+\|u(s)\|_{L^2}^2\big)^2~\big(\in L^1(\tilde\Omega\times[0,T])\big)\quad (by~(\ref{con. linear growth})~and~(\ref{standard uniform bounds 1b})).\nonumber
\end{align}
(\ref{control of convergence of N(varphi)3}) can be handled similarly as (\ref{control of convergence of M(varphi) 2}). We split the infinite sum by an arbitrary $N'\in\mathbb N$, then for summation above $N'$ we see that 
\begin{align}
	&\mathbb{\tilde E}\int_0^T\sum_{k=N'+1}^\infty\big\|\varsigma_k(\tilde u^{\varepsilon}(s))-\varsigma_k(u(s))\big\|_{L^2}^2\big\|u(s)\big\|_{L^2}^2{\rm d}s\nonumber\\
	\lesssim &\left(\sup_{u\in L^2}\frac{\sum_{k>N'}\big\|\varsigma_k(u)\big\|_{L^2}^2}{1+\|u\|_{L^2}^2}\right)\mathbb{\tilde E}\int_0^T\left(2+\|u(s)\|_{L^2}^2+\big\|\tilde u^{\varepsilon}(s)\big\|_{L^2}^2\right)\big\|u(s)\big\|_{L^2}^2{\rm d}s\nonumber\\
	\lesssim &C_u\left(\sup_{u\in L^2}\frac{\sum_{k>N'}\big\|\varsigma_k(u)\big\|_{L^2}^2}{1+\|u\|_{L^2}^2}\right)\nonumber
\end{align}
where on the right-hand-side the expectation term is bounded uniformly in $\varepsilon$ by H\"older's inequality and (\ref{standard uniform bounds 1b}):
\begin{align}
	&\mathbb{\tilde E}\int_0^T\left(2+\|u(s)\|_{L^2}^2+\big\|\tilde u^{\varepsilon}(s)\big\|_{L^2}^2\right)\big\|u(s)\big\|_{L^2}^2{\rm d}s\nonumber\\
	\leq &\mathbb{\tilde E}\int_0^T\left(2+\|u(s)\|_{L^2}^2\right)\big\|u(s)\big\|_{L^2}^2{\rm d}s+C_1^{1/2}\left(\mathbb{\tilde E}\|u\|_{L^2}^4\right)^{1/2}=:C_u.\nonumber
\end{align}
On the other hand, the limit
\begin{align}
	\lim_{\varepsilon\rightarrow0}\mathbb{\tilde E}\int_0^T\sum_{k=1}^{N'}\big\|\varsigma_k(\tilde u^{\varepsilon}(s))-\varsigma_k(u(s))\big\|_{L^2}^2\big\|u(s)\big\|_{L^2}^2{\rm d}s=0\nonumber
\end{align}
follows by Dominated Convergence Theorem:
\begin{align}
    &\sum_{k=1}^{N'}\big\|\varsigma_k(\tilde u^{\varepsilon}(s))-\varsigma_k(u(s))\big\|_{L^2}^2\big\|u(s)\big\|_{L^2}^2\longrightarrow0,~~{\rm d}\mathbb{\tilde P}\times{\rm d}t-a.e.,\quad (by~continuity~of~\varsigma_k~'s);\nonumber\\
	&\sum_{k=1}^{N'}\big\|\varsigma_k(\tilde u^{\varepsilon}(s))-\varsigma_k(u(s))\big\|_{L^2}^2\big\|u(s)\big\|_{L^2}^2\lesssim\big(1+\|u(s)\|_{L^2}^2\big)^2~\big(\in L^1(\tilde\Omega\times[0,T])\big)\quad (by~(\ref{con. linear growth})~and~(\ref{standard uniform bounds 1b})).\nonumber
\end{align}
Hence, putting together we should have
\begin{align}
	\lim_{\varepsilon\rightarrow0}\mathbb{\tilde E}\int_0^T\sum_{k=1}^{\infty}\big\|\varsigma_k(\tilde u^{\varepsilon}(s))-\varsigma_k(u(s))\big\|_{L^2}^2\big\|u(s)\big\|_{L^2}^2{\rm d}s\lesssim C_u\left(\sup_{u\in L^2}\frac{\sum_{k>N'}\big\|\varsigma_k(u)\big\|_{L^2}^2}{1+\|u\|_{L^2}^2}\right).\nonumber
\end{align}
As left-hand-side is independent of $N'$, we could send $N'\rightarrow\infty$, and (\ref{control of convergence of N(varphi)3}) follows by (\ref{con. uniform decay of the tails}).\\ 
For (\ref{control of convergence of N(varphi) 2}), we first write by H\"older's inequality:
\begin{align}
	&\mathbb{\tilde E}\int_0^T\left(1+\big\|\tilde u^{\varepsilon}(s)\big\|_{L^2}^2\right)\big\|\tilde u^{\varepsilon}(s)-u(s)\big\|_{L^2}^2{\rm d}s\nonumber\\
	\lesssim &\mathbb{\tilde E}\int_0^T\left(1+\|u(s)\|_{L^2}^2+\big\|\tilde u^{\varepsilon}(s)\big\|_{L^2}^2\right)^{3/2}\big\|\tilde u^{\varepsilon}(s)-u(s)\big\|_{L^2}{\rm d}s\nonumber\\
	\lesssim &\left[\mathbb{\tilde E}\int_0^T\left(1+\|u(s)\|_{L^2}^2+\big\|\tilde u^{\varepsilon}(s)\big\|_{L^2}^2\right)^{3}{\rm d}s\right]^{1/2}\left(\mathbb{\tilde E}\int_0^T\big\|\tilde u^{\varepsilon}(s)-u(s)\big\|_{L^2}^2{\rm d}s\right)^{1/2}.\nonumber
\end{align}
Then by (\ref{standard uniform bounds 1b}), it suffices to show $\displaystyle\mathbb{\tilde E}\int_0^T\big\|\tilde u^{\varepsilon}(s)-u(s)\big\|_{L^2}^2{\rm d}s\longrightarrow0$. For this, we take arbitrary $\delta>0$ and write by H\"older's inequality:
\begin{align}
	\mathbb{\tilde E}\int_0^T\big\|\tilde u^{\varepsilon}(s)-u(s)\big\|_{L^2}^2{\rm d}s
	&=\int_{\big\{\|\tilde u^{\varepsilon}-u\|_{L^2}\leq\delta\big\}}\big\|\tilde u^{\varepsilon}-u\big\|_{L^2}^2{\rm d}\mathbb{\tilde P}+\int_{\big\{\|\tilde u^{\varepsilon}-u\|_{L^2}>\delta\big\}}\big\|\tilde u^{\varepsilon}-u\big\|_{L^2}^2{\rm d}\mathbb{\tilde P}\nonumber\\
	&\leq\delta^2+\mathbb{\tilde E}\left(1_{\big\{\|\tilde u^{\varepsilon}-u\|_{L^2}>\delta\big\}}\big\|\tilde u^{\varepsilon}-u\big\|_{L^2}^2\right)\nonumber\\
	&\leq\delta^2+\mathbb{\tilde P}\Big(\|\tilde u^{\varepsilon}-u\|_{L^2}>\delta\Big)^{1/2}\left(\mathbb{\tilde E}\big\|\tilde u^{\varepsilon}-u\big\|_{L^2}^4\right)^{1/2}.\nonumber
\end{align}
With (\ref{standard uniform bounds 1b}) and (\ref{strong convergences of u^varepsilon}) in hand, we could send $\varepsilon\rightarrow0$ (possible after subtracting a subsequence) and deduce
\begin{align}
	\lim_{\varepsilon\rightarrow 0}\mathbb{\tilde E}\int_0^T\big\|\tilde u^{\varepsilon}(s)-u(s)\big\|_{L^2}^2{\rm d}s\leq\delta^2,\quad \forall\delta>0.
\end{align}
This finally proves (\ref{control of convergence of N(varphi) 2}). 

\par Now take any bounded continuous function $\phi$ on $Z$ and by (\ref{quadratic variation of N^varepsilon}) we have for $0\leq t_1<t_2\leq T$ that
\begin{align}
	&\mathbb{\tilde E}\Big[\left(\tilde{\mathbf N}^{\varepsilon}_{t_2}(\varphi)-\tilde{\mathbf N}^{\varepsilon}_{t_1}(\varphi)\right)\phi\left(\tilde u^\varepsilon|_{[0,t_1]}\right)\Big]=0,\label{martingale property of N^variepsilon}\\
	&\mathbb{\tilde E}\bigg[\Big(\tilde{\mathbf N}^{\varepsilon}_{t_2}(\varphi)^2-\tilde{\mathbf N}^{\varepsilon}_{t_1}(\varphi)^2\Big)\phi\left(\tilde u^\varepsilon|_{[0,t_1]}\right)\bigg]\nonumber\\
	=&\mathbb{\tilde E}\Bigg[\left(4\sum_{k=1}^N\int_{t_1}^{t_2}\left|\int_{\mathbb{T}^3}\psi_\varepsilon\ast\varsigma_k(\tilde u^\varepsilon(s))\cdot\tilde u^\varepsilon(s)\varphi(s){\rm d}x\right|^2{\rm d}s\right)\phi\left(\tilde u^\varepsilon|_{[0,t_1]}\right)\Bigg].\label{quadratic variation equality of N^varepsilon}
\end{align}
Then by using (\ref{strong convergences of u^varepsilon}), (\ref{L^2-convergence of N^varepsilon}) and (\ref{convergence of quadratic variation of N^varepsilon}), one can take limits in (\ref{martingale property of N^variepsilon}) and (\ref{quadratic variation equality of N^varepsilon}) to see that $\mathbf N(\varphi)$ is a continuous square-integrable martingale (w.r.t. the filtration $\mathscr{\tilde F}_t=\sigma(u(s);0\leq s\leq t)$ ) with quadratic variation
\begin{align}
	\big\langle \mathbf N(\varphi)\big\rangle=4\sum_{k=1}^\infty\int_0^\cdot\left|\int_{\mathbb{T}^3}\varsigma_k(u(s))\cdot u(s)\varphi(s){\rm d}x\right|^2{\rm d}s.
\end{align}
\par To this end, we let $t=T$ in (\ref{LEE for u^varepsilon}), take $L^2(\tilde\Omega)$-inner product with any non-negative $\xi\in L^2(\tilde\Omega)$ on both sides and then send $\varepsilon\rightarrow0$, by (\ref{convergences of I1 I2 I3}) and (\ref{L^2-convergence of N^varepsilon}) we would have
\begin{align}\label{limiting energy equality}
	&\lim_{\varepsilon\rightarrow0}2\nu\mathbb{\tilde E}\,\xi\int_0^T\int_{\mathbb{T}^3}\big|\nabla\tilde u^\varepsilon(t)\big|^2\varphi(t){\rm d}x{\rm d}t\nonumber\\
	=&\mathbb{\tilde E}\,\xi\int_0^T\int_{\mathbb{T}^3}\big|u(t)\big|^2\Big(\partial_t\varphi(t)+\nu\Delta\varphi(t)\Big){\rm d}x{\rm d}t+\mathbb{\tilde E}\,\xi\int_0^T\int_{\mathbb{T}^3}\Big(u(t)^2+2p(t)\Big)u(t)\cdot\nabla\varphi(t){\rm d}x{\rm d}t\nonumber\\
	&+\mathbb{\tilde E}\,\xi\left(\sum_{k=1}^{\infty}\int_0^T\int_{\mathbb{T}^3}\big|\varsigma_k(u(t))\big|^2\varphi(t){\rm d}x{\rm d}t+\mathbf N_T(\varphi)\right).
\end{align}
On the other hand, by (\ref{weak convergences of u^varepsilon}), if we take the test function $\varphi$ to be non-negative, then $\sqrt{\xi\varphi}\nabla\tilde u^\varepsilon\rightharpoonup\sqrt{\xi\varphi}\nabla u$ in $L^2\big(\tilde\Omega\times[0,T]\times\mathbb T^3\big)$. By lower semicontinuity of norm, we have
\begin{align}
	2\nu\mathbb{\tilde E}\,\xi\int_0^T\int_{\mathbb{T}^3}\big|\nabla u(t)\big|^2\varphi(t){\rm d}x{\rm d}t\leq\lim_{\varepsilon\rightarrow0}2\nu\mathbb{\tilde E}\,\xi\int_0^T\int_{\mathbb{T}^3}\big|\nabla\tilde u^\varepsilon(t)\big|^2\varphi(t){\rm d}x{\rm d}t\nonumber,
\end{align}
which, together with (\ref{limiting energy equality}), proves (\ref{LEI-I}).

\end{proof}
\

\section{Proof of Theorem \ref{Thm. LEI. II} and Theorem \ref{Thm. energy supermartingale}}\label{Sec. LEIs of various types}
The derivations of inequality (\ref{LEI. II}) and Theorem \ref{Thm. energy supermartingale} are given by choosing suitable families of test function in (\ref{LEI-I}) and then taking limits. This requires certain continuity of the map $\varphi\longmapsto\mathbf{N}(\varphi)$ which clearly is linear.
\begin{lemma}\label{Lem. continuity of martingale map}
	Let $1\leq\alpha<4$. The martingale $\mathbf{N}(\varphi)$ (see definition \ref{def. martingale suitable weak solutions}) associated to the martingale suitable weak solution given in Theorem \ref{Thm. existence of martingale suitable solution} satisfies the following estimate
	\begin{align}\label{continuity of martingale map}
		\mathbb{\tilde E}\sup_{t\in[0,T]}\big|\mathbf{N}_t(\varphi_1)-\mathbf{N}_t(\varphi_2)\big|^\alpha\lesssim\big\|\varphi_1-\varphi_2\big\|_{L^5([0,T]\times\mathbb T^3)}^\alpha
	\end{align}
	for any non-negative $\varphi_1,\varphi_2\in C^\infty_c\big((0,T)\times\mathbb T^3\big)$.  
\end{lemma}
\begin{proof}[\bf Proof]
	Let $\varphi_1,\varphi_2\in C^\infty_c\big((0,T)\times\mathbb T^3\big)$ be given as in the lemma. Note that by (\ref{L^2-convergence of N^varepsilon}), it suffices to show (\ref{continuity of martingale map}) for $\tilde{\mathbf N}^\varepsilon(\varphi)$. We recall its quadratic variation given by (\ref{quadratic variation of tilde N^varepsilon}). Then for any $\alpha\geq1$ to be chosen later, we estimate by H\"older's inequality for $2<q<6$, by Sobolev imbedding $H^{\frac{3}{2}-\frac{3}{q}}\hookrightarrow L^q$ and interpolation $H^{\frac{3}{2}-\frac{3}{q}}=\big[L^2,H^1\big]_{\frac{3}{q}-\frac{1}{2}}$: 
	\begin{align}\label{martingale map estimate}
		&\ \ \ \ \ \ \mathbb{\tilde E}\sup_{t\in[0,T]}\big|\mathbf{N}^\varepsilon_t(\varphi_1)-\mathbf{N}^\varepsilon_t(\varphi_2)\big|^\alpha\nonumber\\
		&\leq2^\alpha~\mathbb{\tilde E}\left[\sum_{k=1}^N\int_0^T\left(\int_{\mathbb{T}^3}\big|\psi_\varepsilon\ast\sigma_k(\tilde u^\varepsilon(s))\big|\big|\tilde u^\varepsilon(s)\big|\big|\varphi_1(s)-\varphi_2(s)\big|{\rm d}x\right)^2{\rm d}s\right]^{\alpha/2}\nonumber\\
		&\leq2^\alpha~\mathbb{\tilde E}\left(\int_0^T\sum_{k=1}^N\big\|\psi_\varepsilon\ast\sigma_k(\tilde u^\varepsilon(s))\big\|_{L^2}^2\big\|\tilde u^\varepsilon(s)\big\|_{L^q}^2\big\|\varphi_1(s)-\varphi_2(s)\big\|_{L^{\frac{2q}{q-2}}}^2{\rm d}s\right)^{\alpha/2}\nonumber\\
		&\lesssim2^\alpha~\mathbb{\tilde E}\left(\int_0^T\left(1+\big\|\tilde u^\varepsilon(s)\big\|_{L^2}^2\right)\big\|\tilde u^\varepsilon(s)\big\|_{L^2}^{\frac{6}{q}-1}\big\|\tilde u^\varepsilon(s)\big\|_{H^1}^{3(1-\frac{2}{q})}\big\|\varphi_1(s)-\varphi_2(s)\big\|_{L^{\frac{2q}{q-2}}}^2{\rm d}s\right)^{\alpha/2}\nonumber\\
		&\lesssim2^\alpha~\mathbb{\tilde E}\left[\left(1+\big\|\tilde u^\varepsilon\big\|_{L^\infty([0,T];L^2)}^2\right)^{\frac{3}{q}+\frac{1}{2}}\big\|\tilde u^\varepsilon\big\|_{L^2([0,T];H^1)}^{3(1-\frac{2}{q})}\right]^{\alpha/2}\big\|\varphi_1-\varphi_2\big\|_{L^{\frac{4q}{6-q}}\big([0,T];L^{\frac{2q}{q-2}}\big)}^\alpha.
	\end{align}
	Note that by $q<6$, we have $3(1-\frac{2}{q})<2$. Hence, choosing $\alpha\geq1$ such that $\frac{3\alpha}{2}(1-\frac{2}{q})<2$, then by H\"older's inequality and (\ref{standard uniform bounds 1b}), the quantity
	\begin{align}
		\mathbb{\tilde E}\left[\left(1+\big\|\tilde u^\varepsilon\big\|_{L^\infty([0,T];L^2)}^2\right)^{\frac{3}{q}+\frac{1}{2}}\big\|\tilde u^\varepsilon\big\|_{L^2([0,T];H^1)}^{3(1-\frac{2}{q})}\right]^{\alpha/2}\nonumber
	\end{align}
	is bounded uniformly in $\varepsilon$. Now choose $q=10/3$ in (\ref{martingale map estimate}), then we get (\ref{continuity of martingale map}) for any $1\leq\alpha<10/3$.
\end{proof}
\

\par We note here that by taking non-negative test functions in $C^\infty_c\big((a,b)\times\mathbb T^3\big)$ in (\ref{LEI-I}) for arbitrarily fixed $(a,b)\subset[0,T]$, the inequality becomes
\begin{align}\label{LEI on (a,b)}
	&2\nu\mathbb E\,\xi\int_a^b\int_{\mathbb{T}^3}\big|\nabla u(\tau)\big|^2\varphi(\tau){\rm d}x{\rm d}\tau
	\leq\mathbb E\,\xi\left(\sum_{k=1}^{\infty}\int_a^b\int_{\mathbb{T}^3}\big|\varsigma_k(u(\tau))\big|^2\varphi(\tau){\rm d}x{\rm d}\tau+\mathbf N_b(\varphi)\right)+\nonumber\\
	+&\mathbb E\,\xi\int_a^b\int_{\mathbb{T}^3}\big|u(\tau)\big|^2\Big(\partial_\tau\varphi(\tau)+\nu\Delta\varphi(\tau)\Big){\rm d}x{\rm d}\tau+\mathbb E\,\xi\int_a^b\int_{\mathbb{T}^3}\Big(u(\tau)^2+2p(\tau)\Big)u(\tau)\cdot\nabla\varphi(\tau){\rm d}x{\rm d}\tau
\end{align}
And here, the martingale $\mathbf{N}(\varphi)$ satisfies
\begin{align}
	\mathbf{N}_\tau(\varphi)&=0\ \ when\ \tau\in[0,a],\nonumber\\
	\mathbf{N}_\tau(\varphi)&=\mathbf{N}_b(\varphi)\ \ when\ \tau\in[b,T].\nonumber
\end{align}
Moreover, if $\varphi_1$ and $\varphi_2$ agree on some $[c,d]\subset(a,b)$, then
\begin{align}
	\mathbf{N}_d(\varphi_1)-\mathbf{N}_c(\varphi_1)=\mathbf{N}_d(\varphi_2)-\mathbf{N}_c(\varphi_2).\nonumber
\end{align}
Those equalities follow directly for $\mathbf{N}^\varepsilon(\varphi)$ and hence for $\mathbf{N}(\varphi)$ by (\ref{equal Laws}) and (\ref{L^2-convergence of M^varepsilon}). In the following subsections, we use (\ref{LEI on (a,b)}) and Lemma \ref{Lem. continuity of martingale map} to prove Theorem \ref{Thm. LEI. II} and \ref{Thm. energy supermartingale}.

\
\subsection{Local energy inequality (\ref{LEI. II})}\label{Subsec. Local energy inequality of type II}
\begin{proof}[\bf Proof of Theorem \ref{Thm. LEI. II}]
Let $(a,b)\subset[0,T]$ and $a<t<b$ be arbitrarily fixed. We set a family of non-negative test functions $\big\{\varphi_\ell\big\}\subset C^\infty_c\big((a,b)\times\mathbb T^3\big)$ as follows. Let $\chi\in C^\infty(\mathbb R)$ be such that
\begin{equation}\label{def. function chi}
	\left\{\begin{aligned}
		\chi(r)=0,\ \ \ r\leq0,\\
		\chi(r)=1,\ \ \ r\geq1,
	\end{aligned}\right.\ \ \ and\ \ \ 0\leq\chi(r)\leq1\ for\ all\ r\in\mathbb R.
\end{equation}
Then for all sufficiently small $\ell>0$, we have
\begin{equation}
	\left\{\begin{aligned}\label{values of chi_ell}
		&\chi\left(\frac{t-\tau}{\ell}\right)=1,\ \ \ \tau\in[a,t-\ell],\\
		&\chi\left(\frac{t-\tau}{\ell}\right)=0,\ \ \ \tau\in[t,b].
	\end{aligned}\right.\ \
\end{equation}
Now taking any non-negative $\varphi\in C^\infty_c\big((a,b)\times\mathbb T^3\big)$, we set the family of test functions
\begin{align}
	\varphi_\ell(\tau,x):=\chi\left(\frac{t-\tau}{\ell}\right)\varphi(\tau,x).\nonumber
\end{align}
Put $\varphi_\ell$ in (\ref{LEI on (a,b)}) and take limit $\ell\rightarrow0$. Then by (\ref{values of chi_ell}) and Dominated Convergence, we have
\begin{equation}\label{Convergences of bounded terms sec.3.1}
	\left\{\begin{aligned}
		&\int_a^b\int_{\mathbb{T}^3}\big|\nabla u(\tau)\big|^2\varphi_\ell(\tau){\rm d}x{\rm d}\tau\stackrel{L^1(\tilde\Omega)}{\longrightarrow}\int_a^t\int_{\mathbb{T}^3}\big|\nabla u(\tau)\big|^2\varphi(\tau){\rm d}x{\rm d}\tau,\\
		&\sum_{k=1}^{\infty}\int_a^b\int_{\mathbb{T}^3}\big|\varsigma_k(u(\tau))\big|^2\varphi_\ell(\tau){\rm d}x{\rm d}\tau\stackrel{L^1(\tilde\Omega)}{\longrightarrow}
		\sum_{k=1}^{\infty}\int_a^b\int_{\mathbb{T}^3}\big|\varsigma_k(u(\tau))\big|^2\varphi(\tau){\rm d}x{\rm d}\tau,\\
		&\int_a^b\int_{\mathbb{T}^3}\big|u(\tau)\big|^2\nu\Delta\varphi_\ell(\tau){\rm d}x{\rm d}\tau\stackrel{L^1(\tilde\Omega)}{\longrightarrow}\int_a^t\int_{\mathbb{T}^3}\big|u(\tau)\big|^2\nu\Delta\varphi(\tau){\rm d}x{\rm d}\tau,\\
		&\int_a^b\int_{\mathbb{T}^3}\Big(u(\tau)^2+2p(\tau)\Big)u(\tau)\cdot\nabla\varphi_\ell(\tau){\rm d}x{\rm d}\tau\stackrel{L^1(\tilde\Omega)}{\longrightarrow}\int_a^t\int_{\mathbb{T}^3}\Big(u(\tau)^2+2p(\tau)\Big)u(\tau)\cdot\nabla\varphi(\tau){\rm d}x{\rm d}\tau.
	\end{aligned}\right.
\end{equation}
For the term involving $\partial_\tau\varphi_\ell$, we have
\begin{align}
	&\int_a^b\int_{\mathbb{T}^3}\big|u(\tau)\big|^2\partial_\tau\varphi_\ell(\tau){\rm d}x{\rm d}\tau\nonumber\\
	=&\int_a^b\int_{\mathbb{T}^3}\big|u(\tau)\big|^2\chi\left(\frac{t-\tau}{\ell}\right)\partial_\tau\varphi(\tau){\rm d}x{\rm d}\tau
	-\int_a^b\ell^{-1}\chi'\left(\frac{t-\tau}{\ell}\right)\int_{\mathbb{T}^3}\big|u(\tau)\big|^2\partial_\tau\varphi(\tau){\rm d}x{\rm d}\tau.\nonumber
\end{align}
Again by (\ref{values of chi_ell}) and Dominated Convergence, we see that
\begin{align}\label{convergence of the other bounded term sec.3.1}
	\int_a^b\int_{\mathbb{T}^3}\big|u(\tau)\big|^2\chi\left(\frac{t-\tau}{\ell}\right)\partial_\tau\varphi(\tau){\rm d}x{\rm d}\tau\stackrel{L^1(\tilde\Omega)}{\longrightarrow}\int_a^t\int_{\mathbb{T}^3}\big|u(\tau)\big|^2\partial_\tau\varphi(\tau){\rm d}x{\rm d}\tau.
\end{align}
For the other term, it's clear that $\displaystyle\int_a^b\ell^{-1}\chi'\left(\frac{t-\tau}{\ell}\right){\rm d}\tau=1$ and $\displaystyle\ell^{-1}\chi'\left(\frac{\cdot}{\ell}\right)$ is an approximate identity. Denote $\displaystyle\gamma_\ell(\cdot):=\ell^{-1}\chi'\left(\frac{\cdot}{\ell}\right)$ and $\displaystyle{\rm U}(\cdot):=\int_{\mathbb{T}^3}\big|u(\cdot,x)\big|^2\varphi(\cdot,x){\rm d}x$. Since ${\rm U}\in L^1[a,b]$, $\mathbb{\tilde P}$-$a.s.$, we have
\begin{align}
	\int_a^b\big|\gamma_\ell\ast{\rm U}(\tau)-{\rm U}(\tau)\big|{\rm d}\tau\stackrel{a.s.}{\longrightarrow}0,\nonumber
\end{align}
and
\begin{align}
	\int_a^b\big|\gamma_\ell\ast{\rm U}(\tau)-{\rm U}(\tau)\big|{\rm d}\tau\leq2\int_a^b\big|{\rm U}(\tau)\big|{\rm d}\tau\in L^1(\tilde\Omega).\nonumber
\end{align}
Then, by Fubini's Theorem and Dominated Convergence, we have
\begin{align}
	\lim_{\ell\rightarrow0}\int_a^b\mathbb{\tilde E}\big|\gamma_\ell\ast{\rm U}(\tau)-{\rm U}(\tau)\big|{\rm d}\tau=\lim_{\ell\rightarrow0}\mathbb{\tilde E}\int_a^b\big|\gamma_\ell\ast{\rm U}(\tau)-{\rm U}(\tau)\big|{\rm d}\tau=0.\nonumber
\end{align}
Hence, after possibly subtracting a subsequence, we have
\begin{align}
	\lim_{\ell\rightarrow0}\mathbb{\tilde E}\big|\gamma_\ell\ast{\rm U}(\tau)-{\rm U}(\tau)\big|=0\ \ for\ a.e.\ \tau\in[a,b].\nonumber
\end{align}
This gives that 
\begin{align}\label{a.e. convergence term}
	\int_a^b\int_{\mathbb{T}^3}\big|u(\tau)\big|^2\partial_\tau\varphi_\ell(\tau){\rm d}x
	&{\rm d}\tau=\int_a^b\int_{\mathbb{T}^3}\big|u(\tau)\big|^2\chi\left(\frac{t-\tau}{\ell}\right)\partial_\tau\varphi(\tau){\rm d}x{\rm d}\tau-\gamma_\ell\ast{\rm U}(t)\nonumber\\
	&\stackrel{L^1(\tilde\Omega)}{\longrightarrow}\int_a^b\int_{\mathbb{T}^3}\big|u(\tau)\big|^2\partial_\tau\varphi(\tau){\rm d}x{\rm d}\tau-\int_{\mathbb{T}^3}\big|u(t,x)\big|^2\varphi(t,x){\rm d}x\ for\ a.e.\ t\in[a,b].
\end{align}
For the martingale term, we have
\begin{align}
	\mathbf{N}_b(\varphi_\ell)=\mathbf{N}_t(\varphi_\ell).
\end{align}
And by (\ref{continuity of martingale map}), we see that
\begin{align}\label{convergence of martingale term}
	\mathbb{\tilde E}\big|\mathbf{N}_t(\varphi_\ell)-\mathbf{N}_t(\varphi)\big|\lesssim\big\|\varphi_\ell-\varphi\big\|_{L^5([0,t]\times\mathbb T^3)}\leq2\big\|\varphi\big\|_{L^5([t-\ell,t]\times\mathbb T^3)}\longrightarrow0.
\end{align}
Hence, together by (\ref{Convergences of bounded terms sec.3.1}), (\ref{convergence of the other bounded term sec.3.1}), (\ref{a.e. convergence term}) and (\ref{convergence of martingale term}), we see that (\ref{LEI. II}) holds for $a.e.$ $t\in[a,b]$. For arbitrary $t\in[0,T]$, one could set a sequence $\{t_n\}$ of times for which (\ref{LEI. II}) holds such that $t_n\rightarrow t$. And it's straght forward to apply Fatou's lemma to see that (\ref{LEI. II}) holds for $t$. This ends the proof.
\end{proof}

\
\subsection{Super-martingale statement of local energy inequality}\label{Subsec. Super-martingale statement of local energy inequality}
\begin{proof}[\bf Proof of Theorem \ref{Thm. energy supermartingale}]
The proof is quite similar to the previous one. Let $(a,b)\subset[0,T]$ and $a<t<b$ be arbitrarily fixed so that (\ref{LEI. II}) holds. Fixed also an $s\in[a,t)$. We set a family of non-negative test functions $\big\{\phi_\ell\big\}\subset C^\infty_c\big((a,b)\times\mathbb T^3\big)$ by
\begin{align}
	\phi_\ell(\tau,x):=\chi\left(\frac{\tau-s}{\ell}\right)\varphi(\tau,x),\nonumber
\end{align} 
where $\varphi\in C^\infty_c\big((a,b)\times\mathbb T^3;\mathbb R_+\big)$ is arbitrarily given. Now put $\phi_\ell$ in (\ref{LEI. II}) and take limit $\ell\rightarrow0$. Then arguing as in subsection \ref{Subsec. Local energy inequality of type II}, we have the following convergences
\begin{equation}
	\left\{\begin{aligned}\label{Convergences of bounded terms sec.3.2}
	&\int_a^t\int_{\mathbb{T}^3}\big|\nabla u(\tau)\big|^2\phi_\ell(\tau){\rm d}x{\rm d}\tau\stackrel{L^1(\tilde\Omega)}{\longrightarrow}\int_s^t\int_{\mathbb{T}^3}\big|\nabla u(\tau)\big|^2\varphi(\tau){\rm d}x{\rm d}\tau,\\
	&\sum_{k=1}^{\infty}\int_a^b\int_{\mathbb{T}^3}\big|\varsigma_k(u(\tau))\big|^2\phi_\ell(\tau){\rm d}x{\rm d}\tau\stackrel{L^1(\tilde\Omega)}{\longrightarrow}
		\sum_{k=1}^{\infty}\int_a^b\int_{\mathbb{T}^3}\big|\varsigma_k(u(\tau))\big|^2\varphi(\tau){\rm d}x{\rm d}\tau,\\
	&\int_a^t\int_{\mathbb{T}^3}\big|u(\tau)\big|^2\Big(\nu\Delta\phi_\ell(\tau)+\kappa^2\phi_\ell(\tau)\Big){\rm d}x{\rm d}\tau\stackrel{L^1(\tilde\Omega)}{\longrightarrow}\int_s^t\int_{\mathbb{T}^3}\big|u(\tau)\big|^2\Big(\nu\Delta\varphi(\tau)+\kappa^2\varphi(\tau)\Big){\rm d}x{\rm d}\tau,\\
	&\int_a^t\int_{\mathbb{T}^3}\Big(u(\tau)^2+2p(\tau)\Big)u(\tau)\cdot\nabla\phi_\ell(\tau){\rm d}x{\rm d}\tau\stackrel{L^1(\tilde\Omega)}{\longrightarrow}\int_s^t\int_{\mathbb{T}^3}\Big(u(\tau)^2+2p(\tau)\Big)u(\tau)\cdot\nabla\varphi(\tau){\rm d}x{\rm d}\tau.	\end{aligned}\right.
\end{equation}
and
\begin{align}\label{a.e. convergence term sec.3.2}
	\int_a^t\int_{\mathbb{T}^3}\big|u(\tau)\big|^2\partial_\tau\phi_\ell(\tau){\rm d}x
		&\stackrel{L^1(\tilde\Omega)}{\longrightarrow}\int_s^t\int_{\mathbb{T}^3}\big|u(\tau)\big|^2\partial_\tau\varphi(\tau){\rm d}x{\rm d}\tau+\int_{\mathbb{T}^3}\big|u(s,x)\big|^2\varphi(s,x){\rm d}x\ for\ a.e.\ s\in[a,t].
\end{align}
Hence, it remains to deal with the martingale term. We first recall the definition (\ref{def. function chi}) of function $\chi$. And for this time, we have
\begin{equation}
	\left\{\begin{aligned}\label{values of chi_ell}
		&\chi\left(\frac{\tau-s}{\ell}\right)=0,\ \ \ \tau\in[a,s],\\
		&\chi\left(\frac{\tau-s}{\ell}\right)=1,\ \ \ \tau\in[s+\ell,t].
	\end{aligned}\right.\ \ \nonumber
\end{equation}
Then, since $\phi_\ell=\varphi$ on $[s+\ell,t]$, by the properties of $\mathbf{N}(\varphi)$, we have now
\begin{align}
	\mathbf{N}_t(\phi_\ell)
	&=\mathbf{N}_t(\phi_\ell)-\mathbf{N}_{s+\ell}(\phi_\ell)+\mathbf{N}_{s+\ell}(\phi_\ell)\nonumber\\
	&=\mathbf{N}_t(\varphi)-\mathbf{N}_{s+\ell}(\varphi)+\mathbf{N}_{s+\ell}(\phi_\ell).\nonumber
\end{align}
By the continuity of $\tau\mapsto\mathbf{N}_\tau(\varphi)$, $\mathbf{N}_{s+\ell}(\varphi)\longrightarrow\mathbf{N}_{s}(\varphi)$ strongly in $L^2(\tilde\Omega)$. And by (\ref{continuity of martingale map}), we have
\begin{align}
	\mathbb{\tilde E}\big|\mathbf{N}_{s+\ell}(\phi_\ell)\big|\lesssim\big\|\phi_\ell\big\|_{L^5([0,s+\ell]\times\mathbb T^3)}\leq\big\|\varphi\big\|_{L^5([s,s+\ell]\times\mathbb T^3)}\longrightarrow0.\nonumber
\end{align}
Hence,
\begin{align}\label{convergence of martingale term sec.3.2}
	\mathbf{N}_t(\phi_\ell)\stackrel{L^1(\tilde\Omega)}{\longrightarrow}\mathbf{N}_t(\varphi)-\mathbf{N}_{s}(\varphi).
\end{align}
Now, for any non-negative test function $\varphi\in C^\infty_c\left((a,b)\times\mathbb{T}^3\right)$ and non-negative random variable $\xi\in L^\infty(\tilde\Omega)$, combining (\ref{LEI. II}), (\ref{Convergences of bounded terms sec.3.2}), (\ref{a.e. convergence term sec.3.2}) and (\ref{convergence of martingale term sec.3.2}), we have that
\begin{align}\label{a.e. LEI sec.3.2}
	\mathbb{\tilde E}\,\xi\Bigg(\mathcal E_t(u;\varphi)-\mathcal E_s(u;\varphi)-\sum_{k=1}^{\infty}\int_s^t\int_{\mathbb{T}^3}\big|\varsigma_k(u(\tau))\big|^2\varphi(\tau){\rm d}x{\rm d}\tau-\mathbf N_t(\varphi)+\mathbf N_s(\varphi)\Bigg)\leq0
\end{align}
for all $s\in[a,t]$ except on a set $\mathcal T_t\subset[a,b]$ of Lebesgue measure zero. Here, we use the notation
\begin{align}
	\mathcal E_t(u;\varphi):
	&=\int_{\mathbb{T}^3}\big|u(t)\big|^2\varphi(t){\rm d}x+2\nu\int_0^t\int_{\mathbb{T}^3}\big|\nabla u(\tau)\big|^2\varphi(\tau){\rm d}x{\rm d}\tau-\int_0^t\int_{\mathbb{T}^3}\big|u(\tau)\big|^2\Big(\partial_\tau\varphi(\tau)+\nu\Delta\varphi(\tau)\Big){\rm d}x{\rm d}\tau\nonumber\\
	&-\int_0^t\int_{\mathbb{T}^3}\Big(u(\tau)^2+2p(\tau)\Big)u(\tau)\cdot\nabla\varphi(\tau){\rm d}x{\rm d}\tau~,\quad t\in[0,T].\nonumber
\end{align}
If we set $\xi=\mathbf1_{A}$ for arbitrary $A\in\mathscr{\tilde F}_s$ in (\ref{a.e. LEI sec.3.2}), we would see that the process $\mathcal E(\varphi)$ is an $a.e.$ $\{\mathscr{\tilde F}_t\}$-super-martingale. (The conditions (i) and (ii) of remark \ref{Remark a.e. sup-martingale} hold obviously.)
\end{proof}

\
\section{Derivation of Local energy equality}
We should note here that one could immediately obtain an ``energy equality" by defining $-2\mathcal D_t(u;\varphi)$ to be the stochastic process given by (\ref{energy a.e. supermartingale}). But then we lose the structure-function-like formulation (\ref{dissipation distribution D(u)}) for $\mathcal D(u;\varphi)$. A proper way is to mollifier the martingale formulation of the equation, i.e. (\ref{weak solu.}), apply It\^o's formula and take care of the limits. In the following, we denote $e_n(x)=e^{2\pi in\cdot x}$ ($n\in\mathbb Z^3$) and
\begin{align}
	\mathbf e_n^{1}=\left(e_n,0,0\right)^T,~\mathbf e_n^{2}=\left(0,e_n,0\right)^T,~\mathbf e_n^{3}=\left(0,0,e_n\right)^T.\nonumber
\end{align}
\begin{proof}[Proof of Theorem \ref{Thm. LEE}]
Let $\Big\{\left(\Omega, \mathscr{F}, \{\mathscr{F}_t\}_{t\in[0,T]}, \mathbb{P}\right),u\Big\}$ be a martingale suitable weak solution yielded by Theorem \ref{Thm. existence of martingale suitable solution} and $\mathbf M(\varphi)$ is given by (\ref{weak solu.}). Note that each $\mathbf M\big(\mathbf e_n^{i}\big)$ is a one-dimensional martingale. By slightly abusing the notation, we still denote by $\mathbf M(e_n)$ the three-dimensional martingale
\begin{align}
	\Big(\mathbf M\big(\mathbf e_n^{1}\big),\mathbf M\big(\mathbf e_n^{2}\big),\mathbf M\big(\mathbf e_n^{3}\big)\Big)^T.\nonumber
\end{align}
Cross-variation of two components can be calculated:
\begin{align}\label{cross-variation of two components of M(e)}
	\Big\langle\mathbf M(e_n)^i,\mathbf M(e_m)^j\Big\rangle_t=\Big\langle\mathbf M\big(\mathbf e_n^{i}\big),\mathbf M\big(\mathbf e_m^{j}\big)\Big\rangle_t
	&=\frac{\Big\langle\mathbf M\big(\mathbf e_n^{i}+\mathbf e_m^{j}\big)\Big\rangle_t-\Big\langle\mathbf M\big(\mathbf e_n^{i}\big)\Big\rangle_t-\Big\langle\mathbf M\big(\mathbf e_m^{j}\big)\Big\rangle_t}{2}\nonumber\\
	&=\sum_{k=1}^\infty\int_0^t\big\langle\varsigma_k(u(s)),\mathbf e_n^{i}\big\rangle \big\langle\varsigma_k(u(s)),\mathbf e_m^{j}\big\rangle {\rm d}\tau\nonumber\\
	&=\sum_{k=1}^\infty\int_0^t\big\langle\varsigma_k(u(s))^i, e_n\big\rangle \big\langle\varsigma_k(u(s))^j, e_m\big\rangle {\rm d}\tau
\end{align} 
Now, from (\ref{weak solu.}) we have
\begin{align}
	{\widehat u}_n(t)=\widehat u_{0,n}-\int_0^t\left(2\pi in\cdot\widehat{\big(u\otimes u\big)}_n(\tau)+4\pi^2\nu|n|^2{\widehat u}_n(\tau)+2\pi in{\widehat p}_n(\tau)\right){\rm d}\tau+\mathbf M_t(e_n).\nonumber
\end{align}
Applying It\^o's formula to the process ``${\widehat\alpha_\ell}(m){\widehat u}_m(t){\widehat u}_{n-m}(t){\widehat\varphi}_n(t)$" ( ${\widehat\alpha_\ell}(m)$ being Fourier coefficients of the mollifier $\alpha_\ell$ and $\varphi\in C^\infty_c\big((0,T)\times\mathbb T^3\big)$ ) and then summing in $m,n\in\mathbb Z^3$, we should have
\begin{align}\label{mollified LEE 1}
	&\int_{\mathbb T^3}u_{\ell}(t)\cdot u(t)\varphi(t){\rm d}x-\int_{\mathbb T^3}u_{0,\ell}\cdot u_0\varphi(0){\rm d}x\nonumber\\
	=&\int_0^t\int_{\mathbb T^3}\varphi u\cdot\Big(-{\rm div}(u\otimes u)_\ell+\nu\Delta u_\ell-\nabla p_\ell\Big){\rm d}x{\rm d}\tau+\int_0^t\int_{\mathbb T^3}\varphi u_\ell\cdot\Big(-{\rm div}(u\otimes u)+\nu\Delta u-\nabla p\Big){\rm d}x{\rm d}\tau\nonumber\\
	&+\sum_{k=1}^\infty\int_0^t\int_{\mathbb T^3}\alpha_\ell\ast\varsigma_k(u(\tau))\cdot\varsigma_k(u(\tau))\varphi(\tau){\rm d}x{\rm d}\tau+\mathcal N^\ell_t(\varphi)
\end{align}
where $u_\ell:=\alpha_\ell\ast u$, $u_{0,\ell}:=\alpha_\ell\ast u_0$, $(u\otimes u)_\ell:=\alpha_\ell\ast(u\otimes u)$ and
\begin{align}\label{martingale term of mollified LEE}
	\mathcal N^\ell_t(\varphi)
	&:=\sum_{m,n\in\mathbb Z^3}{\widehat\alpha_\ell}(m)\int_0^t{\widehat\varphi}_n(\tau){\widehat u}_{n-m}(\tau)\cdot{\rm d}\mathbf M_\tau(e_m)+\sum_{m,n\in\mathbb Z^3}{\widehat\alpha_\ell}(m)\int_0^t{\widehat\varphi}_n(\tau){\widehat u}_m(\tau)\cdot{\rm d}\mathbf M_\tau(e_{n-m})\nonumber\\
	&=2\sum_{m,n\in\mathbb Z^3}{\widehat\alpha_\ell}(m)\int_0^t{\widehat\varphi}_n(\tau){\widehat u}_{n-m}(\tau)\cdot{\rm d}\mathbf M_\tau(e_m)
\end{align}
Note that it is not hard to show that $\mathcal N^\ell(\varphi)$ is a continuous and square-integrable $\big\{\mathscr F_t\big\}$-martingale. And by (\ref{cross-variation of two components of M(e)}), its quadratic variation can be calculated:
\begin{align}\label{quadratic variation of mathcal-N^ell}
	\big\langle\mathcal N^\ell(\varphi)\big\rangle_t
	&=4\sum_{i,i'=1}^3\sum_{n,m\in\mathbb Z^3}\sum_{n',m'\in\mathbb Z^3}\sum_{k=1}^\infty\int_0^t{\widehat\alpha_\ell}(m){\widehat\alpha_\ell}(m'){\widehat\varphi}_n{\widehat u}_{n-m}^{~i}{\widehat\varphi}_{n'}{\widehat u}_{n'-m'}^{~i'}\big\langle\varsigma_k^i, e_n\big\rangle\big\langle\varsigma_k^{i'}, e_{n'}\big\rangle{\rm d}\tau\nonumber\\
	&\stackrel{*}{=}4\sum_{k=1}^\infty\int_0^t\left(\sum_{i=1}^3\sum_{n,m\in\mathbb Z^3}{\widehat\alpha_\ell}(m){\widehat\varphi}_n{\widehat u}_{n-m}^{~i}\big\langle\varsigma_k(u(s))^i, e_n\big\rangle\right)^2{\rm d}\tau\nonumber\\
	&=4\sum_{k=1}^\infty\int_0^t\left(\int_{\mathbb T^3}\alpha_\ell\ast\varsigma_k(u(\tau))\cdot u(\tau)\varphi(\tau){\rm d}x\right)^2{\rm d}\tau
\end{align} 
where on the equality with ``$*$" the change of orders is justified by Fubini's Theorem. We come back on the energy balance. With the following elementary calculation:
\begin{align}
	&{\bf a)}~u\cdot{\rm div}(u\otimes u)_\ell+u_\ell\cdot{\rm div}(u\otimes u)=u^j\partial_i(u^iu^j)_\ell+u^j_\ell u^i\partial_iu^j=\underbrace{u^j\partial_i(u^iu^j)_\ell- u^iu^j\partial_iu^j_\ell}_{=:{\bf E}^\ell}+{\rm div}\Big[(u_\ell\cdot u)u\Big];\nonumber\\
	&{\bf b)}~u\cdot\Delta u_\ell+u_\ell\cdot\Delta u=\Delta(u_\ell\cdot u)-2\nabla u_\ell:\nabla u~;\nonumber\\
	&{\bf c)}~u\cdot\nabla p_\ell+u_\ell\cdot\nabla p={\rm div}(u_\ell p+u~p_\ell)~,\nonumber
\end{align}
we could rewrite (\ref{mollified LEE 1}) as
\begin{align}\label{mollified LEE}
	&\int_{\mathbb T^3}u_{\ell}(t)\cdot u(t)\varphi(t){\rm d}x-\int_{\mathbb T^3}u_{0,\ell}\cdot u_0\varphi(0){\rm d}x+2\nu\int_0^t\int_{\mathbb T^3}\varphi\nabla u_\ell :\nabla u{\rm d}x{\rm d}\tau+\int_0^t\int_{\mathbb T^3}\varphi{\bf E}^\ell{\rm d}x{\rm d}\tau\nonumber\\
	=&\underbrace{\int_0^t\int_{\mathbb T^3}(u_\ell\cdot u)\big(\partial_\tau\varphi+\nu\Delta\varphi\big){\rm d}x{\rm d}\tau+\int_0^t\int_{\mathbb T^3}\Big((u_\ell\cdot u)u+u_\ell p+u~p_\ell\Big)\cdot\nabla\varphi{\rm d}x{\rm d}\tau}_{=:\mathcal I^\ell(t)}\nonumber\\
	&+\sum_{k=1}^\infty\int_0^t\int_{\mathbb T^3}\alpha_\ell\ast\varsigma_k(u(\tau))\cdot\varsigma_k(u(\tau))\varphi(\tau){\rm d}x{\rm d}\tau+\mathcal N^\ell_t(\varphi)~.
\end{align} 
Verification of convergence of each term above is quite standard. For the first term on the left-hand-side, we write by H\"older's inequality :
\begin{align}
	&\int_0^T\mathbb E\left|\int_{\mathbb T^3}u_{\ell}(t)\cdot u(t)\varphi(t){\rm d}x-\int_{\mathbb T^3}|u(t)|^2\varphi(t){\rm d}x\right|{\rm d}t\nonumber\\
	\leq &\|\varphi\|_{L^\infty}\mathbb E\int_0^T\int_{\mathbb T^3}\big|u_\ell(t)-u(t)\big||u(t)|{\rm d}x{\rm d}t\nonumber\\
	\leq &\|\varphi\|_{L^\infty}\left(\mathbb E\|u\|_{L^2}^2\right)^{1/2}\left(\mathbb E\int_0^T\big\|u_\ell(t)-u(t)\big\|_{L^2}^2{\rm d}t\right)^{1/2}\longrightarrow0
\end{align}
where the convergence is justified easily by Dominated Convergence since $u\in L^2(\Omega\times[0,T]\times\mathbb T^3)$. Hence, for $a.e.$ $t\in[0,T]$ (including $t=T$), we have
\begin{align}\label{LEE convergence 1}
	\int_{\mathbb T^3}u_{\ell}(t)\cdot u(t)\varphi(t){\rm d}x\longrightarrow\int_{\mathbb T^3}|u(t)|^2\varphi(t){\rm d}x\quad in~L^1(\Omega)~.
\end{align}
Similarly, one can show that
\begin{equation}\left\{
\begin{aligned}
	&\int_{\mathbb T^3}u_{0,\ell}\cdot u_0\varphi(0){\rm d}x\longrightarrow\int_{\mathbb T^3}\big|u_0\big|^2\varphi(0){\rm d}x~;\\
	&2\nu\int_0^\cdot\int_{\mathbb T^3}\varphi\nabla u_\ell :\nabla u{\rm d}x{\rm d}\tau\longrightarrow2\nu\int_0^\cdot\int_{\mathbb T^3}\varphi|\nabla u|^2{\rm d}x{\rm d}\tau\quad in~L^1(\Omega;C[0,T])~;\\
	&\sum_{k=1}^\infty\int_0^\cdot\int_{\mathbb T^3}\alpha_\ell\ast\varsigma_k(u(\tau))\cdot\varsigma_k(u(\tau))\varphi(\tau){\rm d}x{\rm d}\tau\longrightarrow\sum_{k=1}^\infty\int_0^\cdot\int_{\mathbb T^3}\big|\varsigma_k(u(\tau))\big|^2\varphi(\tau){\rm d}x{\rm d}\tau\quad in~L^1(\Omega;C[0,T]).
\end{aligned}\right.\label{LEE convergence 2}
\end{equation}
Convergence of $\mathcal I^\ell$ is derived by a uniformly integrability argument. Firstly, by that $u\in L^3([0,T]\times\mathbb T^3)$ and $p=\Delta^{-1}{\rm div}{\rm div}(u\otimes u)\in L^{3/2}([0,T]\times\mathbb T^3)$ almost surely, one sees that
\begin{align}
	\mathcal I^\ell\stackrel{C[0,T]}{\longrightarrow}\int_0^\cdot\int_{\mathbb T^3}|u|^2\big(\partial_\tau\varphi+\nu\Delta\varphi\big){\rm d}x{\rm d}\tau+\int_0^\cdot\int_{\mathbb T^3}\big(|u|^2+2p\big)u\cdot\nabla\varphi{\rm d}x{\rm d}\tau~,\quad\mathbb P-a.s.~.\nonumber
\end{align}
Then $L^1(\Omega)$-convergence follows by the fact that $\big\{\mathcal I^\ell\big\}$ is a uniformly integrable family of $C[0,T]$-valued random variables: for example, by H\"older'a inequality, the Sobolev imbedding $L^4\hookrightarrow H^{3/4}$ and interpolation $H^{3/4}=[L^2,H^1]_{1/4}$ , one have
\begin{align}
	\int_0^T\int_{\mathbb T^3}\big|(u_\ell\cdot u)u\cdot\nabla\varphi\big|{\rm d}x{\rm d}t
	&\leq\|\nabla\varphi\|_{L^\infty}\int_0^T\big\|u_\ell(t)\|_{L^2}\|u(t)\|_{L^4}^2{\rm d}t\nonumber\\
	&\lesssim\|\nabla\varphi\|_{L^\infty}\int_0^T\|u(t)\|_{L^2}^{3/2}\|\nabla u(t)\|_{L^2}^{3/2}{\rm d}t\nonumber\\
	&\lesssim\|\nabla\varphi\|_{L^\infty}\left(\int_0^T\|u(t)\|_{L^2}^6{\rm d}t\right)^{1/4}\|\nabla u\|_{L^2}^{3/2},\nonumber
\end{align}
then again by H\"older's inequality and that $u\in L^r\big(\Omega; L^\infty([0,T];L^2)\big)$ (here we need to set $r\geq10$ in (\ref{standard uniform bounds 1})), one can derive that
\begin{align}\label{LEE uniform integrability}
	\mathbb E\left|\int_0^T\int_{\mathbb T^3}(u_\ell\cdot u)u\cdot\nabla\varphi{\rm d}x{\rm d}t\right|^{10/9}
	&\lesssim\|\nabla\varphi\|_{L^\infty}^{10/9}\mathbb E\left[\left(\int_0^T\|u(t)\|_{L^2}^6{\rm d}t\right)^{5/18}\|\nabla u\|_{L^2}^{5/3}\right]\nonumber\\
	&\lesssim\|\nabla\varphi\|_{L^\infty}^{10/9}\left[\mathbb E\left(\int_0^T\|u(t)\|_{L^2}^6{\rm d}t\right)^{5/3}\right]^{1/6}\big(\mathbb E\|\nabla u\|_{L^2}^2\big)^{5/6}<\infty.
\end{align}
Hence
\begin{align}\label{LEE convergence 3}
	\mathcal I^\ell\longrightarrow\int_0^\cdot\int_{\mathbb T^3}|u|^2\big(\partial_\tau\varphi+\nu\Delta\varphi\big){\rm d}x{\rm d}\tau+\int_0^\cdot\int_{\mathbb T^3}\big(|u|^2+2p\big)u\cdot\nabla\varphi{\rm d}x{\rm d}\tau\quad in~L^1(\Omega;C[0,T])~.
\end{align}
Now, to ensure the convergence of the $\mathbf E^\ell$ term, it remains to show  convergence of the martingale term (\ref{martingale term of mollified LEE}).
In particular, we should show that
\begin{align}\label{martingale convergence of mathcal-N^ell}
	\mathcal N^\ell(\varphi)\stackrel{\mathcal M_{2,c}}{\longrightarrow}~2\sum_{m,n\in\mathbb Z^3}\int_0^\cdot{\widehat\varphi}_n(\tau){\widehat u}_{n-m}(\tau)\cdot{\rm d}\mathbf M_\tau(e_m)=:\mathcal N(\varphi)~.
\end{align}
This is quite straight-forward as by It\^o-isometry and a similar calculation as in (\ref{quadratic variation of mathcal-N^ell}) we have
\begin{align}
	\mathbb E\sup_{t\in[0,T]}\big|\mathcal N^\ell_t(\varphi)-\mathcal N_t(\varphi)\big|^2
	&=4\mathbb E\sum_{k=1}^\infty\int_0^T\left[\int_{\mathbb T^3}\Big(\alpha_\ell\ast\varsigma_k(u(t))-\varsigma_k(u(t))\Big)\cdot u(t)\varphi(t){\rm d}x\right]^2{\rm d}t\nonumber\\
	&\leq4\|\varphi\|_{L^\infty}^2\mathbb E\int_0^T\sum_{k=1}^\infty\big\|\alpha_\ell\ast\varsigma_k(u(t))-\varsigma_k(u(t))\big\|_{L^2}^2\|u(t)\|_{L^2}^2{\rm d}t\longrightarrow0~.\nonumber
\end{align}
Here the convergence is justified by Dominated Convergence since
\begin{align}
	&i)~\sum_{k=1}^\infty\big\|\alpha_\ell\ast\varsigma_k(u(t))-\varsigma_k(u(t))\big\|_{L^2}^2\|u(t)\|_{L^2}^2\longrightarrow0~for~a.e.~(\omega,t)\in\Omega\times[0,T]~;\nonumber\\
	&ii)~\sum_{k=1}^\infty\big\|\alpha_\ell\ast\varsigma_k(u(t))-\varsigma_k(u(t))\big\|_{L^2}^2\|u(t)\|_{L^2}^2\leq4\big(1+\|u(t)\|_{L^2}^2\big)\|u(t)\|_{L^2}^2\quad (~\in L^1(\Omega\times[0,T])~)~.\nonumber
\end{align} 
(Note that $i)$ again follows by Dominated Convergence and for $ii)$ we have used (\ref{con. linear growth}) and (\ref{standard uniform bounds 1}). )
Also, we can calculate quadratic variation of $\mathcal N(\varphi)$ :
\begin{align}
	\big\langle\mathcal N(\varphi)\big\rangle_t
	=4\sum_{k=1}^\infty\int_0^t\left(\int_{\mathbb T^3}\varsigma_k(u(\tau))\cdot u(\tau)\varphi(\tau){\rm d}x\right)^2{\rm d}\tau~,\nonumber
\end{align}
from which continuity on $\varphi$ can be derived very similarly as in (\ref{martingale map estimate}).\par Putting together (\ref{LEE convergence 1}), (\ref{LEE convergence 2}), (\ref{LEE convergence 3}) and (\ref{martingale convergence of mathcal-N^ell}), we see that for $a.e.~t\in[0,T]$ (including $t=T$) the limit
\begin{align}
	\lim_{\ell\rightarrow0}\int_0^t\int_{\mathbb T^3}\varphi{\bf E}^\ell{\rm d}x{\rm d}\tau\nonumber
\end{align}
exists in $L^1(\Omega)$. Finally, following \cite{DR00}, a direct integration by part gives that
\begin{align}
	4\mathcal D^\ell(u)
	&=\int_{\mathbb R^3}\nabla\alpha_{\ell}(y)\cdot\delta_yu\big|\delta_yu\big|^2{\rm d}y\nonumber\\
	&= u^i\partial_i\big(u^ju^j\big)_\ell-\partial_i\big(u^iu^ju^j\big)_\ell+2u^j\partial_i\big(u^iu^j\big)_\ell -2u^iu^j\partial_i u_\ell^j\nonumber\\
	&= \nabla\cdot\Big[u\big(|u|^2\big)_\ell-\big(u|u|^2\big)_\ell\Big]+2{\bf E}^\ell\nonumber
\end{align}
and hence we have
\begin{align}
	4\mathcal D_t(u;\varphi)=\lim_{\ell\rightarrow0}\int_0^t\int_{\mathbb T^3}\Big[u\big(|u|^2\big)_\ell-\big(u|u|^2\big)_\ell\Big]\cdot\nabla\varphi{\rm d}x{\rm d}\tau+\lim_{\ell\rightarrow0}\int_0^t\int_{\mathbb T^3}\varphi{\bf E}^\ell{\rm d}x{\rm d}\tau\nonumber
\end{align} 
where the first limit on the right-hand-side converges to 0 almost surely by the fact that $u\in L^3$ and in $L^1(\Omega)$ by a similar uniform integrability argument as (\ref{LEE uniform integrability}). This complete the proof. 
	
\end{proof}

\
\section{Proof of Proposition \ref{Prop. Leray regularization}, \ref{Prop. LEE for Leray regularization} and \ref{Prop. vorticity bounds for Leray regularised solution}}\label{Sec. proof of Leray regularised solution}
\subsection{Regularity estimates for Leray regularization}
The proof of Proposition \ref{Prop. Leray regularization} is based on standard Galerkin approximation and standard regularity estimates. In this subsection we will fix arbitrarily an $\varepsilon>0$ for system (\ref{regularised sto. NS}). The following lemma is used to estimate the noise term $\mathfrak N^\varepsilon$ in Proposition \ref{Prop. Leray regularization}.
\begin{lemma}[\cite{FG95}, Lemma 2.1]\label{Lem. fractional estimate for Ito integrals} 
Let $K$ and $H$ be two Hilbert spaces. And let $\rm W$ be a cylindrical Wiener process with values in $K$, defined on some stochastic basis. Let $r\geq2$ and $0<\alpha<1/2$ be given. Then for any progressively measurable process $G\in L^r\big(\Omega\times[0,T];L_2(K;H)\big)$, the It\^o integral $I_t(G)=\displaystyle\int_0^tG(s){\rm dW}_s$ satisfies the following estimate:
\begin{align}\label{fractional estimate for Ito integral}
	\mathbb E\|I(G)\|_{{\mathcal W}^{\alpha,r}([0,T];H)}^r\leq C(\alpha,r)\,\mathbb E\int_0^T\|G(s)\|_{L_2(K;H)}^r{\rm d}s
\end{align}
with constant $C>0$.
\end{lemma}
\begin{proof}[\bf Proof of Proposition \ref{Prop. Leray regularization}]
It's quite standard to show by Galerkin approximation that there exists martingale weak solution $u^\varepsilon\in L^2\big(\Omega^\varepsilon;C_{loc}([0,+\infty);L^2)\big)\cap L^2\big(\Omega^\varepsilon;L^2_{loc}([0,+\infty);H^1)\big)$ to system (\ref{regularised sto. NS}) such that (\ref{standard uniform bounds 2}) holds, for which one may refer to framework in \cite{FG95}. For the fractional Sobolev estimates (\ref{fraction Sobolev estimate for noise term}) of noise term, it suffices to verify it under Galerkin approximation:
\begin{align}
	u^{\varepsilon,n}(t)={\rm Q_n}u_0^\varepsilon+\int_0^t\Big[-{\rm Q_n}{\rm P}_L{\rm div}\big((\psi_\varepsilon \ast u^{\varepsilon,n})(s)\otimes u^{\varepsilon,n}(s)\big)+\nu \Delta u^{\varepsilon,n}(s)\Big]{\rm d}s\nonumber\\
	+\sum_{k=1}^N\int_0^t{\rm Q_n}\psi_\varepsilon\ast\varsigma_k(u^{\varepsilon,n}(t)){\rm dB}^k_s\nonumber
\end{align}
where $\{{\rm B}^k\}$ is an independent family of one dimensional Brownian motion on some stochastic basis, ${\rm Q_n}$ denotes finite-dimension projection onto subspace of $L^2$, and ${\rm P}_L$ the Leray projector. To use Lemma \ref{Lem. fractional estimate for Ito integrals}, we rewrite the stochastic integral term as $\displaystyle I_t(G^{\varepsilon,n})=\int_0^tG^{\varepsilon,n}(s){\rm dW}_s$ with cylindrical Wiener process ${\rm W}$ in some Hilbert space $K$. If $\{e_k\}$ is an orthonormal basis of $K$, then the operator $G^{\varepsilon,n}(s): K\longrightarrow L^2_\sigma$ takes the form
\begin{align}
	G^{\varepsilon,n}(s)\left(\sum_ka_ke_k\right)=\sum_{k=1}^Na_k{\rm Q_n}\psi_\varepsilon \ast \varsigma_k(u^\varepsilon(s)).\nonumber
\end{align}
Hence $\displaystyle\|G^{\varepsilon,n}(s)\|_{L_2(K;L^2)}^2=\sum_{k}\|G^{\varepsilon,n}(s)e_k\|_{L^2}^2=\sum_{k=1}^N\left\|{\rm Q_n}\psi_\varepsilon \ast \varsigma_k(u^\varepsilon(s))\right\|_{L^2}^2$. Now, by (\ref{fractional estimate for Ito integral}), we have for any $r\geq2$ and $0<\alpha<1/2$ ,
\begin{align}
	\mathbb E\|I(G^{\varepsilon,n})\|_{{\mathcal W}^{\alpha,r}([0,T];L^2)}^r\leq C(\alpha,r)\,\mathbb E\int_0^T\left(\sum_{k=1}^N\left\|\psi_\varepsilon \ast \varsigma_k(u^{\varepsilon,n}(s))\right\|_{L^2}^2\right)^{r/2}{\rm d}s
\end{align}
This indicates (\ref{fraction Sobolev estimate for noise term}) as we are allowed to pass limit $n\rightarrow\infty$ in the above estimate by Galerkin approximation (for the righthand side one would need a Dominated Convergence argument supported by continuity of $\varsigma_k$ 's and condition(\ref{con. linear growth})).

\par Hence, it remains only to show (\ref{standard uniform bounds 1}). For simplicity, we do a priori estimates via equation (\ref{regularised sto. NS}) directly as it is not hard to justify them under Galerkin approximation. We note that the case $m=0$, i.e.
\begin{align}\label{standard uniform estimate, L^2}
	\mathbb E\sup_{t\in [0,T]}\|u^\varepsilon(t)\|_{L^2}^r\leq C_1(\|u_0\|_{L^2},T,r).
\end{align}
is classical. See Appendix in \cite{FG95} for example. We now show the case $m=1$. For convenience, we denote by $C_\sigma$ the constant for which condition (\ref{con. linear growth}) holds.
\par Applying $\nabla$ on (\ref{regularised sto. NS}) and using It\^o's formula, it is not hard to verify that
\begin{align}\label{equality by Ito formula}
	&e^{-bt}\big(1+\|\nabla u^\varepsilon(t)\|_{L^2}^2\big)^{r/2}+r\nu\int_0^te^{-bs}\big(1+\|\nabla u^\varepsilon(s)\|_{L^2}^2\big)^{\frac{r}{2}-1}\|\nabla^2u^\varepsilon(s)\|_{L^2}^2{\rm d}s\nonumber\\
	&+b\int_0^te^{-bs}\big(1+\|\nabla u^\varepsilon(s)\|_{L^2}^2\big)^{\frac{r}{2}}{\rm d}s\nonumber\\
	=&\big(1+\|\nabla u^\varepsilon_0\|_{L^2}^2\big)^{r/2}-r\int_0^te^{-bs}\big(1+\|\nabla u^\varepsilon(s)\|_{L^2}^2\big)^{\frac{r}{2}-1}\Big\langle\nabla[\psi_\varepsilon\ast u^\varepsilon_i(s)~\partial_i u^\varepsilon(s)],\nabla u^\varepsilon(s)\Big\rangle_{L^2}{\rm d}s\nonumber\\
	&+\frac{r}{2}\int_0^te^{-bs}\big(1+\|\nabla u^\varepsilon(s)\|_{L^2}^2\big)^{\frac{r}{2}-1}\sum_{k=1}^N\big\|\nabla[\psi_\varepsilon\ast\varsigma_k(u^\varepsilon(s))]\big\|_{L^2}^2{\rm d}s\nonumber\\
	&+\frac{r(r-2)}{8}\int_0^te^{-bs}\big(1+\|\nabla u^\varepsilon(s)\|_{L^2}^2\big)^{\frac{r}{2}-2}\sum_{k=1}^N\Big\langle\nabla[\psi_\varepsilon\ast\varsigma_k(u^\varepsilon(s))],\nabla u^\varepsilon(s)\Big\rangle_{L^2}^2{\rm d}s\nonumber\\
	&+r\sum_{k=1}^{N}\int_0^te^{-bs}\big(1+\|\nabla u^\varepsilon(s)\|_{L^2}^2\big)^{\frac{r}{2}-1}\Big\langle\nabla[\psi_\varepsilon\ast\varsigma_k(u^\varepsilon(s))],\nabla u^\varepsilon(s)\Big\rangle_{L^2}{\rm dB}_s^k 
\end{align} 
for any constant $b\in\mathbb R$. By a standard stopping time argument, one can show that
\begin{align}\label{moment estimates 1}
	&e^{-bt}\mathbb E\big(1+\|\nabla u^\varepsilon(t)\|_{L^2}^2\big)^{r/2}+r\nu\mathbb E\int_0^te^{-bs}\big(1+\|\nabla u^\varepsilon(s)\|_{L^2}^2\big)^{\frac{r}{2}-1}\|\nabla^2u^\varepsilon(s)\|_{L^2}^2{\rm d}s\nonumber\\
	&+b~\mathbb E\int_0^te^{-bs}\big(1+\|\nabla u^\varepsilon(s)\|_{L^2}^2\big)^{\frac{r}{2}}{\rm d}s\nonumber\\
	\leq &\big(1+\|\nabla u^\varepsilon_0\|_{L^2}^2\big)^{r/2}+\underbrace{r~\mathbb E\int_0^te^{-bs}\big(1+\|\nabla u^\varepsilon(s)\|_{L^2}^2\big)^{\frac{r}{2}-1}\Big|\Big\langle\nabla[\psi_\varepsilon\ast u^\varepsilon_i(s)~\partial_i u^\varepsilon(s)],\nabla u^\varepsilon(s)\Big\rangle_{L^2}\Big|{\rm d}s}_{:=L^\varepsilon_1(t)}\nonumber\\
	&+\underbrace{\frac{r}{2}\mathbb E\int_0^te^{-bs}\big(1+\|\nabla u^\varepsilon(s)\|_{L^2}^2\big)^{\frac{r}{2}-1}\sum_{k=1}^N\big\|\nabla[\psi_\varepsilon\ast\varsigma_k(u^\varepsilon(s))]\big\|_{L^2}^2{\rm d}s}_{:=L^\varepsilon_2(t)}\nonumber\\
	&+\underbrace{\frac{r(r-2)}{8}\mathbb E\int_0^te^{-bs}\big(1+\|\nabla u^\varepsilon(s)\|_{L^2}^2\big)^{\frac{r}{2}-2}\sum_{k=1}^N\Big\langle\nabla[\psi_\varepsilon\ast\varsigma_k(u^\varepsilon(s))],\nabla u^\varepsilon(s)\Big\rangle_{L^2}^2{\rm d}s}_{:=L^\varepsilon_3(t)}.
\end{align}
By H\"older's inequality, the Sobolev embedding $H^2\hookrightarrow L^\infty$ and standard mollification estimates, we have
\begin{align}\label{triple term estimate}
	\Big|\Big\langle\nabla[\psi_\varepsilon\ast u^\varepsilon_i(s)~\partial_i u^\varepsilon(s)],\nabla u^\varepsilon(s)\Big\rangle_{L^2}\Big|
	&=\Big|\Big\langle\nabla[\psi_\varepsilon\ast u^\varepsilon(s)]\otimes u^\varepsilon(s),\nabla^2 u^\varepsilon(s)\Big\rangle_{L^2}\Big|\nonumber\\
	&\leq\big\|\nabla[\psi_\varepsilon\ast u^\varepsilon(s)]\big\|_{L^\infty}\big\|u^\varepsilon(s)\big\|_{L^2}\big\|\nabla^2 u^\varepsilon(s)\big\|_{L^2}\nonumber\\
	&\leq\big\|\nabla[\psi_\varepsilon\ast u^\varepsilon(s)]\big\|_{H^2}\big\|u^\varepsilon(s)\big\|_{L^2}\big\|\nabla^2 u^\varepsilon(s)\big\|_{L^2}\nonumber\\
		&\lesssim\varepsilon^{-3}\big\|u^\varepsilon(s)\big\|_{L^2}^{2}\big\|\nabla^2 u^\varepsilon(s)\big\|_{L^2}.
\end{align}
Hence, we have for $L^\varepsilon_1$ that
\begin{align}
	L^\varepsilon_1(T)
	\leq & r~\mathbb E\int_0^Te^{-bs}\big(1+\|\nabla u^\varepsilon(s)\|_{L^2}^2\big)^{\frac{r}{2}-1}\Big(C_U\varepsilon^{-3}\big\|u^\varepsilon(s)\big\|_{L^2}^{2}\big\|\nabla^2 u^\varepsilon(s)\big\|_{L^2}\Big){\rm d}s\nonumber\\
	\leq & r~\mathbb E\int_0^Te^{-bs}\big(1+\|\nabla u^\varepsilon(s)\|_{L^2}^2\big)^{\frac{r}{2}-1}\left(\frac{2C_U^2}{\nu\varepsilon^6}\big\|u^\varepsilon(s)\big\|_{L^2}^{4}+\frac{\nu}{2}\big\|\nabla^2 u^\varepsilon(s)\big\|_{L^2}^2\right){\rm d}s\nonumber\\
	\leq &\frac{2rC_U^2}{\nu\varepsilon^6}\mathbb E\int_0^Te^{-bs}\big(1+\|\nabla u^\varepsilon(s)\|_{L^2}^2\big)^{\frac{r}{2}-1}\big\|u^\varepsilon(s)\big\|_{L^2}^{4}{\rm d}s\nonumber\\
	&+\frac{r\nu}{2}\mathbb E\int_0^Te^{-bs}\big(1+\|\nabla u^\varepsilon(s)\|_{L^2}^2\big)^{\frac{r}{2}-1}\|\nabla^2u^\varepsilon(s)\|_{L^2}^2{\rm d}s.\nonumber
\end{align}
Here, $C_U>0$ is an universal constant, and we have used Young's inequality. Then by H\"older's inequality and (\ref{standard uniform estimate, L^2}), we have
\begin{align}
	&\mathbb E\int_0^Te^{-bs}\big(1+\|\nabla u^\varepsilon(s)\|_{L^2}^2\big)^{\frac{r}{2}-1}\big\|u^\varepsilon(s)\big\|_{L^2}^{4}{\rm d}s\nonumber\\
	\leq &\int_0^Te^{-bs}\left[\mathbb E\big(1+\|\nabla u^\varepsilon(s)\|_{L^2}^2\big)^{\frac{r}{2}}\right]^{\frac{r-2}{r}}\left(\mathbb E\big\|u^\varepsilon(s)\big\|_{L^2}^{2r}\right)^{\frac{2}{r}}{\rm d}s\nonumber\\
	\leq &\left(\sup_{t\in[0,T]}\mathbb E\big\|u^\varepsilon(t)\big\|_{L^2}^{2r}\right)^{\frac{2}{r}}\int_0^Te^{-bs}\mathbb E\big(1+\|\nabla u^\varepsilon(s)\|_{L^2}^2\big)^{\frac{r}{2}}{\rm d}s\nonumber\\
	\leq &C_1(\|u_0\|_{L^2},T,2r)^{\frac{2}{r}}\int_0^Te^{-bs}\mathbb E\big(1+\|\nabla u^\varepsilon(s)\|_{L^2}^2\big)^{\frac{r}{2}}{\rm d}s.\nonumber
\end{align}
Hence, we have
\begin{align}\label{L^varepsilon_1 estimate}
	L^\varepsilon_1(T)
	\leq &\frac{2rC_U^2}{\nu\varepsilon^6}C_1(\|u_0\|_{L^2},T,2r)^{\frac{2}{r}}\int_0^Te^{-bs}\mathbb E\big(1+\|\nabla u^\varepsilon(s)\|_{L^2}^2\big)^{\frac{r}{2}}{\rm d}s\nonumber\\
	&+\frac{r\nu}{2}\mathbb E\int_0^Te^{-bs}\big(1+\|\nabla u^\varepsilon(s)\|_{L^2}^2\big)^{\frac{r}{2}-1}\|\nabla^2u^\varepsilon(s)\|_{L^2}^2{\rm d}s.
\end{align}
For $L^\varepsilon_2$, we use standard mollification estimates, (\ref{con. linear growth}) and Poincar\'e inequality to write
\begin{align}\label{L^varepsilon_2 estimate}
	L^\varepsilon_2(T)
	\leq &\frac{r}{2}C_\sigma\varepsilon^{-2}\mathbb E\int_0^Te^{-bs}\big(1+\|\nabla u^\varepsilon(s)\|_{L^2}^2\big)^{\frac{r}{2}-1}\left(1+\big\|u^\varepsilon(s)\big\|_{L^2}^{2}\right){\rm d}s\nonumber\\
	\leq &\frac{r}{2}C_UC_\sigma\varepsilon^{-2}\mathbb E\int_0^Te^{-bs}\big(1+\|\nabla u^\varepsilon(s)\|_{L^2}^2\big)^{\frac{r}{2}}{\rm d}s.
\end{align}
Here, $C_U>0$ is an universal constant. Similarly, for $L^\varepsilon_3$, we have
\begin{align}\label{L^varepsilon_3 estimate}
	L^\varepsilon_3(T)
	\leq &\frac{r(r-2)}{8}\mathbb E\int_0^te^{-bs}\big(1+\|\nabla u^\varepsilon(s)\|_{L^2}^2\big)^{\frac{r}{2}-1}\sum_{k=1}^N\big\|\nabla[\psi_\varepsilon\ast\varsigma_k(u^\varepsilon(s))]\big\|_{L^2}^2{\rm d}s\nonumber\\
	\leq &\frac{r(r-2)}{8}C_UC_\sigma\varepsilon^{-2}\mathbb E\int_0^Te^{-bs}\big(1+\|\nabla u^\varepsilon(s)\|_{L^2}^2\big)^{\frac{r}{2}}{\rm d}s.
\end{align}
Now, inserting (\ref{L^varepsilon_1 estimate}), (\ref{L^varepsilon_2 estimate}) and (\ref{L^varepsilon_3 estimate}) into (\ref{moment estimates 1}), and by choosing
\begin{align}\label{constraint on const. b}
	b\geq\frac{2rC_U^2}{\nu\varepsilon^6}C_1(\|u_0\|_{L^2},T,2r)^{\frac{2}{r}}+\frac{r}{2}C_UC_\sigma\varepsilon^{-2}+\frac{r(r-2)}{8}C_UC_\sigma\varepsilon^{-2},
\end{align}
we get
\begin{align}
	&\sup_{t\in[0,T]}e^{-bt}\mathbb E\big(1+\|\nabla u^\varepsilon(t)\|_{L^2}^2\big)^{r/2}+\frac{r\nu}{2}\mathbb E\int_0^Te^{-bs}\big(1+\|\nabla u^\varepsilon(s)\|_{L^2}^2\big)^{\frac{r}{2}-1}\|\nabla^2u^\varepsilon(s)\|_{L^2}^2{\rm d}s\nonumber\\
	\leq &\big(1+\|\nabla u^\varepsilon_0\|_{L^2}^2\big)^{r/2}.\nonumber
\end{align}
This shows that
\begin{align}\label{not-finished uniform estimate}
	\sup_{t\in[0,T]}\mathbb E\big(1+\|\nabla u^\varepsilon(t)\|_{L^2}^2\big)^{r/2}\leq D\big(\varepsilon,\|u_0^\varepsilon\|_{H^1},T,r\big)
\end{align}
for some constant $D>0$. Now, return on (\ref{equality by Ito formula}), we have that
\begin{align}\label{not-finished uniform estimate 2}
	&\mathbb E\sup_{t\in[0,T]}e^{-bt}\big(1+\|\nabla u^\varepsilon(t)\|_{L^2}^2\big)^{r/2}+r\nu\mathbb E\int_0^Te^{-bs}\big(1+\|\nabla u^\varepsilon(s)\|_{L^2}^2\big)^{\frac{r}{2}-1}\|\nabla^2u^\varepsilon(s)\|_{L^2}^2{\rm d}s\nonumber\\
	&+b\mathbb E\int_0^Te^{-bs}\big(1+\|\nabla u^\varepsilon(s)\|_{L^2}^2\big)^{\frac{r}{2}}{\rm d}s\nonumber\\
	\leq &\big(1+\|\nabla u^\varepsilon_0\|_{L^2}^2\big)^{r/2}+L^\varepsilon_1(T)+L^\varepsilon_2(T)+L^\varepsilon_3(T)\nonumber\\
	&+r~\mathbb E\sup_{t\in[0,T]}\left|\sum_{k=1}^{N}\int_0^te^{-bs}\big(1+\|\nabla u^\varepsilon(s)\|_{L^2}^2\big)^{\frac{r}{2}-1}\Big\langle\nabla[\psi_\varepsilon\ast\varsigma_k(u^\varepsilon(s))],\nabla u^\varepsilon(s)\Big\rangle_{L^2}{\rm dB}_s^k\right|.
\end{align}
By BDG-inequality, standard mollification estimates, (\ref{con. linear growth}), Poincar\'e inequality, H\"older's inequality and (\ref{not-finished uniform estimate}), we have
\begin{align}
	&\mathbb E\sup_{t\in[0,T]}\left|\sum_{k=1}^{N}\int_0^te^{-bs}\big(1+\|\nabla u^\varepsilon(s)\|_{L^2}^2\big)^{\frac{r}{2}-1}\Big\langle\nabla[\psi_\varepsilon\ast\varsigma_k(u^\varepsilon(s))],\nabla u^\varepsilon(s)\Big\rangle_{L^2}{\rm dB}_s^k\right|\nonumber\\
	\leq &\mathbb E\left(\int_0^Te^{-2bs}\big(1+\|\nabla u^\varepsilon(s)\|_{L^2}^2\big)^{r-2}\sum_{k=1}^{N}\Big\langle\nabla[\psi_\varepsilon\ast\varsigma_k(u^\varepsilon(s))],\nabla u^\varepsilon(s)\Big\rangle_{L^2}^2{\rm d}s\right)^{1/2}\nonumber\\
	\leq &\varepsilon^{-2}\mathbb E\left(\int_0^Te^{-2bs}\big(1+\|\nabla u^\varepsilon(s)\|_{L^2}^2\big)^{r-1}\sum_{k=1}^{N}\big\|\varsigma_k(u^\varepsilon(s))\big\|_{L^2}^2{\rm d}s\right)^{1/2}\nonumber\\
	\leq &C_UC_\sigma^{1/2}\varepsilon^{-2}\left(\int_0^Te^{-2bs}\mathbb E\big(1+\|\nabla u^\varepsilon(s)\|_{L^2}^2\big)^{r}{\rm d}s\right)^{1/2}\leq\frac{C_UC_\sigma^{1/2}D^{1/2}}{\varepsilon^2\sqrt{2b}}.\nonumber
\end{align}
Inserting this into (\ref{not-finished uniform estimate 2}), then by (\ref{L^varepsilon_1 estimate}-\ref{L^varepsilon_3 estimate}) and (\ref{constraint on const. b}), we finally get
\begin{align}\label{standard uniform estimate, H^1}
	&\mathbb E\sup_{t\in[0,T]}e^{-bt}\big(1+\|\nabla u^\varepsilon(t)\|_{L^2}^2\big)^{r/2}+p\nu\mathbb E\int_0^Te^{-bs}\big(1+\|\nabla u^\varepsilon(s)\|_{L^2}^2\big)^{\frac{p}{2}-1}\|\nabla^2u^\varepsilon(s)\|_{L^2}^2{\rm d}s\nonumber\\
	\leq &\big(1+\|\nabla u^\varepsilon_0\|_{L^2}^2\big)^{r/2}+\frac{C_UC_\sigma^{1/2}D^{1/2}}{\varepsilon^2\sqrt{2b}}.
\end{align}
This proves (\ref{standard uniform bounds 1}) for $m=1$. For $m\geq2$, one may estimate inductively by using (\ref{standard uniform estimate, L^2}) and (\ref{standard uniform estimate, H^1}). We should not repeat the process here as the arguments are quite similar as above.
\end{proof}

\
\subsection{Local energy balance for Leray regularized Navier-Stokes}
The proof of Proposition \ref{Prop. LEE for Leray regularization} combines a rigorous justification of application of It\^o's formula and another of the martingale term in the formula.
\begin{proof}[Proof of Proposition \ref{Prop. LEE for Leray regularization}]
Still, let $\big\{\big(\Omega^\varepsilon, \mathscr{F}^\varepsilon,\{\mathscr{F}^\varepsilon_t\}_{t\geq0}, \mathbb{P}^\varepsilon\big),u^\varepsilon\big\}$ be a martingale weak solution to system (\ref{regularised sto. NS}), given by Proposition \ref{Prop. Leray regularization}, with initial value $\psi_\varepsilon\ast u_0=:u_0^\varepsilon$. Also $p^\varepsilon:=(-\Delta)^{-1}{\rm divdiv}\big(u^\varepsilon\otimes u^\varepsilon\big)$ is the associated pressure term. As we are arguing with arbitrarily fixed $\varepsilon>0$, in this subsection we would simply denote $\big(\Omega, \mathscr{F}, \mathbb{P}\big)=\big(\Omega^\varepsilon, \mathscr{F}^\varepsilon, \mathbb{P}^\varepsilon\big)$ for notational simplicity. Then it is clear that for eachtest function $\varphi\in \big(C^\infty(\mathbb T^3)\big)^3$, the process
\begin{align}\label{M^varepsilon}
	{\mathbf M}^\varepsilon_t(\varphi):=\langle u^\varepsilon(t),\varphi\rangle-\langle u_0^\varepsilon,\varphi\rangle+\int_0^t\Big\langle{\rm div}\big((\psi_\varepsilon\ast u^\varepsilon)(s)\otimes u^\varepsilon(s)\big)+\nabla p^\varepsilon(s)-\nu\Delta u^\varepsilon(s),\varphi\Big\rangle{\rm d}s
\end{align}
is a continuous square-integrable martingale with quadratic variation
\begin{align}\label{quadratic variation of M^varepsilon}
	\big\langle \mathbf M^\varepsilon(\varphi)\big\rangle_t=\sum_{k=1}^N\int_0^t\big\langle\psi_\varepsilon\ast\varsigma_k(\tilde u^\varepsilon(s)),\varphi\big\rangle^2{\rm d}s.
\end{align}

\par We first show an energy equality for $u^\varepsilon$. For this, fixed arbitrarily an $x_0\in\mathbb T^3$ and take $\tau_{x_0}\alpha_\ell=\alpha_\ell(\cdot-x_0)$ as test function in (\ref{M^varepsilon}), where $\alpha_\ell$ is a standard mollifier. By the smoothness of $u^\varepsilon(t,\cdot)$ for each $t\geq0$, each term on the righthand side converges $\mathbb P$-$a.s.$ if we send $\ell\rightarrow0$. On the other hand, by (\ref{quadratic variation of M^varepsilon}), it can be easily seen that $\big\{\mathbf M^\varepsilon(\tau_{x_0}\alpha_\ell)\big\}_{\ell}$ is Cauchy in space $(\mathcal M_{c,2})^3$ (to be specific, here we mean that take $\varphi_1^{\ell,x_0}=\big(\tau_{x_0}\alpha_\ell,0,0\big)^T,\varphi_2^{\ell,x_0}=\big(0,\tau_{x_0}\alpha_\ell,0\big)^T,\varphi_3^{\ell,x_0}=\big(0,0,\tau_{x_0}\alpha_\ell\big)^T$ in (\ref{M^varepsilon}) seperately and, by slightly abusing the notation, we still denote by $\mathbf M^\varepsilon(\tau_{x_0}\alpha_\ell)$ the three-dimensional martingale $\left(\mathbf M^\varepsilon(\varphi_1^{\ell,x_0}),\mathbf M^\varepsilon(\varphi_2^{\ell,x_0}),\mathbf M^\varepsilon(\varphi_3^{\ell,x_0})\right)^T$ that we get here) :
\begin{align}
	&\mathbb {E}\sup_{0\leq t\leq T}\big|\mathbf M^\varepsilon_t(\tau_{x_0}\alpha_{\ell})-\mathbf M^\varepsilon_t(\tau_{x_0}\alpha_{\ell'})\big|^2\nonumber\\
	=&\mathbb {E}\sup_{0\leq t\leq T}\big|\mathbf M^\varepsilon_t(\tau_{x_0}\alpha_{\ell}-\tau_{x_0}\alpha_{\ell'})\big|^2\nonumber\\
	=&\sum_{k=1}^N\mathbb {E}\int_{0}^{T}\big|\big\langle\psi_\varepsilon\ast\varsigma_k(u^\varepsilon(s)),\tau_{x_0}\alpha_{\ell}-\tau_{x_0}\alpha_{\ell'}\big\rangle\big|^2{\rm d}s\nonumber\\
	=&\sum_{k=1}^N\mathbb {E}\int_{0}^{T}\Big|\alpha_{\ell}\ast\big[\psi_\varepsilon\ast\varsigma_k(u^\varepsilon(s))\big](x_0)-\alpha_{\ell'}\ast\big[\psi_\varepsilon\ast\varsigma_k(u^\varepsilon(s))\big](x_0)\Big|^2{\rm d}s\longrightarrow0,~as~\ell,\ell'\rightarrow0\nonumber
\end{align}
where the limit is justified by Dominated Convergence Theorem. So now we deduce for each $x_0\in\mathbb T^3$ that
\begin{align}
	u^\varepsilon(\cdot,x_0)=u_0^\varepsilon(x_0)
	&+\int_0^\cdot\Big[-{\rm div}\big((\psi_\varepsilon\ast u^\varepsilon)(s,x_0)\otimes u^\varepsilon(s,x_0)\big)-\nabla p^\varepsilon(s,x_0)+\nu\Delta u^\varepsilon(s,x_0)\Big]{\rm d}s+\mathbf{\bar M}^\varepsilon(x_0)\nonumber
\end{align}
where we denote $\displaystyle\mathbf{\bar M}^\varepsilon(x_0)=\lim_{\ell\rightarrow0}\mathbf M^\varepsilon(\tau_{x_0}\alpha_\ell)\in (\mathcal M_{c,2})^3$. Note that all processes under the integrals are progressively measurable. And it is not hard to show that
\begin{align}
	\big\langle \mathbf{\bar M}^\varepsilon(x_0)^i\big\rangle=\sum_{k=1}^N\int_0^\cdot\psi_\varepsilon\ast\varsigma_k(u^\varepsilon(s))^i(x_0)^2{\rm d}s,\quad i=1,2,3.
\end{align}
Here, $\mathbf{\bar M}^\varepsilon_t(x_0)^i$ denotes the $i$-th component of the vector $\mathbf{\bar M}_t^\varepsilon(x_0)$, and $\psi_\varepsilon\ast\varsigma_k(u^\varepsilon(s))^i(x_0)$ the $i$-component of $\psi_\varepsilon\ast\varsigma_k(u^\varepsilon(s))$ valued at point $x_0$.
Also, one can calculate the cross-variation from two point
\begin{align}\label{cross-variation of bar-M^varepsilon}
	\big\langle\mathbf{\bar M}^\varepsilon(x_1)^i,\mathbf{\bar M}^\varepsilon(x_2)^j\big\rangle_t
	&=\frac{\big\langle\mathbf{\bar M}^\varepsilon(x_1)^i+\mathbf{\bar M}^\varepsilon(x_2)^j\big\rangle_t-\big\langle\mathbf{\bar M}^\varepsilon(x_1)^i\big\rangle_t-\big\langle\mathbf{\bar M}^\varepsilon(x_2)^j\big\rangle_t}{2}\nonumber\\
	&=\lim_{\ell\rightarrow0}\frac{\big\langle\mathbf M^\varepsilon(\varphi_i^{\ell,x_1})+\mathbf M^\varepsilon(\varphi_j^{\ell,x_2})\big\rangle_t-\big\langle\mathbf M^\varepsilon(\tau_{x_1}\alpha_\ell)^i\big\rangle_t-\big\langle\mathbf M^\varepsilon(\tau_{x_2}\alpha_\ell)^j\big\rangle_t}{2}\nonumber\\
	&=\sum_{k=1}^N\int_0^t\psi_\varepsilon\ast\varsigma_k(u^\varepsilon(s))^i(x_1)~\psi_\varepsilon\ast\varsigma_k(u^\varepsilon(s))^j(x_2){\rm d}s,\quad i,j=1,2,3.
\end{align}
Now, taking any non-negative $\varphi\in C^\infty_c\left((0,T)\times\mathbb{T}^3\right)$, we are allowed to apply It\^o's formula to the process $t\mapsto\frac{1}{2}\big|u^\varepsilon(t,x)\big|^2\varphi(t,x)$ for arbitrarily fixed $x\in\mathbb T^3$. Similarly as in Appendix \ref{Appen. Intuitive computation}, integrating over $\mathbb T^3$ on both sides and by integration by part, we have $\mathbb P$-$a.s.$ that
\begin{align}\label{LEE-b for u^varepsilon}
	&\ \ \ \ \ \ \int_{\mathbb{T}^3}\big|u^\varepsilon(t)\big|^2\varphi(t){\rm d}x+2\nu\int_0^t\int_{\mathbb{T}^3}\big|\nabla u^\varepsilon(s)\big|^2\varphi(s){\rm d}x{\rm d}s\nonumber\\
	&=\int_0^t\int_{\mathbb{T}^3}\big|u^\varepsilon(s)\big|^2\Big(\partial_\tau\varphi(s)+\nu\Delta\varphi(s)\Big){\rm d}x{\rm d}s+\sum_{k=1}^N\int_0^t\int_{\mathbb{T}^3}|\psi_\varepsilon\ast\varsigma_k(u^\varepsilon(s))|^2\varphi(s){\rm d}x{\rm d}s\nonumber\\
	&\ \ +\int_0^t\int_{\mathbb{T}^3}\left(\big|u^\varepsilon(s)\big|^2\big(\psi_\varepsilon\ast u^\varepsilon(s)\big)+2 p^\varepsilon(s) u^\varepsilon(s)\right)\cdot\nabla\varphi(s){\rm d}x{\rm d}s+\underbrace{2\int_{\mathbb{T}^3}\int_0^t\varphi(s,x) u^\varepsilon(s,x)\cdot{\rm d}\mathbf{\bar M}^\varepsilon_s(x){\rm d}x}_{=:\mathbf N^\varepsilon_t(\varphi)}
	.
\end{align}

\par To this end, we should show that $\mathbf N^\varepsilon(\varphi)\in \mathcal M_{c,2}$. The strategy is to approximate the spatial integral by a ``Riemann sum" which clearly belong to $\mathcal M_{c,2}$ : denote by $\Lambda=\{Q_i: 1\leq i\leq K\}$ a finite partition of $\mathbb T^3$, and from each $Q_i$ we pick an $x_i\in Q_i$, then we define
\begin{align}
	\mathbf N^{\varepsilon,\Lambda}_t(\varphi):=2\sum_{i=1}^K|Q_i|\int_0^t\varphi(s,x_i) u^\varepsilon(s,x_i)\cdot{\rm d}\mathbf{\bar M}^\varepsilon_s(x_i),\quad t\in[0,T].
\end{align}
To see that, for each fixed $t\in[0,T]$, $\mathbf N^{\varepsilon,\Lambda}_t(\varphi)\rightarrow\mathbf N^{\varepsilon}_t(\varphi)$ almost surely as $\displaystyle\|\Lambda\|:=\min_{1\leq i\leq K}|Q_i|\rightarrow0$, it suffices to show that $\displaystyle x\mapsto\int_0^t\varphi(s,x) u^\varepsilon(s,x)\cdot{\rm d}\mathbf{\bar M}^\varepsilon_s(x)$ is continuous. This is indeed the case as by (\ref{cross-variation of bar-M^varepsilon}) we have
\begin{align}
	&\mathbb {E}\sup_{0\leq t\leq T}\left|\int_0^t\varphi(s,x) u^\varepsilon(s,x)\cdot{\rm d}\mathbf{\bar M}^\varepsilon_s(x)-\int_0^t\varphi(s,y) u^\varepsilon(s,y)\cdot{\rm d}\mathbf{\bar M}^\varepsilon_s(y)\right|^2\nonumber\\
	=&\mathbb {E}\int_0^T\sum_{k=1}^N\big|\varphi(s,x) u^\varepsilon(s,x)\cdot[\psi_\varepsilon\ast\varsigma_k(\tilde u^\varepsilon(s))](x)-\varphi(s,y) u^\varepsilon(s,y)\cdot[\psi_\varepsilon\ast\varsigma_k(\tilde u^\varepsilon(s))](y)\big|^2{\rm d}s.\nonumber
\end{align}
On the other hand, let $\Lambda_1$,$\Lambda_2$ be two partitions of $\mathbb T^3$ and $\Lambda'=\{Q_i': 1\leq i\leq K'\}$ be a refinement of them such that
\begin{align}
	\mathbf N^{\varepsilon,\Lambda_1}_t(\varphi)
	&=2\sum_{i=1}^{K'}|Q_i'|\int_0^t\varphi(s,x_i) u^\varepsilon(s,x_i)\cdot{\rm d}\mathbf{\bar M}^\varepsilon_s(x_i)\nonumber\\
	\mathbf N^{\varepsilon,\Lambda_2}_t(\varphi)
	&=2\sum_{i=1}^{K'}|Q_i'|\int_0^t\varphi(s,y_i) u^\varepsilon(s,y_i)\cdot{\rm d}\mathbf{\bar M}^\varepsilon_s(y_i)\nonumber
\end{align}
where $x_i,y_i\in Q_i'$. Then from the previous equality we have
\begin{align}
	&\mathbb {E}\sup_{0\leq t\leq T}\left|\mathbf N^{\varepsilon,\Lambda_1}_t(\varphi)-\mathbf N^{\varepsilon,\Lambda_2}_t(\varphi)\right|^2\nonumber\\
	\lesssim & \sum_{i=1}^{K'}|Q_i'|~\mathbb {E}\int_0^T\sum_{k=1}^N\Big|\varphi(s,x_i) u^\varepsilon(s,x_i)\cdot\big[\psi_\varepsilon\ast\varsigma_k( u^\varepsilon(s))\big](x_i)-\varphi(s,y_i)\tilde u^\varepsilon(s,y_i)\cdot\big[\psi_\varepsilon\ast\varsigma_k( u^\varepsilon(s))\big](y_i)\Big|^2{\rm d}s\nonumber\\
	\lesssim & \sup_{|x-y|\leq\|\Lambda'\|}\mathbb {E}\int_0^T\sum_{k=1}^N\Big|\varphi(s,x) u^\varepsilon(s,x)\cdot\big[\psi_\varepsilon\ast\varsigma_k( u^\varepsilon(s))\big](x)-\varphi(s,y) u^\varepsilon(s,y)\cdot\big[\psi_\varepsilon\ast\varsigma_k( u^\varepsilon(s))\big](y)\Big|^2{\rm d}s\nonumber
\end{align}
where $\displaystyle\|\Lambda'\|:=\min_{1\leq i\leq K}|Q_i'|$. This coupled with a Dominated Convergence argument indicates that $\big\{\mathbf N^{\varepsilon,\Lambda}(\varphi)\big\}_\Lambda$ is Cauchy in space $\mathcal M_{c,2}$. And hence its almost sure limit $\mathbf N^\varepsilon(\varphi)$ belongs to $\mathcal M_{c,2}$. Now from (\ref{cross-variation of bar-M^varepsilon}) we can calculate quadratic variation of $\mathbf N^{\varepsilon,\Lambda}(\varphi)$ :
\begin{align}
	\big\langle \mathbf N^{\varepsilon,\Lambda}(\varphi)\big\rangle_t	=&4\sum_{i,j=1}^K|Q_i||Q_j|\sum_{p,q=1}^3\int_0^t\varphi(s,x_i) u^\varepsilon_p(s,x_i)\varphi(s,x_j) u^\varepsilon_q(s,x_j){\rm d}\big\langle\mathbf{\bar M}^\varepsilon(x_i)^p,\mathbf{\bar M}^\varepsilon(x_j)^q\big\rangle_s\nonumber\\
	=&4\sum_{k=1}^N\int_0^t\left(\sum_{i=1}^K|Q_i|\varphi(s,x_i) u^\varepsilon(s,x_i)\cdot[\psi_\varepsilon\ast\varsigma_k( u^\varepsilon(s))](x_i)\right)^2{\rm d}s\nonumber
\end{align}
(here $u^\varepsilon_p$ denotes $p$-th component of $u^\varepsilon$) which indicates 
\begin{align}
	\big\langle \mathbf N^\varepsilon(\varphi)\big\rangle_t=4\sum_{k=1}^N\int_0^t\left(\int_{\mathbb{T}^3}\psi_\varepsilon\ast\varsigma_k( u^\varepsilon(s))\cdot u^\varepsilon(s)\varphi(s){\rm d}x\right)^2{\rm d}s,\nonumber
\end{align}
i.e. (\ref{quadratic variation of N^varepsilon}). So putting together, we see that for each $\varepsilon$ and test function $\varphi$, $\mathbf N^\varepsilon(\varphi)$ is a continuous square-integrable martingale (w.r.t. the filtration $\{\mathscr{F}_t^\varepsilon:=\sigma(u^\varepsilon(s);0\leq s\leq t)\}_{t\geq0}$ ) with quadratic variation given by (\ref{quadratic variation of N^varepsilon}).
\end{proof}

\
\subsection{Vorticity bounds for Leray regularization} The proof of Proposition \ref{Prop. vorticity bounds for Leray regularised solution} relies on a well-chosen change of dependent variables which is inspired by \cite{Con90}. The difference is that now every thing should be done in the context of stochastic calculus and Leray regularization. And we also have to bound the martingale term and the quadratic variation term arising from It\^o's formula.

\begin{proof}[\bf Proof of Proposition \ref{Prop. vorticity bounds for Leray regularised solution}] Let $\varepsilon>0$ and $0<T<+\infty$ be fixed arbitrarily. And let the initial vorticity $\omega_0$ be given as in the proposition. Take a standard mollifier $\alpha_\ell(\cdot-x_0)$ ($x_0\in\mathbb T^3$) as test function on both sides of (\ref{regularised sto. NS}) and let $\ell\rightarrow0$. Then by the smoothness of $u^\varepsilon(t,\cdot)$ for each $t\geq0$, all the terms converge $\mathbb P$-$a.s.$, and one have $\mathbb P$-$a.s.$ that
\begin{align}
	u^\varepsilon(t,x_0)=u_0^\varepsilon(x_0)
	&+\int_0^t\Big(-{\rm div}\big((\psi_\varepsilon\ast u^\varepsilon)(s,x_0)\otimes u^\varepsilon(s,x_0)\big)-\nabla p^\varepsilon(s,x_0)+\nu\Delta u^\varepsilon(s,x_0)\Big){\rm d}s\nonumber\\
	&+\sum_{k=1}^N\int_0^t\psi_\varepsilon\ast\varsigma_k(\tilde u^\varepsilon(s))(x_0){\rm dB}_s^k,\quad \quad \forall t\in[0,+\infty),\forall x_0\in\mathbb T^3.\nonumber
\end{align}
Note that it's not hard to see that the convergence of the martingale term is also in $\mathcal M_{c,2}$. And it's clear that all process under the integrals are progressively measurable for every fixed $x_0\in\mathbb T^3$. For notational simplicity, we omit the notation of space variables when it is fixed arbitrarily or the functions are integrated over $\mathbb T^3$. Now, applying ``$\nabla\times$'' on both sides, we have
\begin{align}
	\omega^\varepsilon(t)=\omega_0^\varepsilon
	&+\int_0^t\Big(-\psi_\varepsilon\ast u^\varepsilon(s)\cdot\nabla \omega^\varepsilon(s)+\nu\Delta\omega^\varepsilon(s)-\epsilon_{\cdot,j,k}\psi_\varepsilon\ast\partial_ju^\varepsilon_l(s)\partial_lu^\varepsilon_k(s)\Big){\rm d}s+\sum_{k=1}^N\int_0^t\varrho_k^\varepsilon(s){\rm dB}_s^k.\nonumber
\end{align}
Here, $\omega_0^\varepsilon=\psi_\varepsilon\ast\omega_0^\varepsilon$ , $\varrho_k^\varepsilon(s):=\psi_\varepsilon\ast\big[\nabla\times\varsigma_k(\tilde u^\varepsilon(s))\big]$ , and $\epsilon_{i,j,k}$ is the usual antisymmetric tensor. And we have used Einstein's contraction of summation. For arbitrarily fixed $\delta\in (0,1/2]$, we define the function of one variable
\begin{align}\label{def. function h}
	h(r):=r^{1/2}-\frac{r^{\frac{1-\delta}{2}}}{2(1-\delta)}\quad (r>0),
\end{align}
and let $q(y):=h\big(1+|y|^2\big)$, $y\in\mathbb R^3$. Then we take the change of dependent variables:
\begin{align}
	\mathbf w^\varepsilon(t,x):=q\big(\omega^\varepsilon(t,x)\big).\nonumber
\end{align}
By It\^o's formula, we have
\begin{align}
	\mathbf w^\varepsilon(t)=\mathbf w^\varepsilon(0)
	&+\int_0^t\Big(-\psi_\varepsilon\ast u^\varepsilon(s)\cdot\nabla \omega_i^\varepsilon(s)+\nu\Delta\omega_i^\varepsilon(s)-\epsilon_{i,j,k}\psi_\varepsilon\ast\partial_ju^\varepsilon_l(s)\partial_lu^\varepsilon_k(s)\Big)\partial_i q\big(\omega^\varepsilon(s)\big){\rm d}s\nonumber\\
	&+\sum_{k=1}^N\int_0^t\varrho_k^\varepsilon(s)\cdot\nabla q\big(\omega^\varepsilon(s)\big){\rm dB}^k_s+\frac{1}{2}\sum_{k=1}^N\int_0^t\varrho_k^\varepsilon(s)^T\nabla^2 q\big(\omega^\varepsilon(s)\big)\varrho_k^\varepsilon(s){\rm d}s.\nonumber
\end{align}
Since $\partial_m\mathbf w^\varepsilon=\partial_m\omega_l^\varepsilon\partial_l q\big(\omega^\varepsilon\big)$ and $\partial^2_m\mathbf w^\varepsilon=\partial^2_m\omega_l^\varepsilon\partial_l q\big(\omega^\varepsilon\big)+\partial_m\omega_i^\varepsilon(s)\partial_m\omega_j^\varepsilon(s)\partial^2_{i,j}q\big(\omega^\varepsilon(s)\big)$	, we have
\begin{align}
	\mathbf w^\varepsilon(t)=\mathbf w^\varepsilon(0)
	&+\int_0^t\Big(-\psi_\varepsilon\ast u^\varepsilon(s)\cdot\nabla \mathbf w^\varepsilon(s)+\nu\Delta\mathbf w^\varepsilon(s)\Big){\rm d}s-\nu\int_0^t\partial_l\omega_i^\varepsilon(s)\partial_l\omega_j^\varepsilon(s)\partial^2_{i,j}q\big(\omega^\varepsilon(s)\big){\rm d}s\nonumber\\
	&+\int_0^t\epsilon_{i,j,k}\psi_\varepsilon\ast\partial_ju^\varepsilon_l(s)\partial_lu^\varepsilon_k(s)\partial_i q\big(\omega^\varepsilon(s)\big){\rm d}s\nonumber\\
	&+\sum_{k=1}^N\int_0^t\varrho_k^\varepsilon(s)\cdot\nabla q\big(\omega^\varepsilon(s)\big){\rm dB}^k_s+\frac{1}{2}\sum_{k=1}^N\int_0^t\varrho_k^\varepsilon(s)^T\nabla^2 q\big(\omega^\varepsilon(s)\big)\varrho_k^\varepsilon(s){\rm d}s.\nonumber
\end{align}
Integrating over $\mathbb T^3$ on both sides, by integration by part, we have
\begin{align}\label{energy equality for bf{w}^varepsilon}
	&\ \ \ \ \int_{\mathbb T^3}\mathbf w^\varepsilon(t){\rm d}x+\nu\int_0^t\int_{\mathbb T^3}\partial_l\omega_i^\varepsilon(s)\partial_l\omega_j^\varepsilon(s)\partial^2_{i,j}q\big(\omega^\varepsilon(s)\big){\rm d}x{\rm d}s\nonumber\\
	&=\int_{\mathbb T^3}\mathbf w^\varepsilon(0){\rm d}x+\int_0^t\int_{\mathbb T^3}\epsilon_{i,j,k}\psi_\varepsilon\ast\partial_ju^\varepsilon_l(s)\partial_lu^\varepsilon_k(s)\partial_i q\big(\omega^\varepsilon(s)\big){\rm d}x{\rm d}s\nonumber\\
	&\ \ \ \ \ +\sum_{k=1}^N\int_0^t\int_{\mathbb T^3}\varrho_k^\varepsilon(s)\cdot\nabla q\big(\omega^\varepsilon(s)\big){\rm d}x{\rm dB}^k_s+\frac{1}{2}\sum_{k=1}^N\int_0^t\int_{\mathbb T^3}\varrho_k^\varepsilon(s)^T\nabla^2 q\big(\omega^\varepsilon(s)\big)\varrho_k^\varepsilon(s){\rm d}x{\rm d}s.
\end{align}
Now, we estimate each term in (\ref{energy equality for bf{w}^varepsilon}). By the obvious inequality
\begin{align}
	\frac{1-2\delta}{2(1-\delta)}\big(1+|y|^2\big)^{1/2}\leq q(y)\leq\big(1+|y|^2\big)^{1/2},\nonumber
\end{align}
we have
\begin{align}
	\int_{\mathbb T^3}\mathbf w^\varepsilon(t){\rm d}x
	&\geq\frac{1-2\delta}{2(1-\delta)}\int_{\mathbb T^3}\Big(1+\big|\omega^\varepsilon(t)\big|^2\Big)^{1/2}{\rm d}x,\label{estimate bf{w}(t)}\\
	\int_{\mathbb T^3}\mathbf w^\varepsilon(0){\rm d}x
	&\leq\int_{\mathbb T^3}\Big(1+\big|\omega_0^\varepsilon\big|^2\Big)^{1/2}{\rm d}x\leq1+\big\|\omega_0^\varepsilon\big\|_{L^1}\leq1+\big\|\omega_0\big\|_{var}.\label{estimate bf{w}(0)}
\end{align}
Next, for any fixed $y,\eta\in\mathbb R^3$, we decompose $\eta=\eta_{\perp}+\eta_{\parallel}$ with $\eta_{\perp}$ perpendicular to $y$ and $\eta_{\parallel}$ parallel to $y$. Then we have
\begin{align}\label{decomposition of quadratic term of q}
	\eta^{T}\nabla^2q(y)\eta
	&=2|\eta|^2h'\big(1+|y|^2\big)+4(\eta\cdot y)^2h''\big(1+|y|^2\big)\nonumber\\
	&=2|\eta_{\perp}|^2h'\big(1+|y|^2\big)+|\eta_{\parallel}|^2\Big(2h'\big(1+|y|^2\big)+4|y|^2h''\big(1+|y|^2\big)\Big).
\end{align}
Together by
\begin{align}
	h'(\alpha)=\frac{1}{2\alpha^{1/2}}-\frac{1}{4\alpha^{\frac{1+\delta}{2}}}>\frac{1}{4\alpha^{\frac{1+\delta}{2}}}\ \ \ (\alpha>1)\nonumber
\end{align}
and
\begin{align}
	2h'(\alpha)+4(\alpha-1)h''(\alpha)=\frac{1}{\alpha^{3/2}}-\frac{1+\delta}{2}\cdot\frac{1}{\alpha^{\frac{3+\delta}{2}}}+\frac{\delta}{2\alpha^{\frac{1+\delta}{2}}}>\frac{\delta}{2\alpha^{\frac{1+\delta}{2}}},\ \ \ (\alpha>1)\nonumber
\end{align}
we see that
\begin{align}
	\eta^{T}\nabla^2q(y)\eta
	&>\frac{1}{2}\big(1+|y|^2\big)^{-\frac{1+\delta}{2}}|\eta_{\perp}|^2+\frac{\delta}{2}\big(1+|y|^2\big)^{-\frac{1+\delta}{2}}|\eta_{\parallel}|^2\nonumber\\
	&>\frac{\delta}{2}\big(1+|y|^2\big)^{-\frac{1+\delta}{2}}|\eta|^2.\nonumber
\end{align}
Hence, we have that
\begin{align}\label{estimate of enstrophy term of bf{w}}
	\nu\int_0^t\int_{\mathbb T^3}\partial_l\omega_i^\varepsilon(s)\partial_l\omega_j^\varepsilon(s)\partial^2_{i,j}q\big(\omega^\varepsilon(s)\big){\rm d}x{\rm d}s>\frac{\nu\delta}{2}\int_0^t\int_{\mathbb T^3}\Big(1+\big|\omega^\varepsilon(s)\big|^2\Big)^{-\frac{1+\delta}{2}}\big|\nabla\omega^\varepsilon(s)\big|^2{\rm d}x{\rm d}s.
\end{align}
Next, by $|\nabla q|\leq1$, H\"older's inequality and condition (\ref{con. control on vorticity}), we have
\begin{align}\label{estimate of antisymmetric term}
	\left|\int_0^t\int_{\mathbb T^3}\epsilon_{i,j,k}\psi_\varepsilon\ast\partial_ju^\varepsilon_l(s)\partial_lu^\varepsilon_k(s)\partial_i q\big(\omega^\varepsilon(s)\big){\rm d}x{\rm d}s\right|\leq6\int_0^t\big\|\nabla u^\varepsilon(s)\big\|_{L^2}^2{\rm d}s
\end{align}
and
\begin{align}
	\mathbb E\sup_{t\in[0,T]}\left|\sum_{k=1}^N\int_0^t\int_{\mathbb T^3}\varrho_k^\varepsilon(s)\cdot\nabla q\big(\omega^\varepsilon(s)\big){\rm d}x{\rm dB}^k_s\right|
	&\leq\mathbb E\left(\sum_{k=1}^N\int_0^T\left|\int_{\mathbb T^3}\varrho_k^\varepsilon(s)\cdot\nabla q\big(\omega^\varepsilon(s)\big){\rm d}x\right|^2{\rm d}s\right)^{1/2}\nonumber\\
	&\lesssim\mathbb E\left(\int_0^T\sum_{k=1}^N\big\|\nabla\times\varsigma_k(u^\varepsilon(s))\big\|_{L^2}^2{\rm d}s\right)^{1/2}\nonumber\\
	&\lesssim\mathbb E\left[\int_0^T\left(1+\big\|\nabla u^\varepsilon(s)\big\|_{L^2}^2\right){\rm d}s\right]^{1/2}.
\end{align}
Finally, by (\ref{decomposition of quadratic term of q}) and that $\displaystyle|h'(\alpha)|<\frac{1}{2\alpha^{1/2}}$, $\displaystyle|h''(\alpha)|<\frac{1}{4\alpha^{3/2}}$ for $\alpha>1$, we have
\begin{align}
	\big|\eta^{T}\nabla^2q(y)\eta\big|
	&\leq2|\eta|^2\big|h'\big(1+|y|^2\big)\big|+4|\eta|^2|y|^2\big|h''\big(1+|y|^2\big)\big|\nonumber\\
	&\leq2\big(1+|y|^2\big)^{-1/2}|\eta|^2.\nonumber
\end{align}
This and (\ref{con. control on vorticity}) give that
\begin{align}\label{estimate of cross variation term}
	\left|\frac{1}{2}\sum_{k=1}^N\int_0^t\int_{\mathbb T^3}\varrho_k^\varepsilon(s)^T\nabla^2 q\big(\omega^\varepsilon(s)\big)\varrho_k^\varepsilon(s){\rm d}x{\rm d}s\right|
	&\leq\sum_{k=1}^N\int_0^t\int_{\mathbb T^3}\Big(1+\big|\omega^\varepsilon(s)\big|^2\Big)^{-\frac{1}{2}}\big|\varrho_k^\varepsilon(s)\big|^2{\rm d}x{\rm d}s\nonumber\\
	&\leq\int_0^t\sum_{k=1}^N\big\|\nabla\times\varsigma_k(u^\varepsilon(s))\big\|_{L^2}^2{\rm d}s\nonumber\\
	&\lesssim\int_0^t\left(1+\big\|\nabla u^\varepsilon(s)\big\|_{L^2}^2\right){\rm d}s.
\end{align}
Hence, by (\ref{energy equality for bf{w}^varepsilon}), (\ref{estimate bf{w}(t)}), (\ref{estimate bf{w}(0)}), (\ref{estimate of enstrophy term of bf{w}}), (\ref{estimate of antisymmetric term}) and (\ref{estimate of cross variation term}), we have for any $T>0$ that
\begin{align}\label{main estimate}
	&\ \ \ \ \ \frac{1-2\delta}{2(1-\delta)}\mathbb E\sup_{t\in[0,T]}\int_{\mathbb T^3}\Big(1+\big|\omega^\varepsilon(t)\big|^2\Big)^{1/2}{\rm d}x+\frac{\nu\delta}{2}\mathbb E\int_0^T\int_{\mathbb T^3}\Big(1+\big|\omega^\varepsilon(s)\big|^2\Big)^{-\frac{1+\delta}{2}}\big|\nabla\omega^\varepsilon(s)\big|^2{\rm d}x{\rm d}s\nonumber\\
	&\leq1+\big\|\omega_0\big\|_{var}+C_5~\mathbb E\int_0^T\left(1+\big\|\nabla u^\varepsilon(s)\big\|_{L^2}^2\right){\rm d}s.
\end{align}
Here $C_5>0$ is some constant only possibly depending on $T$. 
Then by (\ref{standard uniform bounds 1}), the right hand-side of (\ref{main estimate}) is bounded uniformly in $\varepsilon>0$. This already shows that
\begin{align}
	\mathbb E\sup_{0\leq t\leq T}\|\omega^\varepsilon(t)\|_{L^1}\leq C\big(\|u_0\|_{L^2},\|\omega_0\|_{var},T\big)
\end{align}
for some constant $C>0$.
The uniform bound for the other term in (\ref{vorticity bounds for Leray regularised solution}) follows by (\ref{main estimate}) and direct application of H\"older's inequlity.
\end{proof}
\

\appendix
\renewcommand{\appendixname}{Appendix~\Alph{section}}
\renewcommand{\theequation}{A.\arabic{equation}}
\section{Intuitive computation for local energy equality}\label{Appen. Intuitive computation}
Here, we give an intuitive computation for local energy equality of (\ref{sto. NS}). Assuming that $\mathbb P$-$a.s.$, (\ref{sto. NS}) can be written as 
\begin{align}
	u(t,x)=u_0(x)+\int_0^t\Big(-{\rm div}\big(u(s,x)\otimes u(s,x)\big)-\nabla p(s,x)+\nu\Delta u(s,x)\Big){\rm d}s+\sum_{k}&\int_0^t\varsigma_k(
	u(s))(x){\rm dB}^k_s
	,\nonumber\\
	&\quad\quad\ \,\forall t\in[0,+\infty)\nonumber
\end{align}
for every $x\in\mathbb T^3$. Now, fix an $x\in\mathbb T^3$ arbitrarily. Applying It\^o's formula to the process $t\mapsto \frac{1}{2}|u(t,x)|^2\varphi(t,x)$ with $\varphi\in C^\infty$, we have
\begin{align}
	&\ \ \ \ \ \frac{1}{2}|u(t,x)|^2\varphi(t,x)-\frac{1}{2}|u(0,x)|^2\varphi(0,x)\textcolor{magenta}{}\nonumber\\
	&=\frac{1}{2}\int_0^t|u(s,x)|^2\partial_s\varphi(s,x){\rm d}s+\int_0^t\varphi(s,x)u(s,x)\cdot\Big(-u(s,x)\cdot\nabla u(s,x)+\nu\Delta u(s,x)-\nabla p(s,x)\Big){\rm d}s\nonumber\\
	&\ \ \ \ \ +\sum_k\int_0^t\varsigma_k(
	u(s))(x)\cdot u(s,x)\varphi(s,x){\rm dB}^k_s+\frac{1}{2}\sum_k\int_0^t|\varsigma_k(u(s))(x)|^2\varphi(s,x){\rm d}s.\nonumber\\
\end{align}
Integrating on $\mathbb T^3$ on both sides, we have
\begin{align}
	&\int_{\mathbb T^3}|u(t)|^2\varphi(t){\rm d}x-\int_{\mathbb T^3}|u_0|^2\varphi(0){\rm d}x\nonumber\\
	=&-\underbrace{2\int_0^t\int_{\mathbb T^3}\varphi(s)u(s)\cdot\big(u(s)\cdot\nabla u(s)\big){\rm d}x{\rm d}s}_{I_1(t)}+\int_0^t\int_{\mathbb T^3}|u(s)|^2\partial_s\varphi(s){\rm d}x{\rm d}s\nonumber\\
	&+\underbrace{2\nu\int_0^t\int_{\mathbb T^3}\varphi(s)u(s)\cdot\Delta u(s){\rm d}x{\rm d}s}_{I_2(t)}
	-\underbrace{2\int_0^t\int_{\mathbb T^3}\varphi(s)u(s)\cdot\nabla p(s){\rm d}x{\rm d}s}_{I_3(t)}\nonumber\\
	&+2\sum_k\int_0^t\int_{\mathbb T^3}\varsigma_k(
	u(s))\cdot u(s)\varphi(s){\rm d}x{\rm dB}^k_s+\sum_k\int_0^t\int_{\mathbb T^3}|\varsigma_k(u(s))|^2\varphi(s){\rm d}x{\rm d}s.
\end{align}
Then by integration by part,
\begin{align}
	I_1(t)
	&=\int_0^t\int_{\mathbb T^3}\varphi(s)u(s)\cdot\nabla|u(s)|^2{\rm d}x{\rm d}s=-\int_0^t\int_{\mathbb T^3}|u(s)|^2u(s)\cdot\nabla\varphi(s){\rm d}x{\rm d}s,\nonumber\\
	I_2(t)
	&=-2\nu\int_0^t\int_{\mathbb T^3}\varphi(s)|\nabla u(s)|^2{\rm d}x{\rm d}s-2\nu\int_0^t\int_{\mathbb T^3}\partial_i\varphi(s)u_j(s)\partial_iu_j(s){\rm d}x{\rm d}s\nonumber\\
	&=-2\nu\int_0^t\int_{\mathbb T^3}\varphi(s)|\nabla u(s)|^2{\rm d}x{\rm d}s+\int_0^t\int_{\mathbb T^3}|u(s)|^2\nu\Delta\varphi(s){\rm d}x{\rm d}s,\nonumber\\
	I_3(t)
	&=-\int_0^t\int_{\mathbb T^3}2p(s)u(s)\cdot\nabla\varphi(s){\rm d}x{\rm d}s.\nonumber
\end{align}
Hence, we have
\begin{align}
	&\int_{\mathbb T^3}|u(t)|^2\varphi(t){\rm d}x+2\nu\int_0^t\int_{\mathbb T^3}|\nabla u(s)|^2\varphi(s){\rm d}x{\rm d}s\nonumber\\
	=&\int_{\mathbb T^3}|u_0|^2\varphi(0){\rm d}x+\int_0^t\int_{\mathbb T^3}|u(s)|^2\Big(\partial_s\varphi(s)+\nu\Delta\varphi(s)\Big){\rm d}x{\rm d}s+\int_0^t\int_{\mathbb T^3}\Big(|u(s)|^2+2p(s)\Big)u(s)\cdot\nabla\varphi(s){\rm d}x{\rm d}s\nonumber\\
	&+2\sum_k\int_0^t\int_{\mathbb T^3}\varsigma_k(
	u(s))\cdot u(s)\varphi(s){\rm d}x{\rm dB}^k_s+\sum_k\int_0^t\int_{\mathbb T^3}|\varsigma_k(u(s))|^2\varphi(s){\rm d}x{\rm d}s.\nonumber
\end{align}
This is a local energy balance for (\ref{sto. NS}).
\\
\bibliographystyle{alpha}
\bibliography{Paper2ref}

\end{document}